\documentclass[12pt,reqno,a4paper]{amsart}

\usepackage{amsfonts,amsmath,amssymb,amsthm}
\usepackage{mathrsfs}
\usepackage{mathtools}
\usepackage{xcolor}
\usepackage{tikz}
\usetikzlibrary{arrows.meta,positioning,calc}
\usepackage{float}
\usepackage{placeins}

\usepackage{a4wide}
\allowdisplaybreaks[4]
\usepackage[T1,OT1]{fontenc}
\usepackage[english]{babel}
\usepackage[backref=page]{hyperref}
\hypersetup{
  colorlinks=true,
  linkcolor=blue,
  citecolor=magenta,
  urlcolor=cyan,
  pdftitle={Modular Schur Functions and Macdonald Characters},
  pdfauthor={Naihuan Jing and Ning Liu},
  pdfkeywords={modular Schur functions, Macdonald characters, truncated complete symmetric functions, Schur coefficients, Carter condition, Heisenberg operators, vertex operators}
}

\numberwithin{equation}{section}
\numberwithin{figure}{section}

\theoremstyle{plain}
\newtheorem{theorem}{Theorem}[section]
\newtheorem{mainthm}{Theorem}

\newtheorem{proposition}[theorem]{Proposition}
\newtheorem{lemma}[theorem]{Lemma}
\newtheorem{corollary}[theorem]{Corollary}

\newtheorem{remark}[theorem]{Remark}
\theoremstyle{definition}
\newtheorem{example}[theorem]{Example}

\newcommand{\Sym}{\mathfrak S}
\newcommand{\cP}{\mathcal P}
\newcommand{\eps}{\varepsilon}
\newcommand{\inner}[2]{\left\langle #1,#2\right\rangle}
\newcommand{\emptypar}{\varnothing}
\DeclareMathOperator{\ch}{ch}
\DeclareMathOperator{\Ind}{Ind}

\DeclareMathOperator{\soc}{soc}
\DeclareMathOperator{\Coef}{Coef}

\title[Modular Schur Functions and Macdonald Characters]
{Modular Schur Functions and Macdonald Characters}
\author{Naihuan Jing}
\address{Department of Mathematics, North Carolina State University, Raleigh, NC 27695, USA}
\email{jing@ncsu.edu}
\author{Ning Liu}
\address{School of Mathematics, Hefei University of Technology, Hefei, Anhui 230009, China}
\email{mathliu123@outlook.com}



\begin{document}

\begin{abstract}
Modular Schur functions arose from the problem of determining the
characters of simple polynomial representations of the general linear
group in positive characteristic.  We prove three results about these
functions.  (i) In finitely many variables, we obtain two bases for a
space of symmetric polynomials with bounded degree in each variable,
one consisting of modular Schur polynomials and the other of products
of truncated complete symmetric polynomials.
(ii) We determine the signs of the nonzero Schur coefficients of modular
Schur functions, show that their absolute values are at most the factorial
of the block weight, and characterize equality.
(iii) Over an infinite field of prime characteristic, a modular Schur
polynomial is a simple polynomial character exactly
when its indexing partition satisfies the Carter condition, provided
the number of variables is at least the degree. Thus we prove three conjectures formulated by Walker (1994).
We also use the Heisenberg operators from the Schur coefficient
calculation to prove positivity for Macdonald characters.
After centering the Macdonald parameters at one and negating the
auxiliary variables, the normalized power-sum matrix entries of the
Macdonald creation operator have nonnegative integer coefficients,
which yields nonnegative rational coefficients for normalized
Macdonald characters.  The result for the creation operator confirms
a conjecture of Ben Dali and D'Adderio (2024).
\end{abstract}

\maketitle

\setcounter{tocdepth}{1}
\tableofcontents

\section{Introduction}
\label{sec:introduction}

\subsection{Background and motivation}

Modular Schur functions arise from truncated symmetric powers in
positive characteristic \cite{DotyWalker,Walker}.  To see the
truncation on monomials, fix an integer $r\geq2$.  In the ring of
symmetric functions over $\mathbb Z$, let $h_m$ be the complete
symmetric function of degree $m\geq0$.  Retaining only the monomials
in which each variable has exponent less than $r$ gives the truncated
complete function $\mathfrak h_m$, with $\mathfrak h_0=h_0=1$.
Since $h_1,h_2,\ldots$ freely generate the ring, $h_m\mapsto\mathfrak h_m$
defines an endomorphism $\Theta_r$.  For a partition $\lambda$, the
image of the ordinary Schur function $s_\lambda$ is the modular Schur
function $\mathfrak s_\lambda=\Theta_r(s_\lambda)$.  We study three
questions about these functions.  The operators that describe their
Schur coefficients also lead to a positivity result for Macdonald
characters.

The first question concerns bases in $n\geq1$ variables.  Over an
algebraically closed field of prime characteristic $r$, certain
simple $GL_n$-modules embed in tensor products of truncated symmetric
powers \cite[Section~1]{DotyWalkerTruncated}.  Doty conjectured that the
matrix of composition multiplicities for the partitions in an
$n$ by $n(r-1)$ rectangle is nonsingular.  Doty and Walker showed that
this is equivalent to linear independence of the products of truncated
complete symmetric polynomials \cite[Conjecture~2.6 and Proposition~3.1]{DotyWalker}.
Walker stated the independence conjecture for every integer $r\geq2$
and gave its equivalent formulation for modular Schur polynomials
\cite[Conjecture~1.4 and Proposition~1.5(i)]{Walker}.
Walker credits the two-variable case to Doty and Walker and proves
the three-variable case
\cite[p.~576 and Proposition~2.17]{Walker}.
For prime $r$, Donkin and Geranios showed that every column of the
multiplicity matrix is nonzero \cite[Corollary~6.8, p.~432]{DonkinGeranios}.
Linear independence of its columns would allow the simple characters
in the rectangle to be recovered from the tensor-product characters.

The second question asks for the signs and sizes of the Schur
coefficients of $\mathfrak s_\lambda$.  For a single row $(m)$,
$\mathfrak s_{(m)}=\mathfrak h_m$, and its Schur coefficients belong
to $\{0,1,-1\}$ \cite{FuMei,Grinberg}.  Cheng, Chou, Eu, Fu, and Yao
gave a ribbon interpretation for this row of the transition matrix
and asked for a combinatorial expansion for the other rows
\cite[Theorem~1.3 and Section~4]{ChengEtAl}.  We answered this
question by expressing every entry of the matrix as a signed sum over
pairs of ribbon tableaux \cite[Theorem~8.8]{JingLiuSkewMN}.
For general $\lambda$,
Walker conjectured a sign rule and a factorial bound on the
coefficients \cite[unnumbered conjecture following Corollary~2.10, p.~586]{Walker}.

The third question asks when, for $r=p$ prime, a modular Schur
function is a simple character in characteristic $p$.  The Carter condition
requires the $p$-adic valuations of the hook lengths to be constant
down each column.  For a $p$-regular partition (one with no positive
part repeated $p$ times or more), this condition characterizes
irreducibility of the Specht module \cite{JamesCarter,JamesMurphy}.
Walker proved that the Carter condition suffices for simplicity and,
for $p$-regular partitions, is necessary in sufficiently large rank.
He conjectured necessity for all partitions and observed that when the
conjugate of $\lambda$ is $p$-regular, simplicity forces
$\mathfrak s_\lambda=s_\lambda$
\cite[Theorem~4.4 and the discussion of Conjectures~4.6--4.7]{Walker}.
For odd $p$, we proved that this equality forces $\lambda$ to be a
$p$-core \cite[Theorem~8.11]{JingLiuSkewMN}.  All hook lengths of a
$p$-core have $p$-adic valuation zero, so in odd characteristic a
counterexample could occur only if neither $\lambda$ nor its conjugate
is $p$-regular.

Our fourth question concerns Macdonald characters.  In the
one-parameter setting, Ben Dali and {\fontencoding{T1}\selectfont Do{\l}\k{e}ga} proved
Lassalle's positivity conjecture for Jack characters in Stanley's
multirectangular coordinates \cite{LassalleJack,BenDaliDolega}.
Ben Dali and D'Adderio introduced Macdonald characters as a
two-parameter generalization.  They express these characters using a
product of creation operators and show that character positivity
follows from their conjecture for a single creation operator.  With
Macdonald parameters $q,t$ and auxiliary variable $v$, the conjecture
asserts that its normalized matrix entries in the power-sum basis are
polynomials in $q-1$, $t-1$, and $-v$ with nonnegative integer
coefficients \cite[Definition~1.4, Conjecture~1 and Proposition~4.7]{BD26}.

Our aim is to answer these four questions in the ring of symmetric
functions and its finite-variable specializations.  We construct bases
in finitely many variables from the truncated complete functions.
Heisenberg operators give a formula for the Schur coefficients of
modular Schur functions, and we use these coefficients to determine
which modular Schur functions are simple characters.  The same
Heisenberg operators occur in the first-named author's Hall--Littlewood vertex operators
\cite{JingVertex}, whose modes we use to prove the positivity conjecture
for the Macdonald creation operator.

\subsection{Main results}

We now give precise statements of our four main results.

We begin with two bases for a space of symmetric polynomials with a
degree bound in each variable.
For a partition $\lambda$, write $|\lambda|$ for its size,
$\ell(\lambda)$ for its length, and $\mathfrak h_\lambda$ for
$\prod_{i=1}^{\ell(\lambda)}\mathfrak h_{\lambda_i}$, with empty
product $1$.  In $n$ variables, a product of at most $n$ truncated
complete functions has degree at most $n(r-1)$ in each variable.
Let $V_{n,n(r-1)}$ be the space of symmetric polynomials over
$\mathbb Q$ satisfying this bound, and put
\[
 \mathcal R_{n,r}=\{\lambda:\ell(\lambda)\leq n,\ \lambda_1\leq n(r-1)\},
\]
with first part zero for the empty partition.

\begin{mainthm}[{\hyperref[thm:finite-variable-bases]{Theorem~\ref*{thm:finite-variable-bases}}}]
\label{thm:intro-finite-bases}
For every $n\geq1$ and $r\geq2$, each of the two families
\[
 \{\mathfrak h_\lambda(x_1,\ldots,x_n):\lambda\in\mathcal R_{n,r}\},
 \qquad
 \{\mathfrak s_\lambda(x_1,\ldots,x_n):\lambda\in\mathcal R_{n,r}\}
\]
is a $\mathbb Q$-basis of $V_{n,n(r-1)}$.
\end{mainthm}

This proves Walker's Conjecture~1.4.  For $n\geq2$, the proof of
Theorem~\ref{thm:intro-finite-bases} begins by expressing the truncated
complete symmetric polynomials as coefficients of a product over the
nontrivial $r$th roots of unity.
Relations among these coefficients reduce the question of linear independence to the
invertibility of a multiplication map on sections over projective space.
Lemma~\ref{lem:linear-presentation-multiplication} shows that this
map is an isomorphism.

Our second result gives a sign rule and a sharp bound for the
coefficients of modular Schur functions in the Schur basis.
Recall that an $r$-ribbon is an
edge-connected skew diagram of $r$ boxes containing no $2\times2$
square.  Successively removing such ribbons from $\lambda$ leaves its
$r$-core.  The number removed is its weight $w$, and its $r$-sign
$\eps_r(\lambda)$ is the product of $(-1)^{h}$ over the removals,
where $h$ is the number of rows of a ribbon minus one.
By homogeneity and Proposition~\ref{prop:block-preservation}, only
partitions of the same size and $r$-core can give a nonzero Schur
coefficient.  For a fixed core, the $r$-quotient indexes the partitions
of weight $w$ by $r$-multipartitions, that is, ordered $r$-tuples
of partitions of total size $w$.  We use the runner order
of Subsection~\ref{subsec:abacus-core-quotient} and conjugate each
quotient component.  If $\lambda,\mu$ have the same core and
weight, and $\boldsymbol\alpha,\boldsymbol\beta$ are the resulting
tuples, their normalized Schur coefficient is
\begin{equation}\label{eq:refaudit-intro-coefficients}
 \kappa^{(r,w)}_{\boldsymbol\alpha\boldsymbol\beta}
 =(-1)^w\eps_r(\lambda)\eps_r(\mu)\,
    \inner{\mathfrak s_\lambda}{s_\mu}.
\end{equation}
Here $\inner{\cdot}{\cdot}$ is the Hall inner product for which the
Schur functions are orthonormal.  The coefficient depends only on the
quotients \cite[Corollary~2.10]{Walker}.  For an $r$-tuple
$\boldsymbol\alpha=(\alpha^{(0)},\ldots,\alpha^{(r-1)})$ of total
size $w$, define
\[
 d_{\boldsymbol\alpha}
 =\frac{w!}{\prod_i|\alpha^{(i)}|!}\prod_i f^{\alpha^{(i)}},
\]
where $f^\rho$ counts standard Young tableaux of shape $\rho$ and
$f^\varnothing=1$.

\begin{mainthm}[{\hyperref[thm:universal-refined-bound]{Theorem~\ref*{thm:universal-refined-bound}}
  $+$ \hyperref[eq:factorial-equality-case]{(\ref*{eq:factorial-equality-case})}
  $+$ \hyperref[thm:coefficient-conjecture]{Theorem~\ref*{thm:coefficient-conjecture}}}]
\label{thm:intro-universal-matrices}
For $r\geq2$, $w\geq0$, and any two $r$-multipartitions
$\boldsymbol\alpha,\boldsymbol\beta$ of size $w$, the coefficient
$\kappa^{(r,w)}_{\boldsymbol\alpha\boldsymbol\beta}$ is an integer and
\[
 0\leq\kappa^{(r,w)}_{\boldsymbol\alpha\boldsymbol\beta}
 \leq\min\{d_{\boldsymbol\alpha},d_{\boldsymbol\beta}\}\leq w!.
\]
Equality $\kappa^{(r,w)}_{\boldsymbol\alpha\boldsymbol\beta}=w!$
holds exactly when every component is $\varnothing$ or $(1)$ and
the two sets of nonempty components are disjoint.
\end{mainthm}

To prove Theorem~\ref{thm:intro-universal-matrices}, we pass to quotient
coordinates.  Proposition~\ref{prop:core-independent-matrix} then
identifies the coefficients in \eqref{eq:refaudit-intro-coefficients}
with the Schur coefficients of a product of dual Cauchy kernels.
The Littlewood--Richardson rule gives nonnegativity.  Comparing
with the regular representation of the symmetric group gives the
bound, and the expansion of the kernel determines when equality
holds.  This proves and sharpens Walker's coefficient conjecture.

The third result determines when a modular Schur function is the
character of a simple polynomial representation.  Take $r=p$ prime, let $\mathbb K$
be an infinite field of characteristic $p$, and let $\lambda$ be a
partition of $d$.  We take $n$ to be a positive integer.

\begin{mainthm}[{\hyperref[cor:carter-classification]{Corollary~\ref*{cor:carter-classification}}}]
\label{thm:intro-simple-classification}
For every $n\geq d$, the specialization of $\mathfrak{s}_\lambda$
to $n$ variables is the character of a simple degree-$d$ polynomial
$GL_n(\mathbb K)$-module if and only if $\lambda$ satisfies the
Carter condition.
\end{mainthm}

This proves Walker's Conjecture~4.6.  To prove
Theorem~\ref{thm:intro-simple-classification}, we show that simplicity
forces the Carter condition, since Walker proved the converse.
For odd $p$, Proposition~\ref{prop:restricted-projection} rules out
simplicity when the corresponding Specht module is reducible and gives
a unique candidate for the simple character when it is irreducible.
We reduce the remaining cases to Rouquier blocks, where a Schur
coefficient computed using Proposition~\ref{prop:universal-skew-expansion}
distinguishes the modular Schur function from that candidate.
When $p=2$, truncation sends $s_\lambda$ to the Schur function indexed
by the conjugate partition, so necessity follows directly from the
simplicity criterion for Weyl modules (Proposition~\ref{prop:prime-two}).

Finally, we turn to the Macdonald creation operator and prove positivity
for its normalized coefficients in the power-sum basis.  Work over $\mathbb Q(q,t)$
with $q,t$ algebraically independent, and let $J_\lambda(X;q,t)$ be the integral
form of the Macdonald symmetric function in the alphabet
$X=(x_1,x_2,\ldots)$ \cite[Chapter~VI, (8.3)]{Macdonald}.
For $m\geq1$ write $p_m(X)=\sum_i x_i^m$.  In plethystic notation
$\mathcal P_A$ denotes multiplication by
$\operatorname{Exp}[AX]=\exp(\sum_{m\geq1}p_m[A]p_m(X)/m)$,
where $p_m[A]$ is the $m$th power sum of the alphabet $A$.
For formal parameters $c,u,v$, define
\[
 \boldsymbol\Delta_cJ_\lambda
 =\prod_{(i,j)\in\lambda}(1-cq^{j-1}t^{-(i-1)})J_\lambda,
 \qquad
 \boldsymbol\Gamma(u,v)
 =\boldsymbol\Delta_{1/v}\mathcal P_{uv(1-t)/(1-q)}
     \boldsymbol\Delta_{1/v}^{-1}.
\]
The product runs over row and column indices of boxes in $\lambda$,
and $p_m[uv(1-t)/(1-q)]=(uv)^m(1-t^m)/(1-q^m)$.
Write $\Coef_b$ for extraction of the coefficient of a monomial or
basis element $b$.  Define the algebra automorphism $R$ and the
operators $G_n(v)$ by
\[
 Rp_m=\frac{1-t^m}{1-q^m}p_m\quad(m\geq1),\qquad
 G_n(v)=(-1)^n\Coef_{u^n}\bigl[R^{-1}\boldsymbol\Gamma(u,v)R\bigr]
 \quad(n\geq0).
\]
For a partition $\nu$ with $m_j(\nu)$ parts equal to $j$, put
$p_\nu=\prod_jp_j^{m_j(\nu)}$ and
$z_\nu=\prod_jj^{m_j(\nu)}m_j(\nu)!$.

\begin{mainthm}[{\hyperref[c1:thm:conjecture-one]{Theorem~\ref*{c1:thm:conjecture-one}}}]
\label{thm:intro-c1}
For every $n\geq0$ and all partitions $\mu,\nu$ with
$|\nu|=|\mu|+n$,
\[
 t^{|\mu|}z_\nu\Coef_{p_\nu}\bigl[G_n(v)p_\mu\bigr]
 \in\mathbb N[q-1,t-1,-v].
\]
\end{mainthm}

By \eqref{c1:eq:original-pairing}, this proves Conjecture~1 of
\cite{BD26}.  Our proof of Theorem~\ref{thm:intro-c1} uses the exchange
relation for the Heisenberg operators to obtain nonnegative coefficients for normalized repeated
Hall--Littlewood modes.  These coefficients occur in the expansion
of the Macdonald creation operator and yield nonnegativity.
Their integrality in the Schur basis makes the ordinary Hall pairing
of $t^{|\mu|}G_n(v)p_\mu$ with $p_\nu$ integral.  Applying the result
to a product of creation operators
gives the character positivity in Corollary~\ref{c1:cor:character-positivity}.

\subsection{Organization}

Section~\ref{sec:preliminaries} recalls truncation and finite-variable
specialization.  Section~\ref{sec:finite-variable-bases} uses
Lemma~\ref{lem:linear-presentation-multiplication}, proved in
Appendix~\ref{sec:appendix-lemma}, to show that the two families
$\{\mathfrak h_\lambda\}_{\lambda\in\mathcal R_{n,r}}$ and
$\{\mathfrak s_\lambda\}_{\lambda\in\mathcal R_{n,r}}$ each form a basis
of $V_{n,n(r-1)}$, thereby proving Theorem~\ref{thm:intro-finite-bases}.
Section~\ref{sec:universal-blocks} describes the Heisenberg action in
core--quotient coordinates.  Section~\ref{sec:modular-schur-coefficients}
uses the resulting kernel to prove the sign rule and sharp bounds for
modular Schur coefficients in Theorem~\ref{thm:intro-universal-matrices}.
Section~\ref{sec:conj46} determines which modular Schur functions are
simple characters, proving Theorem~\ref{thm:intro-simple-classification}.
Section~\ref{sec:hl-mode-positivity}
establishes positivity and integrality for the normal-order coefficients
of Hall--Littlewood modes.  Section~\ref{sec:c1-positivity} proves
Theorem~\ref{thm:intro-c1}, which settles Ben Dali and D'Adderio's
Conjecture~1 on the creation operator and yields positivity for
Macdonald characters.

\smallskip
\noindent\begin{minipage}[t]{\linewidth}
The arrows below point from each section to those that use its results.
\par\smallskip
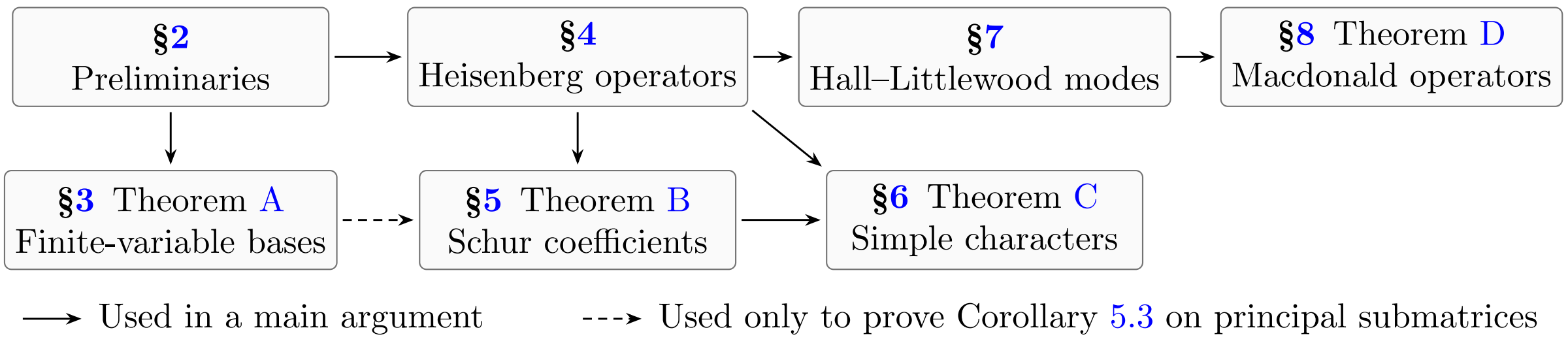
\centering
\begin{tikzpicture}[
  section box/.style={draw=black!55,fill=black!2,rounded corners=2pt,
    minimum width=3.2cm,minimum height=1cm,inner sep=3pt,
    align=center,font=\footnotesize},
  dependency/.style={-{Stealth[length=4.5pt]},line width=.55pt,
    shorten <=1.5pt,shorten >=1.5pt},
  application/.style={dependency,dash pattern=on 2.6pt off 1.8pt},
  key/.style={right,inner xsep=3pt,font=\footnotesize}]
\node[section box] (prelim) at (0,0)
  {\textbf{\S\ref{sec:preliminaries}}\\Preliminaries};
\node[section box] (heis) at (4.12,0)
  {\textbf{\S\ref{sec:universal-blocks}}\\Heisenberg operators};
\node[section box] (modes) at (8.24,0)
  {\textbf{\S\ref{sec:hl-mode-positivity}}\\Hall--Littlewood modes};
\node[section box] (mac) at (12.36,0)
  {\textbf{\S\ref{sec:c1-positivity}}\enspace Theorem~\ref{thm:intro-c1}\\Macdonald operators};
\node[section box] (bases) at (0,-1.65)
  {\textbf{\S\ref{sec:finite-variable-bases}}\enspace Theorem~\ref{thm:intro-finite-bases}\\Finite-variable bases};
\node[section box] (coeff) at (4.12,-1.65)
  {\textbf{\S\ref{sec:modular-schur-coefficients}}\enspace Theorem~\ref{thm:intro-universal-matrices}\\Schur coefficients};
\node[section box] (simple) at (8.24,-1.65)
  {\textbf{\S\ref{sec:conj46}}\enspace Theorem~\ref{thm:intro-simple-classification}\\Simple characters};
\draw[dependency] (prelim) -- (heis);
\draw[dependency] (prelim) -- (bases);
\draw[dependency] (heis) -- (modes);
\draw[dependency] (modes) -- (mac);
\draw[dependency] (heis) -- (coeff);
\draw[dependency] (heis.south east) -- (simple.north west);
\draw[dependency] (coeff) -- (simple);
\draw[application] (bases) -- (coeff);
\draw[dependency] (-1.55,-2.65) -- ++(.7,0)
  node[key]{Used in a main argument};
\draw[application] (4.12,-2.65) -- ++(.7,0)
  node[key]{Used only to prove Corollary~\ref{cor:rectangular-principal-matrices} on principal submatrices};
\end{tikzpicture}
\end{minipage}

\section{Preliminaries}
\label{sec:preliminaries}

We fix an integer $r\geq2$ and collect the notation and
conventions for partitions and symmetric functions.
We also record the finite-variable specialization and the spaces used
in Theorem~\ref{thm:intro-finite-bases}.

\subsection{Partitions and symmetric functions}
\label{subsec:ordinary-symmetric-functions}

Throughout this subsection we use the conventions for symmetric
functions in Macdonald's book \cite[Chapter~I]{Macdonald}.
A \emph{partition} is a weakly decreasing sequence
$\lambda=(\lambda_1,\lambda_2,\ldots)$ of nonnegative integers with
only finitely many nonzero terms.  Its size and length are
\[
 |\lambda|=\sum_i\lambda_i,
 \qquad
 \ell(\lambda)=\#\{i:\lambda_i>0\}.
\]
We write $\lambda\vdash d$ when $|\lambda|=d$, and denote the empty
partition by $\varnothing$.
Let $\mathscr P$ be the set of all partitions, and, for $d\geq0$, set
\[
 \mathscr P_d=\{\lambda\in\mathscr P:|\lambda|=d\}.
\]
For partitions of the same size, $\mu<_{\mathrm{lex}}\lambda$ means
that at the first index where their parts differ, the part of $\mu$
is smaller.  We write $\alpha+\beta$ for componentwise addition,
appending zero parts as necessary.
We draw partitions in English notation, so the diagram of $\lambda$
has $\lambda_a$ boxes in row $a$, counted from the top.  The conjugate
partition $\lambda'$ is obtained by transposing this diagram; in
particular, $\lambda'_b=\#\{a:\lambda_a\geq b\}$.  We write
$\mu\subseteq\lambda$ when the diagram of $\mu$ is contained in that
of $\lambda$.  For $\mu\subseteq\lambda$, the \emph{skew Young diagram}
$\lambda/\mu$ consists of the boxes of $\lambda$ that do not belong to
$\mu$.  A skew Young diagram is a \emph{ribbon} if it is connected
through common edges and contains no $2\times2$ square of boxes.  Its
length is its number of boxes, and its height is one less than the number
of rows it occupies.  Thus a removable $r$-ribbon of $\lambda$ is a
ribbon $\lambda/\mu$ of length $r$ for some partition
$\mu\subseteq\lambda$.

Let $X=\{x_1,x_2,\ldots\}$ be an alphabet, and let $\Lambda$ be the
graded ring of symmetric functions over $\mathbb Z$.  We write
$\Lambda_{\mathbb Q}=\mathbb Q\otimes_{\mathbb Z}\Lambda$.  For
$d\geq0$, let $\Lambda^d$ and $\Lambda_{\mathbb Q}^d$ denote the
homogeneous components of degree $d$ in these two rings, respectively.
For
$m\geq1$, the complete and power-sum symmetric functions are
\[
 h_m(X)=\sum_{i_1\leq\cdots\leq i_m}x_{i_1}\cdots x_{i_m},
 \qquad
 p_m(X)=\sum_{j\geq1}x_j^m,
\]
with $h_0=1$ and $h_m=0$ for $m<0$.  Their ordinary generating
series is
\begin{equation}\label{eq:ordinary-complete-generating}
 H(X;t):=\sum_{m\geq0}h_m(X)t^m
 =\prod_{j\geq1}\frac{1}{1-x_jt}
 =\exp\!\left(\sum_{m\geq1}\frac{p_m(X)t^m}{m}\right).
\end{equation}
For a partition $\eta=(1^{m_1}2^{m_2}\cdots)$, where $m_j$ is the
multiplicity of the part $j$, put
\[
 p_\eta=\prod_{j\geq1}p_j^{m_j},\qquad
 z_\eta=\prod_{j\geq1}j^{m_j}m_j!.
\]
Empty products have value $1$.
The monomial symmetric function $m_\lambda$ is the sum of the
distinct monomials whose positive exponents, rearranged in decreasing
order, form $\lambda$.

For $\mu\subseteq\lambda$, the skew Jacobi--Trudi identity gives
\begin{equation}\label{eq:refaudit-skew-jacobi-trudi}
 s_{\lambda/\mu}
 =\det\bigl(h_{\lambda_i-\mu_j-i+j}\bigr)_{1\leq i,j\leq\ell}
 \qquad(\ell\geq\ell(\lambda)),
\end{equation}
where zero parts are appended to $\mu$ as needed.  We set
$s_{\lambda/\mu}=0$ if $\mu\not\subseteq\lambda$.
The ordinary Schur
function is $s_\lambda=s_{\lambda/\varnothing}$.  We equip $\Lambda$
with the Hall inner product $\inner{\cdot}{\cdot}$, characterized by
\[
 \inner{s_\lambda}{s_\mu}=\delta_{\lambda\mu}.
\]
The involution $\omega$ of $\Lambda$ is determined by
\begin{equation}\label{eq:refaudit-omega}
\omega(p_m)=(-1)^{m-1}p_m;
\qquad \omega(s_\lambda)=s_{\lambda'}.
\end{equation}
For $f\in\Lambda_{\mathbb Q}$, the operator $f^\perp$ is determined by
$\langle f g,h\rangle=\langle g,f^\perp h\rangle$; in particular,
$p_m^\perp=m\,\partial/\partial p_m$.
For alphabets $U,V$, the ordinary and dual Cauchy kernels are
\begin{align}
 C(U,V)&=\prod_{u\in U,\,v\in V}(1-uv)^{-1}
       =\sum_\rho s_\rho(U)s_\rho(V),\label{eq:refaudit-ordinary-cauchy}\\
 D(U,V)&=\prod_{u\in U,\,v\in V}(1+uv)
       =\sum_\rho s_\rho(U)s_{\rho'}(V).\label{eq:refaudit-dual-cauchy}
\end{align}
The sums run over all partitions.

For any specified basis element $b$, the
notation $\Coef_b[F]$ means its coefficient in the basis expansion of $F$.
For formal Laurent series, the subscript specifies a monomial.

For a field $K$ of characteristic zero, write
$\Lambda_K=K\otimes_{\mathbb Z}\Lambda$.
We use square brackets for plethystic substitution: a formal alphabet
$A$ is specified by its power sums $p_m[A]$, and $f[A]$ is obtained
by substituting them into the power-sum expression of $f$.
Thus $p_m[A+B]=p_m[A]+p_m[B]$,
$p_m[AB]=p_m[A]p_m[B]$, and $p_m[-A]=-p_m[A]$.
A monomial $c$ in formal parameters satisfies $p_m[c]=c^m$;
in particular $p_m[q-1]=q^m-1$.  Note that plethystic negation $-A$
differs from replacing every variable by its ordinary negative.
All exponentials and coefficient extractions below are formal;
their action on any fixed homogeneous input is interpreted
coefficientwise in the stated auxiliary variables.

\subsection{Truncation and modular Schur functions}
\label{subsec:modular-schur-functions}
Truncation replaces each factor $(1-x_jt)^{-1}$ in the generating
series \eqref{eq:ordinary-complete-generating} by
$1+x_jt+\cdots+(x_jt)^{r-1}$.  Explicitly, we define
the \emph{truncated complete symmetric functions} by $\mathfrak{h}_0=1$ and
\begin{equation}\label{eq:truncated-complete-generating}
\sum_{m\geq0}\mathfrak{h}_m(X)t^m
=\prod_{j\geq1}\bigl(1+x_jt+\cdots+(x_jt)^{r-1}\bigr)
=\prod_{j\geq1}\frac{1-(x_jt)^r}{1-x_jt}.
\end{equation}
In Grinberg's notation, $\mathfrak{h}_m=G(r,m)$ is a \emph{Petrie symmetric
function} \cite[Definition~2.1(b)]{Grinberg}; in the notation of Fu and
Mei, $\mathfrak{h}_m=h_m^{[r-1]}$ is a \emph{truncated homogeneous symmetric
function} \cite[equation~(1)]{FuMei}.
For $m\geq1$, equivalently,
\[
\mathfrak{h}_m=\sum_{\substack{\nu\vdash m\\\nu_1<r}}m_\nu,
\]
since the largest part $\nu_1$ is the largest exponent occurring in
a monomial represented by $m_\nu$.  Since
$\Lambda=\mathbb Z[h_1,h_2,\ldots]$, the assignment $h_m\mapsto \mathfrak{h}_m$
extends uniquely to a graded ring endomorphism
\[
\Theta_r:\Lambda\longrightarrow\Lambda,
\qquad \Theta_r(h_m)=\mathfrak{h}_m.
\]
The \emph{modular Schur function} indexed by \(\lambda\) is defined by
\begin{equation}\label{eq:modular-schur-definition}
\mathfrak{s}_\lambda:=\Theta_r(s_\lambda).
\end{equation}
Applying $\Theta_r$ to \eqref{eq:refaudit-skew-jacobi-trudi} with
$\mu=\varnothing$ gives, for any $\ell\geq\ell(\lambda)$,
\begin{equation}\label{eq:refaudit-modular-jacobi-trudi}
\mathfrak{s}_\lambda
=\det\bigl(\mathfrak{h}_{\lambda_i-i+j}\bigr)_{1\leq i,j\leq\ell},
\qquad \mathfrak{h}_m=0\quad(m<0).
\end{equation}
For the dual Jacobi--Trudi identity for skew modular Schur functions,
see \cite[Theorem~3.2]{Hamel}.
The action on power sums follows directly from
\eqref{eq:truncated-complete-generating}.  By definition,
$\Theta_r(H(X;t))=\sum_{m\geq0}\mathfrak{h}_m(X)t^m$.  Taking the formal
logarithm in $\Lambda_{\mathbb Q}[[t]]$ gives
\begin{equation}\label{eq:refaudit-truncated-logarithm}
\log\!\left(\sum_{m\geq0}\mathfrak{h}_m t^m\right)
=\sum_{m\geq1}\frac{p_m t^m}{m}
-\sum_{m\geq1}\frac{p_{rm}t^{rm}}{m}.
\end{equation}
To compare \eqref{eq:refaudit-truncated-logarithm} with
\eqref{eq:ordinary-complete-generating}, extend
$\Theta_r$ coefficientwise to $\Lambda_{\mathbb Q}[[t]]$.  Since a
ring homomorphism commutes with the formal logarithm of a series with
constant term $1$, comparison of coefficients yields
\begin{equation}\label{eq:theta-power-sums}
\Theta_r(p_m)=
\begin{cases}
p_m,&r\nmid m,\\
(1-r)p_m,&r\mid m.
\end{cases}
\end{equation}
For $r=2$, the multiplier in \eqref{eq:theta-power-sums} is
$(-1)^{m-1}$, the multiplier in \eqref{eq:refaudit-omega}.  Hence
\begin{equation}\label{eq:prime-two-transpose}
\Theta_2=\omega,
\qquad \mathfrak{s}_\lambda=s_{\lambda'}.
\end{equation}

For every $r\geq2$, all multipliers in
\eqref{eq:theta-power-sums} are nonzero.  Since the power sums are
algebraically independent generators of $\Lambda_{\mathbb Q}$, the inverse
of $\Theta_r$ is the $\mathbb Q$-algebra endomorphism determined by
\[
\Theta_r^{-1}(p_m)=
\begin{cases}
p_m,&r\nmid m,\\
(1-r)^{-1}p_m,&r\mid m.
\end{cases}
\]
In particular, $\Theta_r(p_1)=p_1$.

For $n\geq1$, specialization to $x_1,\ldots,x_n$ means setting all
other variables equal to zero.  We write $f(x_1,\ldots,x_n)$
for the specialization of $f\in\Lambda_{\mathbb Q}$.
The Schur polynomial $s_\lambda(x_1,\ldots,x_n)$ vanishes when
$\ell(\lambda)>n$, and the remaining Schur polynomials form a
$\mathbb Q$-basis of the symmetric polynomial ring.
The elementary symmetric polynomials
\[
 e_i(x_1,\ldots,x_n)
 =\sum_{1\leq j_1<\cdots<j_i\leq n}x_{j_1}\cdots x_{j_i}
 \qquad(1\leq i\leq n)
\]
are algebraically independent generators of that ring.

For $N\geq0$, define the vector space
\begin{equation}\label{eq:finite-variable-space}
 V_{n,N}
 =\{f\in\mathbb Q[x_1,\ldots,x_n]^{\Sym_n}:
                  \deg_{x_i}f\leq N\text{ for every }i\}.
\end{equation}
Its monomial symmetric basis is indexed by the partitions satisfying
$\ell(\lambda)\leq n$ and $\lambda_1\leq N$.
By \cite[Chapter~I, (6.5)]{Macdonald}, the Schur-to-monomial
transition matrix is unitriangular in dominance order.
In $n$ variables, restricting to this rectangle preserves
unitriangularity, since dominance cannot increase the first part
and terms indexed by partitions of length greater than $n$ vanish.
Thus the same partitions index a Schur basis, and
\begin{equation}\label{eq:finite-space-basis}
 V_{n,N}
 =\operatorname{span}_{\mathbb Q}
   \{s_\lambda(x_1,\ldots,x_n):\ell(\lambda)\leq n,\ \lambda_1\leq N\},
 \qquad
 \dim V_{n,N}=\binom{N+n}{n}.
\end{equation}
The count records a weakly decreasing $n$-tuple with entries in
$\{0,\ldots,N\}$, including appended zeros.

\section{Bases in finitely many variables}
\label{sec:finite-variable-bases}

In this section we prove Theorem~\ref{thm:intro-finite-bases},
restated as Theorem~\ref{thm:finite-variable-bases}.
We express the truncated complete symmetric polynomials as
coefficients of a cyclic product and reduce the required linear
independence to the invertibility of a homogeneous substitution map.
We then construct a constant-rank matrix of relations among these
coefficients and use multiplication of sections on projective space
to prove that the map is an isomorphism.

\subsection{Reduction to homogeneous substitution}
\label{subsec:finite-homogeneous-substitution}

Fix integers $n\geq1$ and $r\geq2$, and put
\[
 \mathcal R_{n,r}
 =\{\lambda:\ell(\lambda)\leq n,\ \lambda_1\leq n(r-1)\}.
\]
The empty partition belongs to $\mathcal R_{n,r}$, with first part
understood to be zero.  In the notation of
\eqref{eq:finite-variable-space}, our aim is to obtain two bases of
$V_{n,n(r-1)}$.  The finite-variable generating function for $\mathfrak{h}_j$ is
\begin{equation}\label{eq:finite-truncated-generating}
 \sum_{j=0}^{n(r-1)}\mathfrak h_j(x_1,\ldots,x_n)t^j
 =\prod_{\nu=1}^n
       (1+x_\nu t+\cdots+x_\nu^{r-1}t^{r-1}).
\end{equation}
Write $\mathfrak h_\lambda=\prod_i\mathfrak h_{\lambda_i}$,
including the empty product $1$.

We first work over $\mathbb C$.  Let $\xi$ be a primitive $r$th
root of unity, let $c_0,\ldots,c_n$ be independent parameters,
and write the parameter vector as $\vec c=(c_0,\ldots,c_n)$.
Define
\begin{equation}\label{eq:cyclic-product-coefficients}
 P^{\vec c}(t)=\sum_{i=0}^n c_it^i,\qquad
 Q^{\vec c}(t)=\prod_{s=1}^{r-1}P^{\vec c}(\xi^st).
\end{equation}

For $0\leq j\leq n(r-1)$, define $q_j(\vec c)$ to be the coefficient of
$t^j$ in $Q^{\vec c}(t)$, so that
\[
 Q^{\vec c}(t)=\sum_{j=0}^{n(r-1)}q_j(\vec c)t^j.
\]
Each $q_j$ is homogeneous of ordinary degree $r-1$ in the parameters.
Since the automorphisms of $\mathbb Q(\xi)$ permute the factors in
\eqref{eq:cyclic-product-coefficients}, each $q_j$ belongs to
$\mathbb Q[c_0,\ldots,c_n]$.

For independent formal variables $y_0,\ldots,y_{n(r-1)}$ and an integer
$s\geq0$, define
\begin{equation}\label{eq:homogeneous-substitution}
 \theta_s:\mathbb C[y_0,\ldots,y_{n(r-1)}]_s
       \longrightarrow\mathbb C[c_0,\ldots,c_n]_{s(r-1)},
 \qquad F\longmapsto F(q_0,\ldots,q_{n(r-1)}).
\end{equation}
The subscripts denote ordinary homogeneous degree in the indicated
polynomial rings. The $y$-degree $s$ and the parameter degree $s(r-1)$ are distinct from the degree in the original variables $x_\nu$.  In particular,
\begin{equation}\label{eq:homogeneous-map-dimensions}
 \dim\mathbb C[y_0,\ldots,y_{n(r-1)}]_n
 =\dim\mathbb C[c_0,\ldots,c_n]_{n(r-1)}
 =\binom{nr}{n}.
\end{equation}

\begin{proposition}\label{prop:finite-basis-reduction}
The family $\{\mathfrak h_\lambda(x_1,\ldots,x_n):
\lambda\in\mathcal R_{n,r}\}$ is linearly independent over
$\mathbb Q$ if and only if $\theta_n$ is an isomorphism.
\end{proposition}

\begin{proof}
We compare the monomial basis of the source of $\theta_n$ with the
given family by two successive substitutions.  Append zeros to
each $\lambda\in\mathcal R_{n,r}$ to obtain $n$ parts.  Then
\[
 \left\{\prod_{i=1}^n y_{\lambda_i}:
             \lambda\in\mathcal R_{n,r}\right\}
\]
is exactly the degree-$n$ monomial basis: the $n$ subscripts of any
such monomial, arranged in weakly decreasing order, give a unique
$\lambda$.  The total degree counts the $n$ factors, including
those with subscript zero.

The first substitution, $y_j\mapsto q_j(\vec c)$, is $\theta_n$.
For the second, set $c_0=1$ and
$c_i=(-1)^ie_i(x_1,\ldots,x_n)$ for $1\leq i\leq n$.
Then $P^{\vec c}(t)$ becomes $\prod_{\nu=1}^n(1-x_\nu t)$, so
\[
 \begin{aligned}
 Q^{\vec c}(t)&\longmapsto
   \prod_{\nu=1}^n\prod_{s=1}^{r-1}(1-\xi^s x_\nu t)\\
 &=\prod_{\nu=1}^n
      (1+x_\nu t+\cdots+x_\nu^{r-1}t^{r-1})
  =\sum_{j=0}^{n(r-1)}\mathfrak h_j(x_1,\ldots,x_n)t^j,
 \end{aligned}
\]
where the last equality is \eqref{eq:finite-truncated-generating}.
Comparing coefficients of $t^j$ gives
$q_j(\vec c)\mapsto\mathfrak h_j(x_1,\ldots,x_n)$; in particular,
$q_0=c_0^{r-1}\mapsto1$.  Thus, for every
$\lambda\in\mathcal R_{n,r}$, the two substitutions give
\begin{equation}\label{eq:finite-basis-specialization}
 \prod_{i=1}^n y_{\lambda_i}
 \overset{\theta_n}{\longmapsto}
 \prod_{i=1}^n q_{\lambda_i}(\vec c)
 \longmapsto\mathfrak h_\lambda(x_1,\ldots,x_n).
\end{equation}

\emph{Necessity.}  Suppose the $\mathfrak h_\lambda$ are linearly
independent over $\mathbb Q$.  For coefficients
$a_\lambda\in\mathbb Q$, the second substitution in
\eqref{eq:finite-basis-specialization} gives
\[
 \sum_{\lambda\in\mathcal R_{n,r}}a_\lambda
       \prod_{i=1}^n q_{\lambda_i}(\vec c)=0
 \quad\Longrightarrow\quad
 \sum_{\lambda\in\mathcal R_{n,r}}a_\lambda\mathfrak h_\lambda=0.
\]
By the assumed independence, every $a_\lambda$ is zero.
Thus the matrix of $\theta_n$ in the monomial bases has linearly
independent columns over $\mathbb Q$.  Its entries are rational
because the $q_j$ have rational coefficients, so extending scalars
to $\mathbb C$ preserves its rank.  Hence $\theta_n$ is injective
and, by the equal dimensions in \eqref{eq:homogeneous-map-dimensions},
is an isomorphism.

\emph{Sufficiency.}  Conversely, suppose $\theta_n$ is an
isomorphism.  The products $\prod_{i=1}^n q_{\lambda_i}(\vec c)$ are
then linearly independent.  We must check that the second
substitution in \eqref{eq:finite-basis-specialization} preserves
this independence.

All these products have the same total parameter degree $n(r-1)$.
For nonnegative integers $b_1,\ldots,b_n$ with
$\sum_{i=1}^n b_i\leq n(r-1)$, setting $c_0=1$ sends
\[
 c_0^{\,n(r-1)-\sum_{i=1}^n b_i}
       c_1^{b_1}\cdots c_n^{b_n}
 \longmapsto c_1^{b_1}\cdots c_n^{b_n}.
\]
The exponents on the right determine the exponent of $c_0$ on the
left.  Distinct monomials therefore remain distinct, so this step
is injective on $\mathbb C[c_0,\ldots,c_n]_{n(r-1)}$.
Next, the algebraic independence of $e_1,\ldots,e_n$ makes
$c_i\mapsto(-1)^ie_i$ injective on $\mathbb C[c_1,\ldots,c_n]$.
The combined substitution is consequently injective on the
homogeneous component containing all the products.
By \eqref{eq:finite-basis-specialization}, their images are the
$\mathfrak h_\lambda$, which are therefore linearly independent
over $\mathbb C$, and hence over $\mathbb Q$.
\end{proof}

\subsection{Linear relations and constant rank}
\label{subsec:cyclic-constant-rank}

Assume $n\geq2$.  We construct linear relations among the $q_j$
whose matrix has constant rank for $\vec c\ne0$.

For coefficient calculations, extend the two coefficient sequences
by the conventions $c_i=0$ for $i<0$ or $i>n$, and $q_j=0$ for
$j<0$ or $j>n(r-1)$.  The latter agrees with the degree bound for
$Q^{\vec c}$ in \eqref{eq:cyclic-product-coefficients}.
Define the matrix $\mathcal B(\vec c)$, with rows
$1\leq k\leq nr-1$ and columns $0\leq j\leq n(r-1)$, by
\begin{equation}\label{eq:cyclic-linear-matrix}
 \mathcal B(\vec c)_{k,j}=
 \begin{cases}
 c_{k-j},&r\nmid k,\\
 \bigl(r(k-j)-k\bigr)c_{k-j},&r\mid k.
 \end{cases}
\end{equation}
Write $\mathsf T$ for matrix transpose, and let
$q(\vec c)=(q_0(\vec c),\ldots,q_{n(r-1)}(\vec c))^{\mathsf T}$.
We first derive the two types of coefficient relations encoded
by \eqref{eq:cyclic-linear-matrix}.

For a fixed parameter vector $\vec c$, define the norm polynomial
and, for $f\in\mathbb C(t)$, the cyclic average by
\begin{equation}\label{eq:refaudit-cyclic-norm-average}
 \begin{aligned}
 \mathscr N(P^{\vec c})(t)
   &=\prod_{s=0}^{r-1}P^{\vec c}(\xi^st)
     =P^{\vec c}(t)Q^{\vec c}(t),\\
 \Pi(f)(t)&=\frac1r\sum_{s=0}^{r-1}f(\xi^st).
 \end{aligned}
\end{equation}
For every integer $k\geq0$, evaluating the geometric sum in
\eqref{eq:refaudit-cyclic-norm-average} gives
\begin{equation}\label{eq:cyclic-average-monomials}
 \Pi(t^k)=\frac{t^k}{r}\sum_{s=0}^{r-1}\xi^{sk}
 =\begin{cases}
 t^k,&r\mid k,\\
 0,&r\nmid k.
 \end{cases}
\end{equation}
By linearity, $\Pi$ retains precisely the terms of a polynomial
whose degrees are divisible by $r$, with their coefficients unchanged.

Rotation by $\xi$
permutes the factors of $\mathscr N(P^{\vec c})$, so
$\mathscr N(P^{\vec c})\in\mathbb C[t^r]$.
We have
\begin{equation}\label{eq:cyclic-average-derivative}
 \begin{aligned}
 \Pi\bigl(t(P^{\vec c})'(t)Q^{\vec c}(t)\bigr)
 &=\frac1r\sum_{s=0}^{r-1}
       \xi^st(P^{\vec c})'(\xi^st)Q^{\vec c}(\xi^st)\\
 &=\frac{t}{r}\mathscr N(P^{\vec c})'(t).
 \end{aligned}
\end{equation}
For the second equality, \eqref{eq:cyclic-product-coefficients} and
$\xi^r=1$ give
$Q^{\vec c}(\xi^st)=\prod_{\substack{0\leq j<r\\j\ne s}}P^{\vec c}(\xi^jt)$
for $0\leq s<r$.
Applying the product rule to \eqref{eq:refaudit-cyclic-norm-average},
using the chain rule
$\frac{d}{dt}P^{\vec c}(\xi^st)=\xi^s(P^{\vec c})'(\xi^st)$,
and multiplying by $t/r$ gives the sum in
\eqref{eq:cyclic-average-derivative}.

For $1\leq k\leq nr-1$, the coefficient of $t^k$ in
$P^{\vec c}(t)Q^{\vec c}(t)$ is $\sum_{i=0}^n c_iq_{k-i}$.
If $r\nmid k$, \eqref{eq:refaudit-cyclic-norm-average} and
$\mathscr N(P^{\vec c})\in\mathbb C[t^r]$ give
\[
 \sum_{i=0}^n c_iq_{k-i}=0.
\]
If $r\mid k$, \eqref{eq:cyclic-average-monomials} shows that
$\Pi$ leaves the coefficient of $t^k$ unchanged.
Since $t(P^{\vec c})'(t)=\sum_{i=0}^n i c_it^i$, taking the
coefficient of $t^k$ in \eqref{eq:cyclic-average-derivative} and
using \eqref{eq:refaudit-cyclic-norm-average} gives
\[
 \sum_{i=0}^n i c_iq_{k-i}
 =\frac{k}{r}\sum_{i=0}^n c_iq_{k-i}.
\]
Here differentiating $\mathscr N(P^{\vec c})(t)$ and multiplying
by $t/r$ multiplies its $t^k$ coefficient by $k/r$.
Multiplying by $r$ and rearranging therefore gives
$\sum_{i=0}^n(ri-k)c_iq_{k-i}=0$.
In both cases, setting $j=k-i$ and using
\eqref{eq:cyclic-linear-matrix} yields
\begin{equation}\label{eq:refaudit-cyclic-kernel-vector}
 \mathcal B(\vec c)q(\vec c)=0.
\end{equation}
Equation~\eqref{eq:refaudit-cyclic-kernel-vector} shows that
$q(\vec c)\in\ker\mathcal B(\vec c)$;
Proposition~\ref{prop:cyclic-constant-rank} shows that, for
$\vec c\ne0$, every vector in this kernel is a scalar multiple of
$q(\vec c)$.

\begin{proposition}\label{prop:cyclic-constant-rank}
For every nonzero $\vec c\in\mathbb C^{n+1}$,
$\ker\mathcal B(\vec c)=\mathbb Cq(\vec c)$ and
$\operatorname{rank}\mathcal B(\vec c)=n(r-1)$.
\end{proposition}

\begin{proof}
Fix $\vec c\ne0$. Let $g=(g_0,\ldots,g_{n(r-1)})^{\mathsf T}$ satisfy
$\mathcal B(\vec c)g=0$, and put
$G(t)=\sum_{j=0}^{n(r-1)}g_jt^j$.
We use the two types of rows to derive identities for $G$, then
show that $G$ is a scalar multiple of $Q^{\vec c}$.
For $1\leq k\leq nr-1$,
the $k$th component of $\mathcal B(\vec c)g$ is
\begin{equation}\label{eq:cyclic-matrix-coefficient-action}
 \begin{cases}
 \Coef_{t^k}\bigl[P^{\vec c}G\bigr],&r\nmid k,\\
 \Coef_{t^k}\bigl[rt(P^{\vec c})'G-t(P^{\vec c}G)'\bigr],&r\mid k.
 \end{cases}
\end{equation}
Indeed, expanding the products gives
\[
 \begin{aligned}
 \Coef_{t^k}\bigl[P^{\vec c}G\bigr]
 &=\sum_{j=0}^{n(r-1)}c_{k-j}g_j,\\
 \Coef_{t^k}\bigl[rt(P^{\vec c})'G-t(P^{\vec c}G)'\bigr]
 &=\sum_{j=0}^{n(r-1)}\bigl(r(k-j)-k\bigr)c_{k-j}g_j.
 \end{aligned}
\]
The factors $r(k-j)$ and $k$ in the second line arise from applying
$rt\,d/dt$ to a term of degree $k-j$ in $P^{\vec c}$ and
$t\,d/dt$ to a term of degree $k$ in $P^{\vec c}G$, respectively.
Comparing with \eqref{eq:cyclic-linear-matrix}
proves \eqref{eq:cyclic-matrix-coefficient-action}.

Since $\mathcal B(\vec c)g=0$, the rows with $r\nmid k$ in
\eqref{eq:cyclic-matrix-coefficient-action} force the $t^k$
coefficients of $P^{\vec c}G$ to vanish at all degrees not divisible
by $r$: its degree is at most $nr$, and the omitted endpoints
$0,nr$ are both divisible by $r$. Thus
\begin{equation}\label{eq:cyclic-kernel-invariant-product}
 P^{\vec c}(t)G(t)\in\mathbb C[t^r].
\end{equation}

The rows with $r\mid k$ in \eqref{eq:cyclic-matrix-coefficient-action}
give zero $t^k$ coefficient in $rt(P^{\vec c})'G-t(P^{\vec c}G)'$
for $1\leq k\leq nr-1$. This polynomial also has degree at most
$nr$. At the two omitted endpoints, its constant coefficient is
zero because both terms contain $t$, and its $t^{nr}$ coefficient
is $(rn-nr)c_ng_{n(r-1)}=0$. These calculations include the cases
of smaller degree. Applying \eqref{eq:cyclic-average-monomials}
therefore gives
\begin{equation}\label{eq:cyclic-kernel-average}
 \Pi\bigl(rt(P^{\vec c})'G-t(P^{\vec c}G)'\bigr)=0.
\end{equation}

Since $\vec c\ne0$, $P^{\vec c}$ is a nonzero polynomial, and
\eqref{eq:cyclic-product-coefficients} and
\eqref{eq:refaudit-cyclic-norm-average} show that $Q^{\vec c}$ and
$\mathscr N(P^{\vec c})$ are nonzero as well. Define the rational function
\[
 h=\frac{G}{Q^{\vec c}}
   =\frac{P^{\vec c}G}{\mathscr N(P^{\vec c})}.
\]
It remains to show that $h$ is constant.
By \eqref{eq:cyclic-kernel-invariant-product}, the numerator
$P^{\vec c}G$ is unchanged under $t\mapsto\xi t$; the same holds
for the product $\mathscr N(P^{\vec c})$ in
\eqref{eq:refaudit-cyclic-norm-average}. Hence $h(\xi t)=h(t)$.
Also, $\mathscr N(P^{\vec c})h=P^{\vec c}G$ belongs to
$\mathbb C[t^r]$ by \eqref{eq:cyclic-kernel-invariant-product}, and
$t\,d/dt$ preserves this space. Thus
\eqref{eq:refaudit-cyclic-norm-average} and
\eqref{eq:cyclic-average-monomials} give, for $f\in\mathbb C(t)$,
\begin{equation}\label{eq:cyclic-kernel-averaging}
 \begin{aligned}
 \Pi(hf)&=\frac1r\sum_{s=0}^{r-1}h(\xi^st)f(\xi^st)
          =h\Pi(f),\\
 \Pi\bigl(t(\mathscr N(P^{\vec c})h)'\bigr)
       &=t(\mathscr N(P^{\vec c})h)'.
 \end{aligned}
\end{equation}

Substituting $G=Q^{\vec c}h$ into \eqref{eq:cyclic-kernel-average}
and using \eqref{eq:cyclic-kernel-averaging} and
\eqref{eq:cyclic-average-derivative} in turn, we obtain
\[
 \begin{aligned}
 0&=\Pi\bigl(rt(P^{\vec c})'Q^{\vec c}h
        -t(\mathscr N(P^{\vec c})h)'\bigr)\\
  &=rh\,\Pi\bigl(t(P^{\vec c})'Q^{\vec c}\bigr)
        -t\bigl(\mathscr N(P^{\vec c})h\bigr)'\\
  &=ht\mathscr N(P^{\vec c})'
        -t\bigl(\mathscr N(P^{\vec c})h\bigr)'\\
  &=-t\mathscr N(P^{\vec c})h'.
 \end{aligned}
\]
The last equality uses the product rule. Cancelling the nonzero
factor $t\mathscr N(P^{\vec c})$ in
$\mathbb C(t)$ gives $h'=0$.  In characteristic zero this means
that $h=a$ for some $a\in\mathbb C$. Thus $G=aQ^{\vec c}$, and
comparing coefficients gives $g=a\,q(\vec c)$.
Together with \eqref{eq:refaudit-cyclic-kernel-vector}, this proves
$\ker\mathcal B(\vec c)=\mathbb Cq(\vec c)$.
Since $Q^{\vec c}\ne0$, its coefficient vector $q(\vec c)$ is
nonzero, so the kernel has dimension one. The rank--nullity theorem
now gives $\operatorname{rank}\mathcal B(\vec c)=n(r-1)$.
\end{proof}

We now use Proposition~\ref{prop:cyclic-constant-rank} to construct
an exact sequence of vector bundles.
Readers unfamiliar with vector bundles and sheaf cohomology may
consult Appendix~\ref{subsec:appendix-background} for the background
used below; see also \cite[Chapter~II, \S5]{Hartshorne}.

Regard nonzero parameter vectors up to nonzero scalar multiplication.
They form the projective space $\mathbb P^n_{\mathbb C}$, with
homogeneous coordinates $[c_0:\cdots:c_n]$.
Write $\mathcal O(-1)$ for the tautological line bundle: its fiber,
the vector space attached to $[c_0:\cdots:c_n]$, is the line
$\mathbb C\vec c\subset\mathbb C^{n+1}$.
The rank of a vector bundle is the dimension of its fibers.
Let $\mathcal O(1)$ be the dual of $\mathcal O(-1)$, and for each integer $d$ let
$\mathcal O(d)$ denote the $d$th tensor power of $\mathcal O(1)$,
using dual powers when $d<0$; the trivial line bundle is
$\mathcal O=\mathcal O(0)$.
For $d\ge0$, global sections of $\mathcal O(d)$ correspond to
homogeneous polynomials of degree $d$ in $c_0,\ldots,c_n$.
For an integer $m\ge0$, the notation $\mathcal O(d)^m$ denotes the direct sum
of $m$ copies of this line bundle.

Put $U=\mathbb C^{nr-1}$ and $V=\mathbb C^{n(r-1)+1}$, so that
$U\otimes\mathcal O(-1)=\mathcal O(-1)^{nr-1}$ and
$V\otimes\mathcal O=\mathcal O^{n(r-1)+1}$.
Each entry of $\mathcal B(\vec c)$ is linear in the parameters and
therefore defines a map $\mathcal O(-1)\to\mathcal O$.
Thus its transpose matrix defines a bundle map
$\varphi:U\otimes\mathcal O(-1)\to V\otimes\mathcal O$.
Similarly, the degree-$(r-1)$ polynomials $q_j$ define a bundle map
$\varepsilon:V\otimes\mathcal O\to\mathcal O(r-1)$, sending
the $j$th constant basis section of $V\otimes\mathcal O$ to $q_j$ for
$0\le j\le n(r-1)$. After choosing a representative $\vec c$ of
a projective point, the map on its fiber is represented by the
row vector $q(\vec c)^{\mathsf T}$.
Since $Q^{\vec c}\ne0$ for every $\vec c\ne0$, its coefficients
$q_j(\vec c)$ are not all zero. Hence $\varepsilon$ maps each
fiber onto the one-dimensional target fiber and is surjective.

Transposing \eqref{eq:refaudit-cyclic-kernel-vector} gives
\[
 q(\vec c)^{\mathsf T}\mathcal B(\vec c)^{\mathsf T}
   =\bigl(\mathcal B(\vec c)q(\vec c)\bigr)^{\mathsf T}=0,
\]
so $\varepsilon\varphi=0$ and
$\operatorname{im}\varphi\subseteq\ker\varepsilon$.
By Proposition~\ref{prop:cyclic-constant-rank}, the transpose matrix
also has rank $n(r-1)$ at every point. Thus the kernel and image
of $\varphi$ are vector bundles, and $\operatorname{im}\varphi$
has rank $n(r-1)$. Since $\varepsilon$ is surjective,
$\ker\varepsilon$ has the same rank, namely $n(r-1)+1-1$.
The inclusion is therefore an equality on every fiber, giving
$\operatorname{im}\varphi=\ker\varepsilon$.
Set $E=\ker\varphi$, the kernel bundle of the transpose map.
We obtain the exact sequence
\begin{equation}\label{eq:cyclic-bundle-presentation}
 0\longrightarrow E\longrightarrow
 \mathcal O(-1)^{nr-1}\xrightarrow{\varphi}
 \mathcal O^{n(r-1)+1}\xrightarrow{\varepsilon}
 \mathcal O(r-1)\longrightarrow0.
\end{equation}
Here exactness means that the image of each map equals the kernel
of the next.

The construction has not assumed that the sections $q_j$ are
linearly independent: the surjectivity of $\varepsilon$ uses only
that they have no common zero.
For $n\ge2$, \eqref{eq:cyclic-bundle-presentation} will be used
in the next subsection to prove that the map $\theta_n$ in
\eqref{eq:homogeneous-substitution} is an isomorphism.

\subsection{Proof of Theorem~\ref{thm:intro-finite-bases}}
\label{subsec:multiplication-and-bases}

We first state an isomorphism result for multiplication of sections.
We then apply it to \eqref{eq:cyclic-bundle-presentation} to obtain
the two bases in Theorem~\ref{thm:intro-finite-bases}.

Let $K$ be a field of characteristic zero, let $n,a$ be integers with $n\geq2$ and
$a\geq1$, and write $L=\mathcal O_{\mathbb P^n_K}(a)$.
Let $U,V$ be $K$-vector spaces of dimensions $na+n-1$ and
$na+1$, respectively, and consider the sequence of vector bundles
on $\mathbb P^n_K$
\begin{equation}\label{eq:general-linear-presentation}
 0\longrightarrow E\longrightarrow U\otimes\mathcal O(-1)
 \xrightarrow{\varphi}V\otimes\mathcal O
 \xrightarrow{\varepsilon}L\longrightarrow0.
\end{equation}
For an integer $s\geq0$, write $L^s=L^{\otimes s}$ and let
$\operatorname{Sym}^sV$ denote the $s$th symmetric power of $V$.
The space $H^0(\mathbb P^n_K,L^s)$ consists of global sections of
$L^s$. Multiplication of the sections $\varepsilon(v\otimes1)$, for $v\in V$,
defines
\begin{equation}\label{eq:section-multiplication-map}
 \begin{aligned}
 \alpha_s:\operatorname{Sym}^sV
     &\longrightarrow H^0(\mathbb P^n_K,L^s),\\
 v_1\cdots v_s&\longmapsto
        \varepsilon(v_1\otimes1)\cdots\varepsilon(v_s\otimes1)
        \qquad(v_1,\ldots,v_s\in V).
 \end{aligned}
\end{equation}
For $s=0$, both products are interpreted as $1$.

\begin{lemma}\label{lem:linear-presentation-multiplication}
If \eqref{eq:general-linear-presentation} is exact, then the map
$\alpha_n$ in \eqref{eq:section-multiplication-map} is an isomorphism.
\end{lemma}

The proof of Lemma \ref{lem:linear-presentation-multiplication} is given in Appendix~\ref{subsec:appendix-lemma-proof}. We now prove Theorem~\ref{thm:intro-finite-bases}, restated below.
Return to the integers $n\geq1$ and $r\geq2$,
the rectangle $\mathcal R_{n,r}$, and the space $V_{n,n(r-1)}$.

\begin{theorem}[Finite-variable bases]\label{thm:finite-variable-bases}
For every $n\geq1$ and $r\geq2$, the two families
\[
 \{\mathfrak h_\lambda(x_1,\ldots,x_n):\lambda\in\mathcal R_{n,r}\},
 \qquad
 \{\mathfrak s_\lambda(x_1,\ldots,x_n):\lambda\in\mathcal R_{n,r}\}
\]
are each a $\mathbb Q$-basis of $V_{n,n(r-1)}$.
\end{theorem}

\begin{proof}
For $n\geq2$, apply
Lemma~\ref{lem:linear-presentation-multiplication} over $K=\mathbb C$, with
$a=r-1$, $U=\mathbb C^{nr-1}$, and $V=\mathbb C^{n(r-1)+1}$
to \eqref{eq:cyclic-bundle-presentation}.
If $v_j$ is the $j$th standard basis vector of $V$, then
$\varepsilon(v_j\otimes1)=q_j$. Identify
$\operatorname{Sym}^nV$ with $\mathbb C[y_0,\ldots,y_{n(r-1)}]_n$
by $v_j\mapsto y_j$, and
$H^0(\mathbb P^n_{\mathbb C},L^n)$ with
$\mathbb C[c_0,\ldots,c_n]_{n(r-1)}$.
Then \eqref{eq:section-multiplication-map} sends each product of
the $y_j$ to the same product of the $q_j$, so $\alpha_n$ is exactly
$\theta_n$ in \eqref{eq:homogeneous-substitution}.
Lemma~\ref{lem:linear-presentation-multiplication} makes $\theta_n$
an isomorphism, so
Proposition~\ref{prop:finite-basis-reduction} proves the independence
of the truncated complete products over $\mathbb Q$.
For $n=1$, \eqref{eq:finite-truncated-generating} gives
$\mathfrak h_j(x_1)=x_1^j$ for $0\leq j\leq r-1$, which proves
the independence of $\{\mathfrak h_\lambda(x_1):
\lambda\in\mathcal R_{1,r}\}$ directly.

By \cite[Proposition~1.5(i), pp.~576--577]{Walker}, independence
of these truncated complete products is equivalent to independence
of the modular Schur polynomials indexed by the same rectangle.
Both families belong to $V_{n,n(r-1)}$.  Indeed, each
$\mathfrak h_j(x_1,\ldots,x_n)$ has degree at most $r-1$ in each
variable, and $\mathfrak h_\lambda$ has at most $n$ factors.
For $\mathfrak s_\lambda$, append zero parts to length $n$ and use
\eqref{eq:refaudit-modular-jacobi-trudi} with $\ell=n$; every nonzero term is again a
product of $n$ truncated complete polynomials.
Each family has $|\mathcal R_{n,r}|=\binom{nr}{n}$ elements,
equal to $\dim V_{n,n(r-1)}$ by \eqref{eq:finite-space-basis}.
Independence therefore proves both basis assertions.
\end{proof}

Theorem~\ref{thm:finite-variable-bases} proves
\cite[Conjecture~1.4 and its equivalent form in Proposition~1.5(i),
pp.~576--577]{Walker}.
The matrix formulation in \cite[Proposition~1.5(ii)]{Walker}
will be used in Corollary~\ref{cor:rectangular-principal-matrices}
to relate the bases in Theorem~\ref{thm:finite-variable-bases}
to the modular Schur coefficients.

\section{Heisenberg operators and core--quotient coordinates}
\label{sec:universal-blocks}

We return to the stable ring $\Lambda_{\mathbb Q}$ and an integer
$r\geq2$.  We first fix the abacus conventions and recall
block preservation.  We then describe the Heisenberg action in
core--quotient coordinates, factor the half-vertex operators, and
compute the induced endomorphisms and their reproducing kernels.

\subsection{Abaci, cores, and quotients}
\label{subsec:abacus-core-quotient}

The Schur coefficients of $\mathfrak{s}_\lambda$ are organized by the core and
quotient of $\lambda$.  We describe these objects through an abacus.
Choose an integer $b\geq\ell(\lambda)$ divisible by $r$, append zero
parts to $\lambda$, and form the beta-set
\[
\mathfrak B_b(\lambda)=\{\lambda_j+b-j:1\leq j\leq b\}.
\]
This is a set of $b$ distinct nonnegative integers.  Conversely, if
$q_1>\cdots>q_b$ are its elements, then
$\lambda_j=q_j-b+j$, so the beta-set retains all information about
$\lambda$.  Replacing $b$ by $b+r$ gives
\begin{equation}\label{eq:beta-set-b-plus-r}
\mathfrak B_{b+r}(\lambda)
=\{q+r:q\in \mathfrak B_b(\lambda)\}\sqcup\{0,1,\ldots,r-1\}.
\end{equation}
The $r$-abacus has runners numbered $0,\ldots,r-1$.  Position
$q=ar+i$ lies at level $a$ on runner $i$; positions increase from top
to bottom, so $q-r$ is immediately above $q$ on the same runner.

Placing a bead at each position in $\mathfrak B_b(\lambda)$ gives an
\emph{abacus display} of $\lambda$.  If $q\geq r$ and
the position $q-r$ is unoccupied, replace $q$ by $q-r$ in the beta-set
and reorder its elements.  The corresponding partition is obtained from
$\lambda$ by removing a connected skew diagram of $r$ boxes containing no
$2\times2$ square, that is, an $r$-ribbon.  On the abacus this operation
simply moves the bead at $q$ one level upward.

\begin{figure}[!htbp]
\centering
\begin{tikzpicture}[
  x=1cm,
  y=1cm,
  space/.style={circle,draw=black!70,line width=0.45pt,fill=white,
                text=black!75,minimum size=6.4mm,inner sep=0pt,font=\small},
  bead/.style={circle,draw=black,line width=0.45pt,fill=black,text=white,
               minimum size=6.4mm,inner sep=0pt,font=\small}
]
\path[use as bounding box] (-2.55,1.95) rectangle (9.55,-5.45);
% Column headings.
\node[font=\small\bfseries] at (0.9,1.65) {Before};
\node[font=\small\bfseries] at (6.1,1.65) {After};

% Beta-sets.
\node[font=\small\bfseries,anchor=east] at (-0.45,0.55) {Beta-set};
\node[font=\small] at (0.9,0.55) {$\{5,3,1\}$};
\node[font=\small] at (6.1,0.55) {$\{5,1,0\}$};
\draw[->,thick] (2.25,0.55)--(4.75,0.55);

% Abaci.
\node[font=\small\bfseries,anchor=east] at (-0.45,-1.55) {Abacus};
\foreach \x/\i in {0/0,0.9/1,1.8/2} {
  \draw[densely dashed,black!60,line width=0.45pt] (\x,-0.66)--(\x,-2.52);
  \node[font=\small] at (\x,-0.42) {$\i$};
}
\node[space] at (0,-0.95) {0};
\node[bead]  at (0.9,-0.95) {1};
\node[space] at (1.8,-0.95) {2};
\node[bead]  at (0,-2.15) {3};
\node[space] at (0.9,-2.15) {4};
\node[bead]  at (1.8,-2.15) {5};
\foreach \x/\i in {5.2/0,6.1/1,7/2} {
  \draw[densely dashed,black!60,line width=0.45pt] (\x,-0.66)--(\x,-2.52);
  \node[font=\small] at (\x,-0.42) {$\i$};
}
\node[bead]  at (5.2,-0.95) {0};
\node[bead]  at (6.1,-0.95) {1};
\node[space] at (7,-0.95) {2};
\node[space] at (5.2,-2.15) {3};
\node[space] at (6.1,-2.15) {4};
\node[bead]  at (7,-2.15) {5};
\draw[->,thick] (2.45,-1.55)--(4.65,-1.55)
  node[midway,above,font=\small] {$3\mapsto0$};

% Young diagrams.
\node[font=\small\bfseries,anchor=east] at (-0.45,-4.25) {Young diagram};
\node[font=\small] at (0.75,-3.35) {$(3,2,1)$};
\foreach \x in {0,0.5,1} {
  \draw[fill=black!5] (\x,-3.65) rectangle ({\x+0.5},-4.15);
}
\foreach \x in {0,0.5} {
  \draw[fill=red!25] (\x,-4.15) rectangle ({\x+0.5},-4.65);
}
\draw[fill=red!25] (0,-4.65) rectangle (0.5,-5.15);
\node[font=\small] at (5.95,-3.35) {$(3)$};
\foreach \x in {5.2,5.7,6.2} {
  \draw[fill=black!5] (\x,-3.65) rectangle ({\x+0.5},-4.15);
}
\draw[->,thick] (2.25,-4.25)--(4.75,-4.25)
  node[midway,above,font=\small] {remove};
\end{tikzpicture}
\caption{The bead move $3\mapsto0$ in beta-set, abacus, and Young
diagram notation.  Black circles mark occupied positions, open circles
mark unoccupied positions, and the red boxes form the removed $3$-ribbon.}
\label{fig:three-ribbon-removal}
\end{figure}
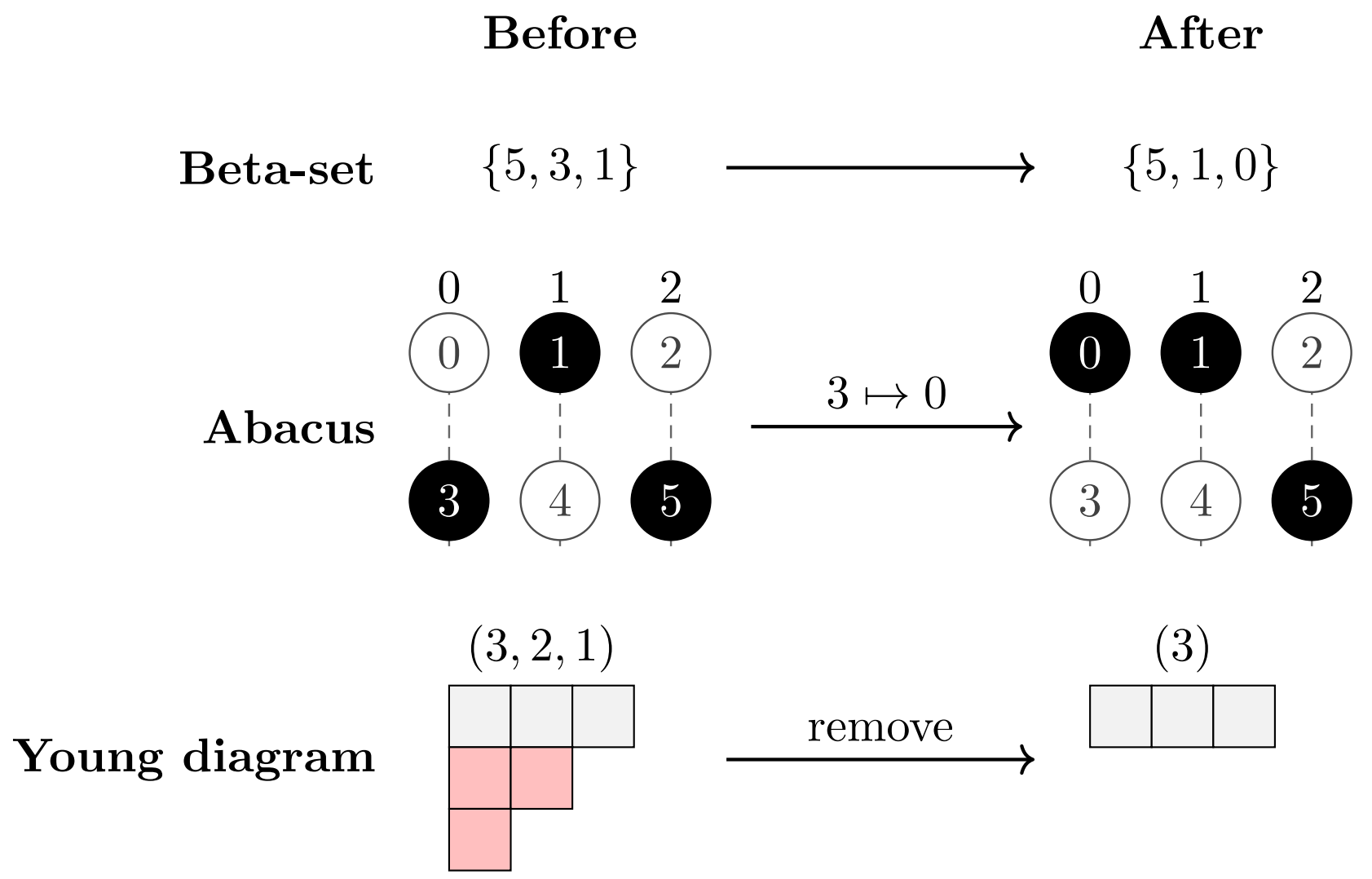

For example, take $r=3$, $\lambda=(3,2,1)$, and $b=3$.  Since
$\mathfrak B_3(\lambda)=\{5,3,1\}$ and $0$ is unoccupied, we may replace
$q=3$ by $q-r=0$, changing the beta-set from $\{5,3,1\}$ to
$\{5,1,0\}$.  Figure~\ref{fig:three-ribbon-removal} compares the
beta-set, abacus, and Young diagram descriptions of this operation.

Slide every bead upward until no such move is possible.  If runner
$i$ contains $b_i$ beads, its final occupied positions must be
$i,i+r,\ldots,i+(b_i-1)r$.  Hence the packed bead configuration is
independent of the order of the moves and is the beta-set of a uniquely determined
partition $\kappa$.  This partition is the \emph{$r$-core} of
$\lambda$, and we write $\kappa=\operatorname{core}_r(\lambda)$. For $\lambda=(3,2,1)$ in Figure~\ref{fig:three-ribbon-removal},
the further move $5\mapsto2$ gives the beta-set
$\{2,1,0\}=\mathfrak B_3(\varnothing)$, so
$\operatorname{core}_3(\lambda)=\varnothing$. A
partition $\lambda$ is itself called an \emph{$r$-core} if
$\operatorname{core}_r(\lambda)=\lambda$, equivalently, if no $r$-ribbon can be removed from $\lambda$.

To describe how the quotient changes with the bead number, temporarily
allow any integer $b\geq\ell(\lambda)$.  For each runner $i$, write
\[
C_i=\{(q-i)/r:q\in \mathfrak B_b(\lambda),\ q\equiv i\pmod r\}
=\{c_{i,1}>\cdots>c_{i,b_i}\}.
\]
Let
\[
Q_b(\lambda)
=\bigl(\alpha^{(0)},\ldots,\alpha^{(r-1)}\bigr),
\]
where the component on runner $i$ is the partition with parts
\[
\alpha^{(i)}_j=c_{i,j}-b_i+j
\qquad(1\leq j\leq b_i).
\]

The tuple $Q_b(\lambda)$ depends on $b$.  By
\cite[Section~1.1.1, p.~442]{Fayers},
\begin{equation}\label{eq:quotient-charge-shift}
Q_{b+1}(\lambda)
=\bigl(\alpha^{(r-1)},\alpha^{(0)},\ldots,
       \alpha^{(r-2)}\bigr).
\end{equation}

For $\lambda=(3,2,1)$ in Figure~\ref{fig:three-ribbon-removal},
\[
 Q_3(\lambda)=\bigl((1),\emptypar,(1)\bigr),
 \qquad
 Q_4(\lambda)=\bigl((1),(1),\emptypar\bigr),
\]
where $\mathfrak B_4(\lambda)=\{6,4,2,0\}$.
For a fixed bead number, \eqref{eq:quotient-charge-shift} applies
simultaneously to every partition represented with that number of beads.

Applying \eqref{eq:quotient-charge-shift} $r$ times gives
$Q_{b+r}(\lambda)=Q_b(\lambda)$.  We therefore define the ordered
\emph{$r$-quotient} of $\lambda$ by
\begin{equation}\label{eq:standard-quotient-tuple}
Q_b(\lambda)
=\bigl(\lambda^{(0)},\ldots,\lambda^{(r-1)}\bigr)
\qquad\bigl(b\geq\ell(\lambda),\ r\mid b\bigr).
\end{equation}
Each $\lambda^{(i)}$ is independent of the choice of $b$ under
these conditions, since any two such bead numbers differ by a multiple of $r$.

Define the $r$-weight by
$w(\lambda)=\sum_{i=0}^{r-1}|\lambda^{(i)}|$.

Consider the core--quotient map
\[
 \lambda\longmapsto
 \bigl(\operatorname{core}_r(\lambda);
       \lambda^{(0)},\ldots,\lambda^{(r-1)}\bigr)
\]
from partitions to pairs consisting of an $r$-core and an ordered
$r$-tuple of partitions.
Proposition~\ref{prop:core-quotient-correspondence} is standard; see
\cite[Chapter~I, \S1, Example~8(b), (c), and~(g)]{Macdonald}.

\begin{proposition}
\label{prop:core-quotient-correspondence}
For every integer $r\geq2$, the core--quotient map is a bijection.
\end{proposition}

For every partition $\lambda$, put
$\kappa=\operatorname{core}_r(\lambda)$.  Then
\begin{equation}\label{eq:core-weight-size}
 |\lambda|=|\kappa|+r\sum_{i=0}^{r-1}|\lambda^{(i)}|,
 \qquad
 w(\lambda)=\sum_{i=0}^{r-1}|\lambda^{(i)}|
 =\frac{|\lambda|-|\kappa|}{r}.
\end{equation}
Moreover, $w(\lambda)$ is the number of $r$-ribbons removed in
passing from $\lambda$ to its core.

With the bead number divisible by $r$, conjugation reverses the
runner order and conjugates each quotient component
\cite[Section~1.1.1, p.~442]{Fayers}:
\begin{equation}\label{eq:quotient-conjugation}
(\lambda')^{(i)}=(\lambda^{(r-1-i)})'
\qquad(0\leq i<r).
\end{equation}

Label the beads of $\mathfrak B_b(\lambda)$ by $1,\ldots,b$ in
decreasing order of their positions.  Pack every runner upward,
keeping the order of its beads fixed.  The \emph{$r$-sign}
$\eps_r(\lambda)$ is the sign of the permutation obtained by reading
the labels in decreasing order of the final bead positions.

Under
\eqref{eq:beta-set-b-plus-r}, every old bead is shifted down by one
level and one new bead is added at the top of each runner;
the old beads undergo the same permutation as before this enlargement.
The new beads remain above all old beads, and their mutual order is
unchanged by packing.  Hence $\eps_r(\lambda)$ is independent of the
admissible choice of $b$.

Put $w=w(\lambda)$.  If $r$-ribbons of heights $h_1,\ldots,h_w$ are removed from
$\lambda$ to its core, then
\begin{equation}\label{eq:r-sign-height}
\eps_r(\lambda)=(-1)^{h_1+\cdots+h_w}.
\end{equation}
Indeed, the bead move corresponding to an $r$-ribbon of height $h$ crosses
exactly $h$ labelled beads.  Since $\eps_r(\lambda)$ is defined by
the initial and packed bead orders, $\eps_r(\lambda)$ and hence the
parity in \eqref{eq:r-sign-height} are independent of the removal
sequence.

For $\lambda=(3,2,1)$ in Figure~\ref{fig:three-ribbon-removal},
label the beads at $5,3,1$ by $A,B,C$ in that order.  After packing, their positions in decreasing order are
$2,1,0$, carrying the labels $A,C,B$.  Thus one transposition has
occurred and $\eps_3((3,2,1))=-1$.

\subsection{Combinatorial blocks and block preservation}
\label{subsec:combinatorial-blocks}

Let $\kappa$ be any partition and let $w\geq0$.  Define
\[
\mathcal B_{\kappa,w}^{(r)}
=\{\lambda:\operatorname{core}_r(\lambda)=\kappa,
 \ w(\lambda)=w\}
\]
and call it the \emph{$r$-block} indexed by $(\kappa,w)$.  If $\kappa$
is not an $r$-core, then $\mathcal B_{\kappa,w}^{(r)}=\varnothing$,
because every value of $\operatorname{core}_r$ is an $r$-core.  When
$|\lambda|=d$ is fixed,
\eqref{eq:core-weight-size} shows that the core determines the weight.

We recall \cite[Corollary~2.8, p.~586]{Walker}.

\begin{proposition}\label{prop:block-preservation}
For every $r\geq2$ and every partition $\lambda$, the Schur expansion
of $\mathfrak{s}_\lambda$ involves only partitions with the same $r$-core as
$\lambda$.
\end{proposition}

By Proposition~\ref{prop:block-preservation}, $\Theta_r$ preserves
the Schur span indexed by each $\mathcal B_{\kappa,w}^{(r)}$.
Equivalently, its matrix in each degree is block diagonal when the
Schur basis is ordered by $r$-core.

\begin{example}[The blocks for $r=3$ and $|\lambda|=5$]
\label{ex:block-preservation-p3-d5}
There are seven partitions of $5$.  Packing their $3$-abaci gives the
following complete list of nonempty blocks:
\[
\begin{array}{c@{\qquad}c@{\qquad}l}
\text{$3$-core }\kappa&\text{weight }w&
   \mathcal B_{\kappa,w}^{(3)}\\ \hline
(2)&1&\{(5),(2,2,1),(2,1,1,1)\},\\
(1,1)&1&\{(4,1),(3,2),(1,1,1,1,1)\},\\
(3,1,1)&0&\{(3,1,1)\}.
\end{array}
\]
Order both the rows and columns by listing the blocks from top to
bottom and the partitions in each block from left to right.  Put
\[
 K=\begin{pmatrix}
  0&1&-1\\
  1&0&1\\
  -1&1&0
 \end{pmatrix}.
\]
Substituting \eqref{eq:truncated-complete-generating} into
\eqref{eq:refaudit-modular-jacobi-trudi} gives
\begin{equation}\label{eq:p3-degree-five-block-matrix}
 \left(\inner{\mathfrak s_\lambda}{s_\mu}\right)_{
    \lambda,\mu\vdash5}
 =\operatorname{diag}(K,K,(1)).
\end{equation}
The two weight-one blocks have the same matrix $K$, while the zero
off-diagonal blocks agree with Proposition~\ref{prop:block-preservation}.
\end{example}

\subsection{Heisenberg modes and the operators
\texorpdfstring{$\mathcal D_{r,u}$}{D(r,u)}}
\label{subsec:power-sum-scaling}

We use the conventions for cores, quotients, weights, and the $r$-sign
from Subsection~\ref{subsec:abacus-core-quotient}.
For $m\geq1$, recall from
Subsection~\ref{subsec:ordinary-symmetric-functions} that $p_m$ denotes
the power sum function. For an integer $r\geq2$, put
\[
 a_{-m}^{(r)}=p_{rm},
 \qquad
 a_m^{(r)}=p_{rm}^{\perp}
 \qquad(m\geq1).
\]
For operators $A,B$, write $[A,B]=AB-BA$.
Since $p_n^{\perp}=n\,\partial/\partial p_n$, these operators satisfy
\begin{equation}\label{eq:heisenberg-relations}
 [a_m^{(r)},a_{-n}^{(r)}]=rm\,\delta_{mn},
 \qquad
 [a_m^{(r)},a_n^{(r)}]=[a_{-m}^{(r)},a_{-n}^{(r)}]=0
 \qquad(m,n\geq1).
\end{equation}
Thus they give a Heisenberg action of central charge $r$ on
$\Lambda_{\mathbb Q}$.  Define the degree-zero operator
\begin{equation}\label{eq:divisible-power-sum-counting-operator}
\mathcal N_r
=\sum_{m\geq1}\frac{1}{rm}a_{-m}^{(r)}a_m^{(r)}
=\sum_{m\geq1}p_{rm}\frac{\partial}{\partial p_{rm}}.
\end{equation}
The sum is finite on every homogeneous component.  If
$\eta=(1^{m_1}2^{m_2}\cdots)$ is a partition, then
\begin{equation}\label{eq:number-operator-eigenvalue}
\mathcal N_r(p_\eta)=\left(\sum_{r\mid j}m_j\right)p_\eta.
\end{equation}
For $u\in\mathbb Q$, define $\mathcal D_{r,u}$ by
\begin{equation}\label{eq:power-sum-scaling-operator}
\mathcal D_{r,u}(p_\eta)=u^{\sum_{r\mid j}m_j}p_\eta,
\end{equation}
where $u^0=1$.  Since the exponent in
\eqref{eq:power-sum-scaling-operator} is additive under multiplication
of power sums, \eqref{eq:power-sum-scaling-operator} is multiplicative.
As the power sums freely generate
$\Lambda_{\mathbb Q}$ as a $\mathbb Q$-algebra,
$\mathcal D_{r,u}$ is the algebra endomorphism determined by
\begin{equation}\label{eq:power-sum-scaling-generators}
\mathcal D_{r,u}(p_m)=
\begin{cases}
u p_m,&r\mid m,\\
p_m,&r\nmid m.
\end{cases}
\end{equation}
Thus $\mathcal D_{r,u}\mathcal D_{r,v}=\mathcal D_{r,uv}$, and
$\mathcal D_{r,u}$ is invertible when $u\neq0$.

Comparing \eqref{eq:theta-power-sums} with
\eqref{eq:power-sum-scaling-generators} gives
$\mathcal D_{r,1-r}=\Theta_r$, so
$\mathcal D_{r,1-r}(s_\lambda)=\mathfrak s_\lambda$.
We study the family $\mathcal D_{r,u}$ to compute this action
uniformly in core--quotient coordinates.

For $d\geq0$, the eigenvalues in
\eqref{eq:number-operator-eigenvalue} belong to
$\{0,1,\ldots,\lfloor d/r\rfloor\}$.  Choose an interpolation
polynomial $F_{d,u}\in\mathbb Q[t]$ satisfying $F_{d,u}(k)=u^k$ on
this set.  Then
\begin{equation}\label{eq:d-polynomial-in-n}
\mathcal D_{r,u}=F_{d,u}(\mathcal N_r)
\qquad\text{on }\Lambda_{\mathbb Q}^d.
\end{equation}

\begin{proposition}
\label{prop:d-core-span-preservation}
For every integer $r\geq2$ and every $u\in\mathbb Q$, the operator
$\mathcal D_{r,u}$ preserves the Schur span of each fixed $r$-core in
every degree.
\end{proposition}

\begin{proof}
By the Murnaghan--Nakayama rule, for $m\geq1$,
$p_{rm}^{\perp}$ removes ribbons of
length $rm$ and multiplication by $p_{rm}$ adds ribbons of that
length.  On an $r$-abacus, either operation moves a bead within one
runner and therefore leaves the $r$-core unchanged.  Every summand of
$\mathcal N_r$ in \eqref{eq:divisible-power-sum-counting-operator}
consequently preserves the span of the Schur functions
with any prescribed $r$-core.
For each $d$, relation~\eqref{eq:d-polynomial-in-n}
expresses $\mathcal D_{r,u}$ as a polynomial in $\mathcal N_r$ on
$\Lambda_{\mathbb Q}^d$.  Hence $\mathcal D_{r,u}$ preserves the same
fixed-core Schur spans in degree $d$, and therefore in every degree.
\end{proof}

\subsection{Core--quotient coordinates and half-vertex operators}
\label{subsec:heisenberg-intertwiner}

Lam used a map of Heisenberg modules from Fock space to symmetric
functions to study ribbon tableau generating functions
\cite[Theorem~9]{LamRibbon}.  Here we identify each fixed-core Schur
space with a tensor product of $r$ symmetric-function spaces.  We use the quotient
components themselves for half-vertex calculations and their
conjugates for the modular Schur coefficients.

For integers $r\geq2$ and $w\geq0$, let $\mathsf P_r(w)$ be the set
of $r$-multipartitions of total size $w$:
\[
 \mathsf P_r(w)
 =\left\{
 \boldsymbol\alpha=(\alpha^{(0)},\ldots,\alpha^{(r-1)})
 \in\mathscr P^r:
 \sum_{i=0}^{r-1}|\alpha^{(i)}|=w
 \right\}.
\]
For an $r$-tuple of alphabets
$\mathbf X=(X^{(0)},\ldots,X^{(r-1)})$, write
\[
 s_{\boldsymbol\alpha}(\mathbf X)
 =\prod_{i=0}^{r-1}s_{\alpha^{(i)}}(X^{(i)}),
 \qquad
 p_m^{(i)}=p_m(X^{(i)}),
 \qquad
 P_m=\sum_{i=0}^{r-1}p_m^{(i)}.
\]
Define the product Hall inner product
$\langle\cdot,\cdot\rangle_{\otimes}$ on
$\Lambda_{\mathbb Q}^{\otimes r}$ by declaring the product Schur
functions to be an orthonormal basis.
Write $P_m^{\perp}$ for the adjoint of multiplication by $P_m$
with respect to $\langle\cdot,\cdot\rangle_{\otimes}$.
Thus $P_m^{\perp}=m\sum_{i=0}^{r-1}\partial/\partial p_m^{(i)}$.
We write $\Lambda_{\mathbb Q}^{\otimes r}[w]$ for its homogeneous
component of total degree $w$; it is zero for $w<0$.

Fix an $r$-core $\delta$.  For $w\geq0$, recall that
$\mathcal B_{\delta,w}^{(r)}$ is the $r$-block indexed by $(\delta,w)$, and put
\[
\mathcal F_{\delta,w}
=\operatorname{span}_{\mathbb Q}
  \{s_\lambda:\lambda\in\mathcal B_{\delta,w}^{(r)}\},
\qquad
\mathcal F_\delta=\bigoplus_{w\geq0}\mathcal F_{\delta,w}.
\]
For $\lambda\in\mathcal B_{\delta,w}^{(r)}$, put
$\boldsymbol\alpha(\lambda)=((\lambda^{(0)})',\ldots,
(\lambda^{(r-1)})')$.  Define
\begin{align}
 \Phi_{\delta,w}&:\mathcal F_{\delta,w}
 \longrightarrow\Lambda_{\mathbb Q}^{\otimes r}[w],
 \label{eq:signed-quotient-coordinate-map}\\*
 \Phi_{\delta,w}(s_\lambda)
 &=\eps_r(\lambda)s_{\boldsymbol\alpha(\lambda)},
 \label{eq:signed-quotient-coordinate-basis}
\end{align}
and let $\Phi_\delta=\bigoplus_{w\geq0}\Phi_{\delta,w}$.
We grade $\mathcal F_\delta$ by $r$-weight $w$; its original
symmetric-function degree is $|\delta|+rw$.
Proposition~\ref{prop:core-quotient-correspondence} identifies
$\mathcal B_{\delta,w}^{(r)}$ with $\mathsf P_r(w)$.  Componentwise conjugation is a size-preserving bijection
and is realized by $\omega^{\otimes r}$, an isometry for the product
Hall inner product.  Multiplication by $\eps_r(\lambda)\in\{\pm1\}$
also preserves inner products.  Hence
$\Phi_\delta$ is a graded isometric isomorphism.

We also use the isometric map without conjugating the quotient
components,
\begin{equation}\label{eq:untransposed-quotient-map}
 \Psi_\delta:\mathcal F_\delta\longrightarrow\Lambda_{\mathbb Q}^{\otimes r},
 \qquad
 \Psi_\delta(s_\lambda)
 =\eps_r(\lambda)\prod_{i=0}^{r-1}s_{\lambda^{(i)}}(X^{(i)}).
\end{equation}
Thus $\Phi_\delta=(\omega^{\otimes r})\Psi_\delta$.

\begin{theorem}[Heisenberg intertwining]
\label{thm:quotient-coordinate-action}
For every $m\geq1$, the following identities hold as linear operators
on $\Lambda_{\mathbb Q}^{\otimes r}$:
\begin{align}
 \Phi_\delta\,p_{rm}\,\Phi_\delta^{-1}
 &=(-1)^{m-1}P_m,
 \label{eq:heisenberg-creation-intertwining}\\
 \Phi_\delta\,p_{rm}^{\perp}\,\Phi_\delta^{-1}
 &=(-1)^{m-1}P_m^{\perp}.
 \label{eq:heisenberg-annihilation-intertwining}
\end{align}
\end{theorem}

\begin{proof}
The adjointness identity and skew core--quotient formula in
\cite[Propositions~2.2--2.3, pp.~582--583]{Walker} give, for partitions
$\lambda,\mu$ with $r$-core $\delta$,
\begin{equation}\label{eq:refaudit-quotient-pairing}
\begin{aligned}
 \eps_r(\lambda)\eps_r(\mu)
 \langle p_{rm}s_\lambda,s_\mu\rangle
 &=\left\langle p_m,
      \prod_{i=0}^{r-1}s_{\mu^{(i)}/\lambda^{(i)}}\right\rangle\\
 &=\left\langle
      P_m\prod_{i=0}^{r-1}s_{\lambda^{(i)}}(X^{(i)}),
      \prod_{i=0}^{r-1}s_{\mu^{(i)}}(X^{(i)})
    \right\rangle_{\otimes}.
\end{aligned}
\end{equation}
The second equality in \eqref{eq:refaudit-quotient-pairing}
can be read directly from the derivation
$p_m^\perp=m\,\partial/\partial p_m$.
Pairing $p_m$ with the product of skew Schur functions sums over
the choice of one factor to pair with $p_m$, while every other
factor pairs with $1$.
Here $\langle1,s_{\mu^{(j)}/\lambda^{(j)}}\rangle
=\delta_{\lambda^{(j)},\mu^{(j)}}$, and skew Schur adjointness gives
$\langle p_m,s_{\mu^{(i)}/\lambda^{(i)}}\rangle
=\langle p_m s_{\lambda^{(i)}},s_{\mu^{(i)}}\rangle$.
Summing over $i$ is exactly the action of $P_m=\sum_i p_m^{(i)}$
in the product Hall pairing.
Thus the equality of matrix entries in
\eqref{eq:refaudit-quotient-pairing} gives
\begin{equation}\label{eq:untransposed-heisenberg-intertwining}
 \Psi_\delta\,p_{rm}\,\Psi_\delta^{-1}=P_m.
\end{equation}
Since $\Psi_\delta$ is an isometry, taking Hall adjoints in
\eqref{eq:untransposed-heisenberg-intertwining} gives
\begin{equation}\label{eq:untransposed-heisenberg-adjoint}
 \Psi_\delta p_{rm}^{\perp}\Psi_\delta^{-1}=P_m^{\perp}.
\end{equation}
Now $\Phi_\delta=(\omega^{\otimes r})\Psi_\delta$, and
$\omega(p_m)=(-1)^{m-1}p_m$.  This proves
\eqref{eq:heisenberg-creation-intertwining}.  Since $\Phi_\delta$ is
an isometry, taking Hall adjoints gives
\eqref{eq:heisenberg-annihilation-intertwining}.
\end{proof}

Multiplication by $P_m$ raises total degree by $m$, and
$P_m^{\perp}$ lowers total degree by $m$.  Hence, for every $w\geq0$,
\eqref{eq:heisenberg-creation-intertwining} maps degree $w$ to
degree $w+m$, while \eqref{eq:heisenberg-annihilation-intertwining}
maps degree $w$ to degree $w-m$ and is zero when $w<m$.

For an auxiliary alphabet $A$, write
\[
 \operatorname{Exp}[A]
 =\exp\left(\sum_{m\geq1}\frac{p_m[A]}m\right),
 \qquad
 D(A,X)=\exp\left(
   \sum_{m\geq1}\frac{(-1)^{m-1}}m p_m[A]p_m(X)\right).
\]
Define the creation half-vertex operator
\begin{equation}\label{eq:creation-half-vertex}
 \Gamma_{r,+}(A)
 =\exp\left(\sum_{m\geq1}\frac{p_m[A]}m p_{rm}\right),
\end{equation}
and let $\Gamma_{r,-}(A)=\Gamma_{r,+}(A)^\perp$, where the adjoint
acts only on the original symmetric-function variables.  On
$\Lambda_{\mathbb Q}(X)$, define the ordinary half-vertex operators
\begin{equation}\label{bridge:eq:half-fields}
 \mathcal H_+(A)
 =\exp\left(\sum_{m\geq1}\frac{p_m[A]}m p_m\right),
 \qquad
 \mathcal H_-(A)
 =\exp\left(\sum_{m\geq1}p_m[A]
                  \frac{\partial}{\partial p_m}\right).
\end{equation}
Thus $\mathcal H_+(A)$ is multiplication by
$C(A,X)=\operatorname{Exp}[AX]$, and $\mathcal H_-(A)$ sends $f[X]$ to
$f[X+A]$.  The operators $\mathcal H_\pm(A)$ are the standard creation
and annihilation factors with the normalizations in
\eqref{bridge:eq:half-fields}; see \cite[Section~4]{DLT}
and \cite{JingVertex}.  The creation operators take values in
the degree completions, the annihilation operators are locally finite,
and $\Phi_\delta$ and $\Psi_\delta$ extend degreewise to the
completions.  A superscript $(i)$ indicates action on the $i$th tensor
factor, and a kernel in an operator identity denotes multiplication
by that kernel.
All adjoints on the tensor space below use
$\langle\cdot,\cdot\rangle_{\otimes}$ and leave auxiliary alphabets
unchanged.

For two auxiliary alphabets $A,B$, the commutator of the exponents
in $\mathcal H_-(B)\mathcal H_+(A)$ is the scalar
$\sum_{m\geq1}p_m[A]p_m[B]/m$.  Consequently,
\begin{equation}\label{bridge:eq:half-commutation}
 \mathcal H_-(B)\mathcal H_+(A)
 =C(A,B)\mathcal H_+(A)\mathcal H_-(B),
 \qquad C(A,B)=\operatorname{Exp}[AB].
\end{equation}

\begin{corollary}
\label{cor:half-vertex-factorization}
For every auxiliary alphabet $A$,
\begin{equation}\label{bridge:eq:half-fields-core-quotient}
 \Psi_\delta\Gamma_{r,\pm}(A)\Psi_\delta^{-1}
 =\prod_{i=0}^{r-1}\mathcal H_\pm^{(i)}(A).
\end{equation}
\end{corollary}

\begin{proof}
Conjugation by $\Psi_\delta$ preserves products of operators and hence
their formal exponentials.  By
\eqref{eq:untransposed-heisenberg-intertwining} and
\eqref{eq:creation-half-vertex},
\begin{equation}\label{eq:refaudit-creation-factorization-proof}
\begin{aligned}
 \Psi_\delta\Gamma_{r,+}(A)\Psi_\delta^{-1}
 &=\exp\left(\sum_{m\geq1}\frac{p_m[A]}m P_m\right)\\
 &=\exp\left(\sum_{i=0}^{r-1}\sum_{m\geq1}
                  \frac{p_m[A]}m p_m^{(i)}\right)\\
 &=\prod_{i=0}^{r-1}\exp\left(
                  \sum_{m\geq1}\frac{p_m[A]}m p_m^{(i)}\right)\\
 &=\prod_{i=0}^{r-1}\mathcal H_+^{(i)}(A).
\end{aligned}
\end{equation}
In \eqref{eq:refaudit-creation-factorization-proof}, the third equality
uses commutativity on distinct tensor factors, and the last uses
\eqref{bridge:eq:half-fields}.

Since $\Psi_\delta$ is an isometry, taking adjoints in
\eqref{eq:refaudit-creation-factorization-proof} gives
\begin{equation}\label{eq:refaudit-annihilation-factorization-proof}
\begin{aligned}
 \Psi_\delta\Gamma_{r,-}(A)\Psi_\delta^{-1}
 &=\left(\Psi_\delta\Gamma_{r,+}(A)\Psi_\delta^{-1}\right)^\perp\\
 &=\left(\prod_{i=0}^{r-1}\mathcal H_+^{(i)}(A)\right)^\perp
 =\prod_{i=0}^{r-1}\mathcal H_-^{(i)}(A).
\end{aligned}
\end{equation}
By \eqref{bridge:eq:half-fields} and
$p_m^\perp=m\,\partial/\partial p_m$, the adjoint of
$\mathcal H_+^{(i)}(A)$ is $\mathcal H_-^{(i)}(A)$.
The reversed order of the factors in
\eqref{eq:refaudit-annihilation-factorization-proof} gives the same
product because operators on distinct tensor factors commute.
\end{proof}

The algebra involution $\omega_X$ acts only on $X$, leaving $A$
unchanged.  Since $\omega_X C(A,X)=D(A,X)$, for
$f\in\Lambda_{\mathbb Q}(X)$ we have
\begin{equation}\label{eq:refaudit-omega-kernel-conjugation}
 (\omega_X\mathcal H_+(A)\omega_X)f
 =\omega_X\bigl(C(A,X)\omega_X f\bigr)
 =D(A,X)f.
\end{equation}
Applying
\eqref{eq:refaudit-omega-kernel-conjugation} to each tensor
factor in \eqref{bridge:eq:half-fields-core-quotient}, and using
$\Phi_\delta=\omega^{\otimes r}\Psi_\delta$, gives
\begin{equation}\label{eq:creation-half-vertex-factorization}
 \Phi_\delta\Gamma_{r,+}(A)\Phi_\delta^{-1}
 =\prod_{i=0}^{r-1}D(A,X^{(i)}).
\end{equation}
Since $\Phi_\delta$ is an isometry and
$\Gamma_{r,-}(A)=\Gamma_{r,+}(A)^\perp$, taking adjoints in
\eqref{eq:creation-half-vertex-factorization} yields
\begin{equation}\label{eq:annihilation-half-vertex-factorization}
\begin{aligned}
 \Phi_\delta\Gamma_{r,-}(A)\Phi_\delta^{-1}
 &=\left(\Phi_\delta\Gamma_{r,+}(A)\Phi_\delta^{-1}\right)^\perp\\
 &=\left(\prod_{i=0}^{r-1}D(A,X^{(i)})\right)^\perp.
\end{aligned}
\end{equation}

\subsection{The induced action and its reproducing kernel}
\label{subsec:general-action-spectra}

We first compute the action of $\mathcal D_{r,u}$ in quotient
coordinates, then its reproducing kernel.  The key is to separate
the sums $P_m$ from differences between the power sums of the
individual alphabets.

On $\Lambda_{\mathbb Q}^{\otimes r}$ define
\begin{equation}\label{eq:quotient-number-operator}
 \widehat{\mathcal N}_r
 =\sum_{m\geq1}\frac{1}{rm}P_mP_m^{\perp}.
\end{equation}
For $m\geq1$ and $0\leq a<r-1$, put
\[
 v_{m,a}=p_m^{(a)}-p_m^{(r-1)}.
\]
For all $m,n\geq1$ and $0\leq a<r-1$,
\[
 P_m^{\perp}(P_n)=rm\,\delta_{mn},
 \qquad
 P_m^{\perp}(v_{n,a})
 =m\delta_{mn}-m\delta_{mn}=0.
\]
For each $m$, the change from the generators
$p_m^{(0)},\ldots,p_m^{(r-1)}$ to
$P_m,v_{m,0},\ldots,v_{m,r-2}$ is invertible over $\mathbb Q$; indeed,
\[
 p_m^{(r-1)}
 =\frac1r\left(P_m-\sum_{a=0}^{r-2}v_{m,a}\right),
 \qquad
 p_m^{(a)}=v_{m,a}+p_m^{(r-1)}.
\]
Hence the elements $P_m$ for $m\geq1$, together with the elements
$v_{m,a}$ for $m\geq1$ and $0\leq a<r-1$, freely generate
$\Lambda_{\mathbb Q}^{\otimes r}$.  If
\[
 M=\prod_{j\geq1}P_j^{c_j}
   \prod_{j\geq1}\prod_{a=0}^{r-2}v_{j,a}^{d_{j,a}}
\]
is a monomial in these generators, with nonnegative integer exponents, all but finitely many zero, then the derivation property of $P_m^{\perp}$ gives
\[
 \frac{1}{rm}P_mP_m^{\perp}(M)=c_mM,
 \qquad
 \widehat{\mathcal N}_r(M)=\left(\sum_{m\geq1}c_m\right)M.
\]
Thus $\widehat{\mathcal N}_r$ counts the total number of factors
$P_m$, with multiplicity, and does not count any of the factors
$v_{m,a}$.

For $u\in\mathbb Q$, let $\widehat{\mathcal D}_{r,u}$ be the algebra
endomorphism of $\Lambda_{\mathbb Q}^{\otimes r}$ defined by
\begin{equation}\label{eq:quotient-d-action}
 \widehat{\mathcal D}_{r,u}(P_m)=uP_m,
 \qquad
 \widehat{\mathcal D}_{r,u}(v_{m,a})=v_{m,a}.
\end{equation}
We use $u^0=1$, including when $u=0$.
For $M$, \eqref{eq:quotient-d-action} gives
\[
 \widehat{\mathcal D}_{r,u}(M)
 =u^{\sum_{m\geq1}c_m}M.
\]
Equivalently, for $m\geq1$ and $0\leq j<r$,
\begin{equation}\label{eq:quotient-d-generators}
 \widehat{\mathcal D}_{r,u}(p_m^{(j)})
 =p_m^{(j)}+\frac{u-1}{r}P_m.
\end{equation}

\begin{theorem}
\label{thm:d-under-phi}
For every $r\geq2$, every $r$-core $\delta$, and every
$u\in\mathbb Q$,
\begin{align}
 \Phi_\delta\mathcal N_r\Phi_\delta^{-1}
 &=\widehat{\mathcal N}_r,
 \label{eq:n-under-phi}\\
 \Phi_\delta\mathcal D_{r,u}\Phi_\delta^{-1}
 &=\widehat{\mathcal D}_{r,u}.
 \label{eq:d-under-phi}
\end{align}
\end{theorem}

\begin{proof}
By \eqref{eq:divisible-power-sum-counting-operator},
\[
 \mathcal N_r
 =\sum_{m\geq1}\frac{1}{rm}p_{rm}p_{rm}^{\perp}.
\]
The two signs in
\eqref{eq:heisenberg-creation-intertwining} and
\eqref{eq:heisenberg-annihilation-intertwining} cancel, proving
\eqref{eq:n-under-phi}.  On every finite homogeneous component,
$\mathcal D_{r,u}$ and $\widehat{\mathcal D}_{r,u}$ assign the
eigenvalue $u^k$ to the $k$th eigenspaces of $\mathcal N_r$ and
$\widehat{\mathcal N}_r$, respectively.
Since \eqref{eq:n-under-phi} identifies these eigenspaces under
$\Phi_\delta$, this proves
\eqref{eq:d-under-phi}.
\end{proof}

The operator on the right of \eqref{eq:d-under-phi} depends only on
$r,u$.  Thus its action in quotient coordinates is the same for
every $r$-core $\delta$.

For another $r$-tuple of alphabets
$\mathbf Y=(Y^{(0)},\ldots,Y^{(r-1)})$, let
$P_m(\mathbf X)=\sum_i p_m(X^{(i)})$ and define
$P_m(\mathbf Y)$ analogously.  For a degree-preserving linear operator
$T$ on $\Lambda_{\mathbb Q}^{\otimes r}$, a series
$K_T(\mathbf X,\mathbf Y)$ is called its \emph{reproducing kernel}
for the product Hall inner product if
\begin{equation}\label{eq:operator-reproducing-kernel}
 (Tf)(\mathbf X)
 =\left\langle
   K_T(\mathbf X,\mathbf Y),f(\mathbf Y)
  \right\rangle_{\otimes,\mathbf Y}
 \qquad
 \bigl(f\in\Lambda_{\mathbb Q}^{\otimes r}\bigr),
\end{equation}
where $\langle\cdot,\cdot\rangle_{\otimes,\mathbf Y}$ pairs only
the $\mathbf Y$ variables.  The product Cauchy kernel
\begin{equation}\label{eq:refaudit-product-cauchy-kernel}
 \prod_{j=0}^{r-1}C(X^{(j)},Y^{(j)})
 =\sum_{\boldsymbol\alpha\in\mathscr P^r}
   s_{\boldsymbol\alpha}(\mathbf X)
   s_{\boldsymbol\alpha}(\mathbf Y)
\end{equation}
is the reproducing kernel of the identity operator by
\eqref{eq:refaudit-ordinary-cauchy} and
\eqref{eq:operator-reproducing-kernel}.  Define
\[
 \mathcal C_{r,u}(\mathbf X,\mathbf Y)
 :=\widehat{\mathcal D}_{r,u}\left(
     \prod_{j=0}^{r-1}C(X^{(j)},Y^{(j)})\right),
\]
where $\widehat{\mathcal D}_{r,u}$ acts only on the $\mathbf X$
variables and leaves $\mathbf Y$ unchanged.  Since
$\widehat{\mathcal D}_{r,u}$ is a degree-preserving algebra endomorphism,
applying \eqref{eq:quotient-d-generators} to the exponential form of
$\prod_{j=0}^{r-1}C(X^{(j)},Y^{(j)})$ gives
\begin{equation}\label{eq:general-d-kernel}
\begin{aligned}
 \mathcal C_{r,u}(\mathbf X,\mathbf Y)
 &=\exp\left(
   \sum_{m\geq1}\frac1m\sum_{j=0}^{r-1}
   \left[p_m(X^{(j)})+\frac{u-1}{r}P_m(\mathbf X)\right]
   p_m(Y^{(j)})\right)\\
 &=\exp\left(
  \sum_{m\geq1}\frac1m
  \left[
   \sum_{j=0}^{r-1}p_m(X^{(j)})p_m(Y^{(j)})
   +\frac{u-1}{r}P_m(\mathbf X)P_m(\mathbf Y)
  \right]
 \right).
\end{aligned}
\end{equation}
The second equality in \eqref{eq:general-d-kernel} uses
$\sum_{j=0}^{r-1}p_m(Y^{(j)})=P_m(\mathbf Y)$.
The series $\mathcal C_{r,u}(\mathbf X,\mathbf Y)$ lies in the
graded completion of the tensor product in the $\mathbf X$ and
$\mathbf Y$ variables.  For every
$f\in\Lambda_{\mathbb Q}^{\otimes r}$, the action of
$\widehat{\mathcal D}_{r,u}$ on $\mathbf X$ commutes with the Hall
pairing in $\mathbf Y$.  Equations
\eqref{eq:refaudit-product-cauchy-kernel} and
\eqref{eq:operator-reproducing-kernel} therefore give
\begin{align*}
 &\left\langle
   \mathcal C_{r,u}(\mathbf X,\mathbf Y),f(\mathbf Y)
  \right\rangle_{\otimes,\mathbf Y}\\
 &\qquad=\widehat{\mathcal D}_{r,u}\left(
  \left\langle\prod_{j=0}^{r-1}C(X^{(j)},Y^{(j)}),f(\mathbf Y)
  \right\rangle_{\otimes,\mathbf Y}\right)\\
 &\qquad=(\widehat{\mathcal D}_{r,u}f)(\mathbf X).
\end{align*}
Thus $\mathcal C_{r,u}$ is the reproducing kernel of
$\widehat{\mathcal D}_{r,u}$ for the product Hall inner product.

\section{Schur-positive kernels and modular Schur coefficients}
\label{sec:modular-schur-coefficients}

We use the reproducing kernel \eqref{eq:general-d-kernel} to study
the signs and sizes of modular Schur coefficients.  We first
determine which parameter values give Schur-positive kernels after
composition with the grading involution.  We then specialize to
modular Schur functions and identify their normalized Schur
coefficients with coefficients of the specialized kernel.  The dual
Cauchy identity and the Littlewood--Richardson rule give a
nonnegative expansion, while a comparison with representations of
symmetric groups gives the factorial bound in
Theorem~\ref{thm:coefficient-conjecture}.  Together with the sharper
bound \eqref{eq:refined-coefficient-bound} and the equality criterion
\eqref{eq:factorial-equality-case}, these results prove
Theorem~\ref{thm:intro-universal-matrices}.  We also obtain the
support condition \eqref{eq:coefficient-support-condition} and the
nonsingularity of specified principal submatrices in
Corollary~\ref{cor:rectangular-principal-matrices}.

\subsection{Schur positivity of the reproducing kernels}
\label{subsec:universal-matrices}
\label{subsec:integral-specializations}

Fix an integer $r\geq2$ and two $r$-tuples of alphabets
$\mathbf X=(X^{(0)},\ldots,X^{(r-1)})$ and
$\mathbf Y=(Y^{(0)},\ldots,Y^{(r-1)})$.
Let $\sigma$ be the grading involution on
$\Lambda_{\mathbb Q}^{\otimes r}$,
\begin{equation}\label{eq:grading-involution}
 \sigma(f)=(-1)^{\deg f}f
 \qquad\text{for homogeneous }f.
\end{equation}
In total quotient degree $w\geq0$, this supplies the factor
$(-1)^w$ in the normalization \eqref{eq:def-a-lambda-mu} below.
Fix $u\in\mathbb Q$ and put
\[
 \tau=\frac{1-u}{r}.
\]
Extend $\sigma$ degreewise to the completion and let it act only on
the $\mathbf X$ variables. Define
\[
 \mathcal K_{r,\tau}(\mathbf X,\mathbf Y)
 :=\sigma\!\left(\mathcal C_{r,u}(\mathbf X,\mathbf Y)\right).
\]
By the reproducing identity \eqref{eq:operator-reproducing-kernel},
$\mathcal K_{r,\tau}$ is the reproducing kernel of
$\sigma\widehat{\mathcal D}_{r,u}$.
For $0\leq i,j<r$, abbreviate
\[
 D_{ij}=D(X^{(i)},Y^{(j)}).
\]
If $a\in\mathbb Q$, its rational power is interpreted in the graded
completion by
\[
 D_{ij}^{a}
 =\exp\left(
   a\sum_{m\geq1}\frac{(-1)^{m-1}}m
   p_m(X^{(i)})p_m(Y^{(j)})
  \right).
\]
Since $(u-1)/r=-\tau$, the kernel $\mathcal C_{r,u}$ in
\eqref{eq:general-d-kernel} factorizes as
\[
 \mathcal C_{r,u}(\mathbf X,\mathbf Y)
 =\prod_i C(X^{(i)},Y^{(i)})^{1-\tau}
  \prod_{i\neq j}C(X^{(i)},Y^{(j)})^{-\tau}.
\]
For $m\geq1$ and $0\leq i<r$, \eqref{eq:grading-involution} gives
$\sigma(p_m(X^{(i)}))=(-1)^m p_m(X^{(i)})$.
Since $\sigma$ is an algebra automorphism, for $a\in\mathbb Q$ and
$0\leq i,j<r$ we have
\[
\begin{aligned}
 \sigma\!\left(C(X^{(i)},Y^{(j)})^a\right)
 &=\exp\left(
    a\sum_{m\geq1}\frac{(-1)^m}{m}
    p_m(X^{(i)})p_m(Y^{(j)})\right)\\
 &=D_{ij}^{-a}.
\end{aligned}
\]
Thus the exponents $1-\tau$ and $-\tau$ become $\tau-1$ and
$\tau$, respectively, and applying $\sigma$ to
\eqref{eq:general-d-kernel} gives
\begin{equation}\label{eq:kernel-rt}
 \mathcal K_{r,\tau}(\mathbf X,\mathbf Y)
 =\prod_{i\neq j}D_{ij}^{\tau}
  \prod_iD_{ii}^{\tau-1}.
\end{equation}

\begin{theorem}[Schur-positive specializations]
\label{thm:schur-positive-specializations}
Let $r\geq2$ and $\tau\in\mathbb Q$.  All coefficients of
$\mathcal K_{r,\tau}$ in the product Schur basis are nonnegative if and
only if $\tau\in\mathbb Z_{>0}$.
\end{theorem}

\begin{proof}
If $\tau=Q\in\mathbb Z_{>0}$, then the exponents $Q$ and $Q-1$ in
\eqref{eq:kernel-rt} are nonnegative integers. By
\eqref{eq:refaudit-dual-cauchy}, each factor $D_{ij}$ has nonnegative
integral Schur coefficients, and multiplication preserves this
property by the Littlewood--Richardson rule. Hence all
product Schur coefficients of $\mathcal K_{r,\tau}$ are nonnegative
integers.

For the converse, we exhibit a negative coefficient whenever
$\tau\notin\mathbb Z_{>0}$.  If $\tau<0$, then for $i\neq j$ the coefficient of
$s_{(1)}(X^{(i)})s_{(1)}(Y^{(j)})$ is $\tau<0$.  If $0\leq \tau<1$, the
coefficient of $s_{(1)}(X^{(i)})s_{(1)}(Y^{(i)})$ is $\tau-1<0$.

It remains to consider $\tau\geq1$ with $\tau\notin\mathbb Z$.  Put
$n=\lfloor \tau\rfloor$ and $k=n+2$, and fix $i\neq j$.  Consider the
coefficient of $s_{(k)}(X^{(i)})s_{(k)}(Y^{(j)})$, with the empty
partition in every other component.  Since those other alphabets
have degree zero, every factor in \eqref{eq:kernel-rt} except
$D_{ij}^{\tau}$ must contribute its constant term.  Specialize
$X^{(i)}$ and $Y^{(j)}$ to single variables $x$ and $y$.
Then $s_{(k)}(x)=x^k$, while every other Schur function of degree
$k$ vanishes, and similarly for $y$.  Hence the coefficient is
\[
 \Coef_{x^ky^k}\bigl[(1+xy)^\tau\bigr]=\binom{\tau}{k}.
\]
Since
\[
 \binom{\tau}{n+2}
 =\frac{\tau(\tau-1)\cdots(\tau-n)(\tau-n-1)}{(n+2)!}<0,
\]
$\mathcal K_{r,\tau}$ is not Schur-positive.
\end{proof}

Since $\tau=(1-u)/r$, Theorem~\ref{thm:schur-positive-specializations}
shows that $\mathcal K_{r,\tau}$ is Schur-positive exactly when
\[
 u=1-rQ
 \qquad\text{for some }Q\in\mathbb Z_{>0}.
\]
For the remainder of this subsection, $Q$ is a positive integer.  For
$X=\sum_iX^{(i)}$ and $Y=\sum_iY^{(i)}$, set
\begin{equation}\label{eq:kernel-rq}
\begin{split}
 \mathcal K_{r,Q}(\mathbf X,\mathbf Y)
 &=D(X,Y)^Q\prod_iD(X^{(i)},Y^{(i)})^{-1}\\
 &=\prod_{i\neq j}D(X^{(i)},Y^{(j)})^Q
   \prod_iD(X^{(i)},Y^{(i)})^{Q-1}.
\end{split}
\end{equation}
The operator $\sigma\widehat{\mathcal D}_{r,1-rQ}$ and its reproducing
kernel $\mathcal K_{r,Q}$ are independent of the choice of the
$r$-core $\delta$.

\subsection{Walker's coefficient conjecture}
\label{subsec:coefficient-conjecture-proof}
\label{subsec:positive-expansions-bounds}

For the rest of this section, specialize $\mathcal K_{r,Q}$ to $Q=1$.
The parameter $u=1-rQ$ becomes $1-r$, and the diagonal factors in
\eqref{eq:kernel-rq} disappear.  We now identify the coefficients of
this kernel with the normalized modular Schur coefficients in each
fixed-core space, and then use the kernel to prove nonnegativity
and bounds.

Recall from \eqref{eq:theta-power-sums} and
\eqref{eq:power-sum-scaling-generators} that
\begin{equation}\label{eq:theta-as-d}
 \Theta_r=\mathcal D_{r,1-r}
 \qquad\text{on }\Lambda_{\mathbb Q}
\end{equation}
for every $r\geq2$.  Put
\begin{equation}\label{eq:specialized-kernel-operator}
 \mathcal T_r=\sigma\widehat{\mathcal D}_{r,1-r},
\end{equation}
and
\[
 \mathcal K_r(\mathbf X,\mathbf Y)
 :=\mathcal K_{r,1}(\mathbf X,\mathbf Y)
 =\prod_{i\neq j}D(X^{(i)},Y^{(j)}).
\]
Fix an integer $w\geq0$.  For
$\boldsymbol\alpha,\boldsymbol\beta\in\mathsf P_r(w)$, write
\begin{equation}\label{eq:abstract-kappa}
 \kappa^{(r,w)}_{\boldsymbol\alpha\boldsymbol\beta}
 =\Coef_{s_{\boldsymbol\alpha}(\mathbf X)s_{\boldsymbol\beta}(\mathbf Y)}
   \bigl[\mathcal K_r(\mathbf X,\mathbf Y)\bigr],
 \qquad
 A_{r,w}=\bigl(\kappa^{(r,w)}_{\boldsymbol\alpha\boldsymbol\beta}\bigr).
\end{equation}
The reproducing kernel of $\mathcal T_r$ is $\mathcal K_r$, so $A_{r,w}$ is
its matrix on $\Lambda_{\mathbb Q}^{\otimes r}[w]$ in the product
Schur basis, with rows indexing outputs and columns indexing inputs.
Since
$\mathcal K_r(\mathbf X,\mathbf Y)=\mathcal K_r(\mathbf Y,\mathbf X)$,
the matrix $A_{r,w}$ is symmetric.

We now identify $A_{r,w}$ with Walker's coefficient matrix.  Suppose
that $\lambda$ and $\mu$ lie in the same $r$-block of weight $w$.
Using the $r$-sign $\eps_r$, set
\begin{equation}\label{eq:def-a-lambda-mu}
 a_{\lambda\mu}
 =(-1)^w\eps_r(\lambda)\eps_r(\mu)
  \inner{\mathfrak{s}_\lambda}{s_\mu}.
\end{equation}

Fix $r\geq2$, an $r$-core $\delta$, and $w\geq0$.
By Proposition~\ref{prop:core-quotient-correspondence}, conjugating
each component of the $r$-quotient gives the bijection from
$\mathcal B_{\delta,w}^{(r)}$ to $\mathsf P_r(w)$,
\[
 \lambda\longmapsto
 \boldsymbol\alpha(\lambda)
 =\bigl((\lambda^{(0)})',\ldots,(\lambda^{(r-1)})'\bigr).
\]

\begin{proposition}
\label{prop:core-independent-matrix}
For all $\lambda,\mu\in\mathcal B_{\delta,w}^{(r)}$,
\begin{equation}\label{eq:core-independent-entry}
 a_{\lambda\mu}
 =\kappa^{(r,w)}_{
    \boldsymbol\alpha(\lambda),\boldsymbol\alpha(\mu)}.
\end{equation}
\end{proposition}

\begin{proof}
On $\Lambda_{\mathbb Q}^{\otimes r}[w]$,
\eqref{eq:grading-involution} says that $\sigma$ acts by $(-1)^w$.
Thus \eqref{eq:d-under-phi}, \eqref{eq:theta-as-d}, and
\eqref{eq:specialized-kernel-operator} give
\[
 \Phi_{\delta,w}\Theta_r\Phi_{\delta,w}^{-1}
 =(-1)^w\mathcal T_r
 \qquad\text{on }\Lambda_{\mathbb Q}^{\otimes r}[w].
\]
Using \eqref{eq:signed-quotient-coordinate-basis}, compare the entry
with input $s_{\boldsymbol\alpha(\lambda)}$ and output
$s_{\boldsymbol\alpha(\mu)}$ to obtain
\[
 \eps_r(\mu)\eps_r(\lambda)
 \inner{\mathfrak{s}_\lambda}{s_\mu}
 =(-1)^w
 \kappa^{(r,w)}_{
  \boldsymbol\alpha(\mu),\boldsymbol\alpha(\lambda)}.
\]
Multiplying by $(-1)^w$ and using \eqref{eq:def-a-lambda-mu} gives
the transposed version of \eqref{eq:core-independent-entry}.
The matrix $A_{r,w}$ is symmetric, so the entries agree in the stated
order.
\end{proof}

By \eqref{eq:core-independent-entry}, indexing both rows and columns
of $(a_{\lambda\mu})_{\lambda,\mu\in\mathcal B_{\delta,w}^{(r)}}$
by $\lambda\mapsto\boldsymbol\alpha(\lambda)$ gives $A_{r,w}$.

The core-independence of the normalized coefficients is
\cite[Corollary~2.10, p.~586]{Walker};
\eqref{eq:core-independent-entry} identifies them with the
coefficients of $\mathcal K_r$.

We now relate the finite-variable bases to these coefficient
matrices.  Fix $n\geq1$ and $r\geq2$, and use
the rectangle $\mathcal R_{n,r}$ from
Subsection~\ref{subsec:finite-homogeneous-substitution}.  Define
\[
 M_{\mathcal R}
 =\bigl(\langle\mathfrak s_\lambda,s_\mu\rangle\bigr)_{
                         \lambda,\mu\in\mathcal R_{n,r}}.
\]
For an $r$-core $\delta$ and $w\geq0$ with
$\mathcal R_{n,r}\cap\mathcal B_{\delta,w}^{(r)}\neq\varnothing$,
set
\[
 J_{\delta,w}^{n,r}
 =\{\boldsymbol\alpha(\lambda):
             \lambda\in\mathcal R_{n,r}\cap\mathcal B_{\delta,w}^{(r)}\}.
\]
Let $M_{\delta,w}^{\mathcal R}$ be the submatrix of $M_{\mathcal R}$
indexed by this intersection, and put
$D=\operatorname{diag}(\eps_r(\lambda))$ in the same order.

\begin{corollary}\label{cor:rectangular-principal-matrices}
For every $n\geq1$ and $r\geq2$, the matrix $M_{\mathcal R}$ is
nonsingular over $\mathbb Q$.  For every $r$-core $\delta$ and
$w\geq0$ with $J_{\delta,w}^{n,r}\neq\varnothing$, the principal
submatrix
$A_{r,w}[J_{\delta,w}^{n,r},J_{\delta,w}^{n,r}]$ is nonsingular.
\end{corollary}

\begin{proof}
Theorem~\ref{thm:finite-variable-bases} and
\eqref{eq:finite-space-basis} give two bases of $V_{n,n(r-1)}$,
indexed respectively by $\mathfrak s_\lambda$ and $s_\mu$.
On specializing the stable Schur expansion to $n$ variables, terms
with $\ell(\mu)>n$ vanish, whereas coefficients with
$\ell(\mu)\leq n$ remain unchanged, since those Schur polynomials
are linearly independent.  Membership in $V_{n,n(r-1)}$ and its Schur
basis in \eqref{eq:finite-space-basis} then force all surviving coefficients
with $\mu_1>n(r-1)$ to vanish.  Hence $M_{\mathcal R}$ is exactly the
change-of-basis matrix between the two finite-variable bases, and
is nonsingular.

Homogeneity and Proposition~\ref{prop:block-preservation} make
$M_{\mathcal R}$ block diagonal by $r$-core and weight.
Every nonempty block $M_{\delta,w}^{\mathcal R}$ is therefore
nonsingular.  Using the same order on the quotient labels,
Proposition~\ref{prop:core-independent-matrix} and
\eqref{eq:def-a-lambda-mu} give
\begin{equation}\label{eq:rectangular-block-conjugation}
 A_{r,w}[J_{\delta,w}^{n,r},J_{\delta,w}^{n,r}]
       =(-1)^wD M_{\delta,w}^{\mathcal R}D.
\end{equation}
Since $D^2=I$, \eqref{eq:rectangular-block-conjugation} implies that
$A_{r,w}[J_{\delta,w}^{n,r},J_{\delta,w}^{n,r}]$ is nonsingular.
\end{proof}

The index set $J_{\delta,w}^{n,r}$ depends on the core as well as
on $n,r,w$; the conclusion applies to these specified principal
submatrices.

We next record an expansion that will be used to compute a specific
Schur coefficient in Subsection~\ref{subsec:rouquier-nonvanishing}.
For a tuple of partitions $\boldsymbol\gamma$, write
$|\boldsymbol\gamma|=\sum_i|\gamma^{(i)}|$.
Proposition~\ref{prop:universal-skew-expansion} expresses the
coefficients in \eqref{eq:core-independent-entry} using skew Schur
functions and the quotient convention of
\eqref{eq:quotient-conjugation}.

\begin{proposition}
\label{prop:universal-kernel-expansions}
\label{prop:universal-skew-expansion}
For $\boldsymbol\alpha,\boldsymbol\beta\in\mathsf P_r(w)$, one has
the signed expansion
\begin{equation}\label{eq:universal-skew-expansion}
 \kappa^{(r,w)}_{\boldsymbol\alpha\boldsymbol\beta}
 =\sum_{\boldsymbol\gamma\in\mathscr P^r}
  (-1)^{|\boldsymbol\gamma|}
 \inner{
  \prod_i s_{\alpha^{(i)}/\gamma^{(i)}}
 }{
  \prod_i s_{(\beta^{(i)})'/(\gamma^{(i)})'}
 }.
\end{equation}
\end{proposition}

\begin{proof}
By Proposition~\ref{prop:core-quotient-correspondence}, choose
$\lambda,\mu\in\mathcal B_{\delta,w}^{(r)}$
with $\boldsymbol\alpha=\boldsymbol\alpha(\lambda)$ and
$\boldsymbol\beta=\boldsymbol\alpha(\mu)$.
By \cite[Proposition~2.9, p.~586]{Walker}, with the sum over all
partitions $\nu$ with $r$-core $\delta$,
\[
 \inner{\mathfrak{s}_\lambda}{s_\mu}
 =(-1)^w\eps_r(\lambda)\eps_r(\mu)
  \sum_\nu(-1)^{w(\nu)}
  \inner{
   \prod_{j=0}^{r-1}s_{(\lambda')^{(j)}/(\nu')^{(j)}}
  }{
   \prod_{i=0}^{r-1}s_{\mu^{(i)}/\nu^{(i)}}
  }.
\]
Set $\gamma^{(i)}=(\nu^{(i)})'$.
Proposition~\ref{prop:core-quotient-correspondence} identifies the
sum with one over $\boldsymbol\gamma\in\mathscr P^r$, and
$w(\nu)=|\boldsymbol\gamma|$.
By \eqref{eq:quotient-conjugation}, reindexing the first product by
$j=r-1-i$ gives $\prod_i s_{\alpha^{(i)}/\gamma^{(i)}}$,
while the second is
$\prod_i s_{(\beta^{(i)})'/(\gamma^{(i)})'}$.
Multiplying by $(-1)^w\eps_r(\lambda)\eps_r(\mu)$ and using
\eqref{eq:def-a-lambda-mu} and \eqref{eq:core-independent-entry} gives
\eqref{eq:universal-skew-expansion}.
\end{proof}

In \eqref{eq:universal-skew-expansion}, only tuples with
$\gamma^{(i)}\subseteq\alpha^{(i)}$ and
$\gamma^{(i)}\subseteq\beta^{(i)}$ for every $i$ contribute.
This follows from the vanishing convention for skew Schur functions.

To prove nonnegativity, we now expand $\mathcal K_r$ directly in the
product Schur basis using \eqref{eq:refaudit-dual-cauchy}.
For $m\geq1$ and partitions
$\rho^{(1)},\ldots,\rho^{(m)},\tau$, define the generalized
Littlewood--Richardson coefficient by
\[
 c^\tau_{\rho^{(1)},\ldots,\rho^{(m)}}
 :=\inner{\prod_{a=1}^m s_{\rho^{(a)}}}{s_\tau}.
\]
For $m=2$ this is the ordinary Littlewood--Richardson coefficient;
for general $m$, its nonnegativity follows by iterating the
Littlewood--Richardson rule.

Set
\[
 \mathcal I_r=\{(i,j):0\leq i,j<r,\ i\neq j\},
\]
and let $\mathscr P^{\mathcal I_r}$ be the set of arrays
$\boldsymbol\rho=(\rho^{(i,j)})_{(i,j)\in\mathcal I_r}$ of ordinary
partitions.  Thus $(\rho^{(i,j)})_{j\neq i}$ denotes the $r-1$
partitions with fixed first index, and
$((\rho^{(i,j)})')_{i\neq j}$ denotes the conjugates of the $r-1$
partitions with fixed second index.

The coefficient $\kappa^{(r,w)}_{\boldsymbol\alpha\boldsymbol\beta}$
also has the nonnegative expansion
\begin{equation}\label{eq:abstract-positive-lr}
 \kappa^{(r,w)}_{\boldsymbol\alpha\boldsymbol\beta}
 =\sum_{\boldsymbol\rho\in\mathscr P^{\mathcal I_r}}
  \prod_i c^{\alpha^{(i)}}_{(\rho^{(i,j)})_{j\neq i}}
  \prod_j c^{\beta^{(j)}}_{((\rho^{(i,j)})')_{i\neq j}}.
\end{equation}
Only arrays satisfying
\begin{equation}\label{eq:refaudit-lr-array-degrees}
 \sum_{j\neq i}|\rho^{(i,j)}|=|\alpha^{(i)}|,
 \qquad
 \sum_{i\neq j}|\rho^{(i,j)}|=|\beta^{(j)}|
\end{equation}
for every $i$ and $j$ can contribute.  In particular, every entry of
$A_{r,w}$ is a nonnegative integer.

To prove \eqref{eq:abstract-positive-lr}, apply
\eqref{eq:refaudit-dual-cauchy} to every factor of $\mathcal K_r$:
\[
 D(X^{(i)},Y^{(j)})
 =\sum_{\rho^{(i,j)}\in\mathscr P}
   s_{\rho^{(i,j)}}(X^{(i)})
   s_{(\rho^{(i,j)})'}(Y^{(j)})
 \qquad(i\ne j).
\]
For a fixed array $\boldsymbol\rho$, the factors involving
$X^{(i)}$ are $\prod_{j\ne i}s_{\rho^{(i,j)}}(X^{(i)})$;
their $s_{\alpha^{(i)}}$ coefficient is the corresponding
Littlewood--Richardson coefficient in the first product of
\eqref{eq:abstract-positive-lr}.  Similarly, the factors involving
$Y^{(j)}$ give the second product, with each partition conjugated.
Since the alphabets are independent, these coefficients multiply.
The conditions
\eqref{eq:refaudit-lr-array-degrees} bound the sizes of all the
partitions in an array, so only finitely many arrays contribute.

Equation~\eqref{eq:abstract-positive-lr} proves nonnegative
integrality.  We next bound the coefficients by dimensions of
induced representations of symmetric groups, obtaining a bound
that implies the required factorial bound.
For an ordinary partition $\rho$, let $f^\rho$ be the number of
standard Young tableaux of shape $\rho$, with $f^\varnothing=1$.
For $\boldsymbol\alpha\in\mathsf P_r(w)$, define
\begin{equation}\label{eq:universal-dimension}
 d_{\boldsymbol\alpha}
 =\frac{w!}{\prod_i|\alpha^{(i)}|!}\prod_i f^{\alpha^{(i)}}.
\end{equation}

\begin{theorem}
\label{thm:universal-refined-bound}
For all $\boldsymbol\alpha,\boldsymbol\beta\in\mathsf P_r(w)$,
\begin{equation}\label{eq:refined-coefficient-bound}
 0\leq
 \kappa^{(r,w)}_{\boldsymbol\alpha\boldsymbol\beta}
 \leq\min\{d_{\boldsymbol\alpha},d_{\boldsymbol\beta}\}
 \leq w!.
\end{equation}
\end{theorem}

\begin{proof}
We prove \eqref{eq:refined-coefficient-bound} by comparing
coefficients with \eqref{eq:refaudit-dual-cauchy}, interpreting the
resulting coefficients as dimensions of homomorphism spaces between
induced representations of symmetric groups, and using an embedding
into the regular representation.

For $X=\sum_iX^{(i)}$ and $Y=\sum_iY^{(i)}$, put
\[
 b_{\boldsymbol\alpha\boldsymbol\beta}
 =\Coef_{s_{\boldsymbol\alpha}(\mathbf X)s_{\boldsymbol\beta}(\mathbf Y)}
       [D(X,Y)].
\]
The factorization
\[
 D(X,Y)=\mathcal K_r(\mathbf X,\mathbf Y)
             \prod_iD(X^{(i)},Y^{(i)})
\]
has Schur-positive factors with constant term $1$.  Taking the
constant term from each diagonal factor recovers the coefficient
of $\mathcal K_r$; all remaining contributions are nonnegative by
the Littlewood--Richardson rule.  Hence
$0\leq\kappa^{(r,w)}_{\boldsymbol\alpha\boldsymbol\beta}
\leq b_{\boldsymbol\alpha\boldsymbol\beta}$.
Applying \eqref{eq:refaudit-dual-cauchy} to $D(X,Y)$ gives
\[
 b_{\boldsymbol\alpha\boldsymbol\beta}
 =\sum_{\tau\vdash w}
   c^\tau_{\alpha^{(0)},\ldots,\alpha^{(r-1)}}
   c^{\tau'}_{\beta^{(0)},\ldots,\beta^{(r-1)}}.
\]

For $m\geq0$, let $\Sym_m$ be the symmetric group on
$\{1,\ldots,m\}$, with $\Sym_0$ trivial.
For $\rho\vdash m$, let $S^\rho$ be the irreducible complex Specht
module of dimension $f^\rho$, with $S^\varnothing=\mathbb C$.
Write $\boxtimes$ for external tensor products and $\Ind$ for
induction, and write $\operatorname{Hom}_{\Sym_w}$ for equivariant
complex-linear maps. Set
\[
 H_{\boldsymbol\alpha}=\prod_i\Sym_{|\alpha^{(i)}|},
 \qquad
 U_{\boldsymbol\alpha}=\boxtimes_i S^{\alpha^{(i)}},
 \qquad
 I_{\boldsymbol\alpha}
 =\Ind_{H_{\boldsymbol\alpha}}^{\Sym_w}U_{\boldsymbol\alpha}.
\]
The Young subgroup acts on successive intervals of the indicated
sizes. Define the same objects for other tuples of total size $w$;
the subgroup index is $w!/\prod_i|\alpha^{(i)}|!$, so the dimension
formula for induction gives
$\dim I_{\boldsymbol\alpha}=d_{\boldsymbol\alpha}$ as in
\eqref{eq:universal-dimension}.
Put $\boldsymbol\beta'=((\beta^{(0)})',\ldots,(\beta^{(r-1)})')$.
The characteristic map preserves induction products and inner
products and sends irreducible characters to Schur functions
\cite[Chapter~I, (7.3), (7.5)--(7.6), pp.~113--114]{Macdonald}.
For $\tau\vdash w$, the multiplicities of $S^\tau$ in
$I_{\boldsymbol\alpha}$ and $I_{\boldsymbol\beta'}$ are
$c^\tau_{\alpha^{(0)},\ldots,\alpha^{(r-1)}}$ and
$c^{\tau'}_{\beta^{(0)},\ldots,\beta^{(r-1)}}$, respectively.
The multiplicity of $S^\tau$ in $I_{\boldsymbol\beta'}$ is obtained
by applying $\omega$ to $\prod_i s_{\beta^{(i)}}$.
Complete reducibility over $\mathbb C$ therefore gives
\[
 b_{\boldsymbol\alpha\boldsymbol\beta}
 =\dim\operatorname{Hom}_{\Sym_w}
      (I_{\boldsymbol\alpha},I_{\boldsymbol\beta'}).
\]

Every irreducible complex representation of a finite group occurs
in its regular representation.  Thus the irreducible
$H_{\boldsymbol\beta}$-module
$U_{\boldsymbol\beta'}$ embeds in the regular representation
$\mathbb C[H_{\boldsymbol\beta}]$.
Induction is exact over $\mathbb C$, so
\[
 I_{\boldsymbol\beta'}\hookrightarrow
 \Ind_{H_{\boldsymbol\beta}}^{\Sym_w}\mathbb C[H_{\boldsymbol\beta}]
 \cong\mathbb C[\Sym_w].
\]
Composing with the inclusion
$I_{\boldsymbol\beta'}\hookrightarrow\mathbb C[\Sym_w]$ gives an
injection on equivariant homomorphism spaces.  Consequently,
\begin{equation}\label{eq:refaudit-hom-regular-bound}
 b_{\boldsymbol\alpha\boldsymbol\beta}
 \leq\dim\operatorname{Hom}_{\Sym_w}
       (I_{\boldsymbol\alpha},\mathbb C[\Sym_w])
 =d_{\boldsymbol\alpha}.
\end{equation}
The equality in \eqref{eq:refaudit-hom-regular-bound} follows because
each irreducible representation occurs in the regular representation
with multiplicity equal to its dimension. Interchanging $\mathbf X$ and $\mathbf Y$ gives
$b_{\boldsymbol\alpha\boldsymbol\beta}\leq d_{\boldsymbol\beta}$.
Finally, $f^\rho\leq|\rho|!$ for every partition $\rho$, so
\eqref{eq:universal-dimension} gives
$d_{\boldsymbol\alpha},d_{\boldsymbol\beta}\leq w!$.
\end{proof}

The absence of diagonal factors in $\mathcal K_r$ gives further
information about the coefficients.  For the same
$\boldsymbol\alpha,\boldsymbol\beta\in\mathsf P_r(w)$, it imposes
\begin{equation}\label{eq:coefficient-support-condition}
 \kappa^{(r,w)}_{\boldsymbol\alpha\boldsymbol\beta}\neq0
 \quad\Longrightarrow\quad
 |\alpha^{(i)}|+|\beta^{(i)}|\leq w
 \qquad(0\leq i<r).
\end{equation}

Indeed, if $\kappa^{(r,w)}_{\boldsymbol\alpha\boldsymbol\beta}\neq0$,
choose an array $\boldsymbol\rho$ contributing a nonzero term to
\eqref{eq:abstract-positive-lr}.  Put
$m_{ij}=|\rho^{(i,j)}|$.  Then
\[
 |\alpha^{(i)}|=\sum_{j\neq i}m_{ij},
 \qquad
 |\beta^{(i)}|=\sum_{j\neq i}m_{ji},
 \qquad
 w=\sum_{k\neq l}m_{kl}.
\]
The entries in row $i$ and column $i$ are disjoint because the
diagonal entry is absent.  Their sum is therefore at most $w$, proving
\eqref{eq:coefficient-support-condition}.

For a tuple $\boldsymbol\alpha$, write
$\operatorname{supp}(\boldsymbol\alpha)=\{i:\alpha^{(i)}\ne\varnothing\}$.
Equality with $w!$ in \eqref{eq:refined-coefficient-bound} is
characterized by
\begin{equation}\label{eq:factorial-equality-case}
 \kappa^{(r,w)}_{\boldsymbol\alpha\boldsymbol\beta}=w!
 \quad\Longleftrightarrow\quad
 \left\{
 \begin{array}{l}
  \alpha^{(i)},\beta^{(i)}\in\{\varnothing,(1)\}
  \text{ for every }i,\\
  \operatorname{supp}(\boldsymbol\alpha)\cap
  \operatorname{supp}(\boldsymbol\beta)=\varnothing.
 \end{array}
 \right.
\end{equation}
Suppose that equality holds.  Equation
\eqref{eq:refined-coefficient-bound} forces
$d_{\boldsymbol\alpha}=d_{\boldsymbol\beta}=w!$, and
\eqref{eq:universal-dimension} gives
\[
 \frac{d_{\boldsymbol\alpha}}{w!}
 =\prod_i\frac{f^{\alpha^{(i)}}}{|\alpha^{(i)}|!}=1.
\]
For $|\rho|\geq2$, the standard tableaux form a proper subset of
all bijective fillings, so $f^\rho/|\rho|!<1$; this ratio is $1$
for $\rho=\varnothing$ or $(1)$.  Thus the product can equal $1$
only when every component of $\boldsymbol\alpha$ is empty or $(1)$,
and the same argument applies to $\boldsymbol\beta$.

When every component of $\boldsymbol\alpha$ and
$\boldsymbol\beta$ is $\varnothing$ or $(1)$, let
$S=\operatorname{supp}(\boldsymbol\alpha)$,
$T=\operatorname{supp}(\boldsymbol\beta)$.
Both supports have cardinality $w$.  By
\eqref{eq:refaudit-lr-array-degrees}, a contributing array has
exactly one nonempty entry in each row indexed by $S$ and each
column indexed by $T$, and every such entry is $(1)$.
The nonzero summands in \eqref{eq:abstract-positive-lr} therefore
correspond to bijections $\pi:S\to T$ satisfying $\pi(i)\ne i$
for $i\in S\cap T$, because diagonal entries are excluded.
Each contributing Littlewood--Richardson coefficient is $1$.
If $S\cap T=\varnothing$, all $w!$ bijections contribute.
If $S\cap T\ne\varnothing$, choose $i\in S\cap T$.
There is a bijection $S\to T$ fixing $i$, obtained by extending
$i\mapsto i$ with any bijection $S\setminus\{i\}\to T\setminus\{i\}$.
Thus strictly fewer than $w!$ bijections contribute. This proves both directions of
\eqref{eq:factorial-equality-case}.  When $w=0$, the unique empty
bijection contributes $1=0!$ and the support conditions are vacuous.

We now use Proposition~\ref{prop:core-independent-matrix} to return
to modular Schur coefficients.  Let $r\geq2$, and let
$\lambda,\mu$ lie in the same $r$-block of
weight $w$.  Put
\[
 \boldsymbol\alpha=\boldsymbol\alpha(\lambda)
 =\bigl((\lambda^{(0)})',\ldots,(\lambda^{(r-1)})'\bigr),
 \qquad
 \boldsymbol\beta=\boldsymbol\alpha(\mu)
 =\bigl((\mu^{(0)})',\ldots,(\mu^{(r-1)})'\bigr).
\]

Theorem~\ref{thm:coefficient-conjecture} proves the conjecture following \cite[Corollary~2.10, p.~586]{Walker}.

\begin{theorem}[Walker's coefficient conjecture]
\label{thm:coefficient-conjecture}
For $r\geq2$ and $\lambda,\mu$ in the same $r$-block of weight $w$,
\begin{equation}\label{eq:positive-lr}
 a_{\lambda\mu}
 =\sum_{\boldsymbol\rho\in\mathscr P^{\mathcal I_r}}
  \prod_i c^{\alpha^{(i)}}_{(\rho^{(i,j)})_{j\neq i}}
  \prod_j c^{\beta^{(j)}}_{((\rho^{(i,j)})')_{i\neq j}}.
\end{equation}
The coefficient $a_{\lambda\mu}$ is a nonnegative integer and
satisfies $a_{\lambda\mu}\leq w!$.
\end{theorem}

\begin{proof}
Proposition~\ref{prop:core-independent-matrix} gives
\[
 a_{\lambda\mu}
 =\kappa^{(r,w)}_{\boldsymbol\alpha\boldsymbol\beta}.
\]
Thus \eqref{eq:positive-lr} follows from
\eqref{eq:abstract-positive-lr}, whose summands are nonnegative
integers.  The factorial bound follows from
Theorem~\ref{thm:universal-refined-bound}.
\end{proof}

For the same $\lambda,\mu$, Proposition~\ref{prop:core-independent-matrix}
transfers the refined bound \eqref{eq:refined-coefficient-bound},
the necessary condition \eqref{eq:coefficient-support-condition},
and the equality criterion \eqref{eq:factorial-equality-case} to
$a_{\lambda\mu}$. Conjugation preserves component sizes and tableau
dimensions, and fixes $\varnothing$ and $(1)$.
In particular, $a_{\lambda\mu}=w!$ exactly when every quotient
component of $\lambda$ and $\mu$ is empty or $(1)$ and the sets of
their nonempty components are disjoint.
Both supports have cardinality $w$, so they can be disjoint among
the $r$ runners if and only if $r\geq2w$.
Proposition~\ref{prop:core-quotient-correspondence} realizes every
such ordered pair of supports in the block with any prescribed
$r$-core.  Thus, in every weight-$w$ block, the bound $w!$ is
attained if and only if $r\geq2w$.  When $r\geq2w$, choosing $w$
runners for the first support and then $w$ of the remaining runners
for the second gives the number of ordered pairs attaining the bound:
\begin{equation}\label{eq:number-factorial-equality-pairs}
 \binom rw\binom{r-w}{w}.
\end{equation}
This includes $w=0$, when the
unique coefficient is $1=0!$.

\section{Simple modular Schur characters}
\label{sec:conj46}

We prove that a modular Schur function is a simple polynomial
character exactly when its indexing partition satisfies the Carter
condition, recalled in Subsection~\ref{subsec:restricted-projection}.
For odd primes, Lemma~\ref{lem:initial-reduction} reduces any
counterexample not satisfying the Carter condition to a doubly singular
partition.  The operators
of Section~\ref{sec:universal-blocks} and categorical divided powers
allow us to transfer simplicity along
adjacent blocks to a Rouquier block.  There we determine the only
possible simple character and use the coefficient formula of
Proposition~\ref{prop:universal-skew-expansion} to show that the
remaining virtual character is nonzero.  Subsections
\ref{subsec:affine-intertwiners}--\ref{subsec:rouquier-nonvanishing}
establish these steps; Subsection~\ref{subsec:carter-classification}
assembles them and treats characteristic two separately.  The resulting
classification in Corollary~\ref{cor:carter-classification} proves
\cite[Conjecture~4.6]{Walker}.

We now set $r=p$, where $p$ is prime.  Throughout this section,
$\mathbb K$ is an arbitrary infinite field of characteristic $p$, and
$\mathfrak{s}_\lambda$ is taken with $r=p$.

Using the English convention for Young diagrams adopted in
Subsection~\ref{subsec:ordinary-symmetric-functions}, the nodes of
$\lambda$ are the pairs $(a,b)$ with $1\leq b\leq\lambda_a$, and
\begin{equation}\label{eq:hook-residue-conventions}
h_\lambda(a,b)=\lambda_a-b+\lambda'_b-a+1,
\qquad
\operatorname{res}(a,b)=\overline{b-a}\in\mathbb Z/p\mathbb Z
\end{equation}
are the hook length and residue of the node, where the bar denotes
the residue class modulo $p$.  We write
$v_p(m)$ for the $p$-adic valuation of a positive integer $m$.

\subsection{Heisenberg operators and the affine action}
\label{subsec:affine-intertwiners}

We first show that the operators $\mathcal D_{p,u}$ from
Subsection~\ref{subsec:power-sum-scaling} commute with adding and
removing nodes of a fixed residue.  For $\Theta_p=\mathcal D_{p,1-p}$,
this supplies the commutation hypothesis of
Theorem~\ref{thm:divided-power-simple-characters}, which we will apply
in Proposition~\ref{prop:adjacent-simplicity}.  On the
level-one Fock space $\Lambda_{\mathbb Q}$ with Schur basis, define,
for $i\in\mathbb Z/p\mathbb Z$,
\begin{equation}\label{eq:refaudit-chevalley-operators}
\begin{aligned}
e_i(s_\lambda)
&=\sum_{\substack{\nu\subset\lambda\\
\lambda/\nu\text{ is a removable }i\text{-node}}}s_\nu,\\
f_i(s_\lambda)
&=\sum_{\substack{\mu\supset\lambda\\
\mu/\lambda\text{ is an addable }i\text{-node}}}s_\mu,
\end{aligned}
\end{equation}
where an $i$-node is a box $(a,b)$ with $\overline{b-a}=i$.
A removable node of $\lambda$ is a box in its Young diagram whose
deletion leaves the Young diagram of a partition; an addable node is
a box outside the diagram whose addition gives the Young diagram of
a partition.  Over $\mathbb Q$, put
$e_i^{(c)}=e_i^c/c!$ and $f_i^{(c)}=f_i^c/c!$ for $c\geq0$.

To separate the node operations by residue, we recall the affine
weight of a partition.  For each $i\in\mathbb Z/p\mathbb Z$, let
$n_i(\lambda)$ be the number of $i$-nodes of $\lambda$.  Let
$\Lambda_0$ be the weight of the empty partition and let $\alpha_i$
be the affine simple roots of type $A_{p-1}^{(1)}$.  The affine
weight of $\lambda$, also assigned to $s_\lambda$, is
\begin{equation}\label{eq:fock-affine-weight}
\operatorname{wt}(\lambda)
=\Lambda_0-\sum_{i\in\mathbb Z/p\mathbb Z}n_i(\lambda)\alpha_i.
\end{equation}
For the argument below, the $\alpha_i$ can be treated as independent
symbols recording the residue counts.  Thus two partitions have
the same weight exactly when all these counts agree; in particular,
they have the same size, since $|\lambda|=\sum_i n_i(\lambda)$.

For a weight $\varpi$, the weight space
$\mathcal W_\varpi\subseteq\Lambda_{\mathbb Q}$ is the
$\mathbb Q$-span of the Schur functions $s_\lambda$ with
$\operatorname{wt}(\lambda)=\varpi$.  Write
$\operatorname{pr}_\varpi$ for the projection that retains precisely
these terms in a Schur expansion; it is zero if no partition has
weight $\varpi$.  Among partitions of a fixed size, equality of
weights is equivalent to equality of $p$-cores
\cite[Section~5.3, equation~(8)]{HTY}.
In Proposition~\ref{prop:d-chevalley-commutation}, we use these
projections to select the residue-$i$ terms after adding or removing
one node without restricting its residue.

\begin{proposition}
\label{prop:d-chevalley-commutation}
For every $u\in\mathbb Q$, every residue $i$, and every $c\geq0$,
\begin{equation}\label{eq:d-chevalley-commutation}
\mathcal D_{p,u}e_i^{(c)}=e_i^{(c)}\mathcal D_{p,u},
\qquad
\mathcal D_{p,u}f_i^{(c)}=f_i^{(c)}\mathcal D_{p,u}.
\end{equation}
\end{proposition}

\begin{proof}
It suffices to prove the assertion for $c=1$.  By
\eqref{eq:power-sum-scaling-generators}, the algebra endomorphism
$\mathcal D_{p,u}$ fixes $p_1$ and scales the other power-sum
variables independently of $p_1$.  It therefore commutes with both
multiplication by $p_1$ and $p_1^\perp=\partial/\partial p_1$.
The Pieri rule and its Hall adjoint give
\[
p_1s_\lambda=\sum_{\substack{\mu\supset\lambda\\
|\mu/\lambda|=1}}s_\mu,
\qquad
p_1^\perp s_\lambda=\sum_{\substack{\nu\subset\lambda\\
|\lambda/\nu|=1}}s_\nu.
\]
Thus these operators add or remove one node of any residue.

Since $\mathcal D_{p,u}$ preserves degree,
Proposition~\ref{prop:d-core-span-preservation} shows that it
preserves every fixed-degree $p$-core block, and hence every weight
space.  It follows that $\mathcal D_{p,u}$ commutes with each
$\operatorname{pr}_\varpi$: acting on a Schur expansion cannot
move terms from one weight space to another.

Now let $F\in\mathcal W_\varpi$.  By
\eqref{eq:fock-affine-weight}, adding an $i$-node changes the weight
to $\varpi-\alpha_i$, whereas removing an $i$-node changes it to
$\varpi+\alpha_i$.  Different residues give different target
weights.  Consequently, projecting the expansions for $p_1F$ and
$p_1^\perp F$ onto these target spaces retains exactly the terms
in \eqref{eq:refaudit-chevalley-operators}:
\begin{equation}\label{eq:chevalley-weight-projections}
f_iF=\operatorname{pr}_{\varpi-\alpha_i}(p_1F),
\qquad
e_iF=\operatorname{pr}_{\varpi+\alpha_i}(p_1^\perp F).
\end{equation}
Using \eqref{eq:chevalley-weight-projections} and the commutations
just proved, we obtain
\[
\begin{aligned}
\mathcal D_{p,u}(f_iF)
&=\mathcal D_{p,u}\bigl(\operatorname{pr}_{\varpi-\alpha_i}(p_1F)\bigr)\\
&=\operatorname{pr}_{\varpi-\alpha_i}\bigl(p_1\mathcal D_{p,u}(F)\bigr)\\
&=f_i\mathcal D_{p,u}(F).
\end{aligned}
\]
The last equality uses $\mathcal D_{p,u}(F)\in\mathcal W_\varpi$.
Replacing $p_1$ by $p_1^\perp$ and $\varpi-\alpha_i$ by
$\varpi+\alpha_i$ proves the assertion for $e_i$.  Since the Schur
basis is the union of the bases of the weight spaces, both
commutations hold on all of $\Lambda_{\mathbb Q}$.  Taking powers
and dividing by $c!$ over $\mathbb Q$ proves
\eqref{eq:d-chevalley-commutation} for $c\geq1$; for $c=0$ it is
immediate.
\end{proof}

The Heisenberg action from Subsection~\ref{subsec:power-sum-scaling}
commutes with the derived affine Lie algebra action generated by the
operators in \eqref{eq:refaudit-chevalley-operators}
\cite[Section~2.3 and Lemma~2.4(2)]{HongYacobi}.

By \eqref{eq:theta-as-d} and \eqref{eq:d-chevalley-commutation},
$\Theta_p=\mathcal D_{p,1-p}$ commutes with all divided powers of
the Chevalley operators.

Proposition~\ref{cor:walker-conjecture-47} is
\cite[Theorem~8.11]{JingLiuSkewMN}.  It gives the symmetric-function
case of \cite[Conjecture~4.7]{Walker} for odd primes.  We will use it in
Lemma~\ref{lem:initial-reduction} to exclude the $p$-restricted case.

\begin{proposition}
\label{cor:walker-conjecture-47}
For every odd prime $p$ and every partition $\lambda$,
$\mathfrak{s}_\lambda=s_\lambda$ implies that $\lambda$ is a $p$-core.
\end{proposition}

\subsection{Categorical divided powers}
\label{subsec:divided-powers}

The commutation in Proposition~\ref{prop:d-chevalley-commutation}
relates symmetric functions.  We now use operations on polynomial
representations to transfer simplicity between the characters
connected by the divided powers, under the assumptions
\eqref{eq:extremal-string-identities}.  Lemma~\ref{lem:integral-character}
identifies the integral character group, and
Lemma~\ref{lem:endomorphism-field} shows that every endomorphism
of a simple object is scalar.
Theorem~\ref{thm:divided-power-simple-characters} then gives the
corresponding criterion for two characters to be simple
simultaneously.

We use strict polynomial functors to describe polynomial
representations uniformly as the dimension varies.  For an integer
$d\geq0$, such a functor $M$ assigns a finite-dimensional vector space
$M(V)$ to each finite-dimensional $\mathbb K$-vector space $V$, and
a linear map $M(f):M(V)\to M(W)$ to each linear map $f:V\to W$,
compatibly with identities and composition.  Since $\mathbb K$ is
infinite, the condition that $M$ be homogeneous of degree $d$ means
that, in chosen bases, the entries of $M(f)$ are homogeneous
degree-$d$ polynomials in the entries of $f$, with the same degree
for all $V,W$.  For example, $M(V)=V^{\otimes2}$ and
$M(f)=f^{\otimes2}$ give a functor of degree two.  See
\cite[Sections~4.1--4.3]{HTY} for further background.

Let $\cP_d$ be the category of homogeneous degree-$d$ strict polynomial functors
over $\mathbb K$, and put
$\cP=\bigoplus_{d\geq0}\cP_d$.  For every positive integer $n\geq d$,
evaluation $M\mapsto M(\mathbb K^n)$ identifies $\cP_d$ with the
category of finite-dimensional degree-$d$ polynomial representations
of $GL_n(\mathbb K)$
\cite[Lemma~3.4]{FriedlanderSuslin}.  Thus objects and maps
in a fixed degree may be understood through these representations;
the functor notation lets us change degree without choosing a new
rank at each step.

For $M\in\cP_d$ and a positive integer $n\geq d$, each
$g\in GL_n(\mathbb K)$ acts on $M(\mathbb K^n)$ as $M(g)$.
The formal character of this representation records the dimensions
of its simultaneous eigenspaces for diagonal matrices as the
coefficients of the corresponding weight monomials
\cite[Sections~3.2 and~3.4]{Green}.  These symmetric
polynomials are compatible under setting extra variables to zero
and thus define the formal character of $M$ in $\Lambda$.
Its coefficients are integer dimensions, even though $\mathbb K$
has characteristic $p$.  For $\eta\vdash d$, write $W_\eta$ and $L(\eta)$
for the Weyl and simple strict polynomial functors indexed by $\eta$.
Write $l_\eta$ for the formal character
of $L(\eta)$; a \emph{simple character} means one of the $l_\eta$.

Let $K_0(\cP)$ be the Grothendieck group: it is generated by the
classes $[M]$ of objects of $\cP$, with the relations
\[
[M]=[M']+[M'']
\quad\text{for every short exact sequence}\quad
0\longrightarrow M'\longrightarrow M\longrightarrow M''
\longrightarrow0.
\]
These relations retain the multiplicities of the simple composition
factors and forget how they are joined by extensions.  Write
$\ch[M]$ for the formal character of $M$.  Its additivity on
short exact sequences gives the character map
\[
\ch:K_0(\cP)\longrightarrow\Lambda.
\]
The character-group and endomorphism-ring assertions below hold for
every prime $p$.  When $p$ is odd, for $i\in\mathbb Z/p\mathbb Z$,
let $E_i,F_i$ be the residue-$i$
Chevalley functors on $\cP$ supplied by \cite[Theorem~4]{HTY}.
These functors act on representations and lower or raise the degree
by one, respectively.  They are \emph{exact}, meaning that they
preserve short exact sequences, so they induce additive maps on
$K_0(\cP)$.  The next lemma identifies the integer character group
and relates these maps to the operators in
\eqref{eq:refaudit-chevalley-operators}.

\begin{lemma}\label{lem:integral-character}
The map $\ch$ is an isomorphism of $\mathbb Z$-modules, with
$\ch[W_\eta]=s_\eta$ and $\ch[L(\eta)]=l_\eta$ for every partition
$\eta$.  For odd $p$, the endomorphisms induced by $E_i,F_i$
correspond under this isomorphism to $e_i,f_i$.
\end{lemma}

\begin{proof}
Every object of $\cP_d$ has finite length, since evaluation gives
a finite-dimensional representation.  By
\cite[Theorem~2(i) and the proof of Corollary~2]{HTY}, the classes
$[L(\eta)]$, $\eta\vdash d$, form a $\mathbb Z$-basis of
$K_0(\cP_d)$, and for $\eta\vdash d$ one has
\begin{equation}\label{eq:refaudit-weyl-simple-triangularity}
[W_\eta]=[L(\eta)]+
 \sum_{\substack{\nu\vdash d\\\nu<_{\mathrm{lex}}\eta}}
 m_{\eta\nu}[L(\nu)]
\qquad(m_{\eta\nu}\in\mathbb Z_{\geq0}).
\end{equation}
The change of basis in \eqref{eq:refaudit-weyl-simple-triangularity}
is unitriangular over $\mathbb Z$, so the $[W_\eta]$ also form a
$\mathbb Z$-basis.  Their characters are $s_\eta$, by the
identification with Weyl modules in \cite[Proposition~2(ii)]{HTY}.
Thus $\ch$ sends a $\mathbb Z$-basis to the Schur basis and is an
isomorphism; $\ch[L(\eta)]=l_\eta$ holds by definition.

For odd $p$, \cite[Proposition~7]{HTY} identifies
the actions of $E_i,F_i$ with $e_i,f_i$ after tensoring $K_0(\cP)$
with $\mathbb C$.
Exactness makes $E_i,F_i$ act on $K_0(\cP)$ over $\mathbb Z$,
and \eqref{eq:refaudit-chevalley-operators} shows that $e_i,f_i$
preserve $\Lambda$.  Since $\Lambda$ embeds in
$\Lambda\otimes_{\mathbb Z}\mathbb C$, the operator identities
already hold over $\mathbb Z$.
\end{proof}

Lemma~\ref{lem:integral-character} relates the functors $E_i,F_i$
to the operators $e_i,f_i$ on characters.  We next show in
Lemma~\ref{lem:endomorphism-field} that every endomorphism of a
simple object of $\cP$ is a scalar multiple of the identity, as
required in the proof of
Theorem~\ref{thm:divided-power-simple-characters}.

Write $S_\eta^{\mathrm{pol}}$ for the Schur functor indexed by $\eta$.
For a functor $M\in\cP$ and a finite-dimensional
$\mathbb K$-vector space $V$, its Kuhn dual is
$M^\#(V)=M(V^*)^*$, where $*$ denotes $\mathbb K$-linear duality.
This duality reverses maps and sends injections to surjections.
Write $\soc M$ for the socle of $M$, the sum of its simple subobjects.
In \cite[Definition~4(2)--(3) and Theorem~2(i)]{HTY}, the Weyl
functor is defined by Kuhn duality, and the socle of the Schur
functor is shown to be simple.  In our notation,
\[
W_\eta=(S_\eta^{\mathrm{pol}})^\#,
\qquad
L(\eta)=\soc S_\eta^{\mathrm{pol}}.
\]

\begin{lemma}\label{lem:endomorphism-field}
For every partition $\eta$,
\[
\operatorname{End}_{\cP}(L(\eta))=\mathbb K.
\]
\end{lemma}

\begin{proof}
The socle inclusion gives $u:L(\eta)\hookrightarrow
S_\eta^{\mathrm{pol}}$.  Dualizing $u$ and using
$L(\eta)^\#\cong L(\eta)$ from \cite[Theorem~2(ii)]{HTY} gives
a surjection $q:W_\eta\twoheadrightarrow L(\eta)$.  Thus we have
\[
W_\eta\xrightarrow{q}L(\eta)\xrightarrow{u}S_\eta^{\mathrm{pol}}.
\]
Composition defines a $\mathbb K$-linear map
\[
\Phi:\operatorname{End}_{\cP}(L(\eta))
 \longrightarrow
 \operatorname{Hom}_{\cP}(W_\eta,S_\eta^{\mathrm{pol}}),
\qquad
\Phi(a)=u\circ a\circ q.
\]
The map $\Phi$ is injective: since $u$ is injective and $q$ is
surjective, $u\circ a\circ q=0$ forces $a=0$.
By \cite[Theorem~2(iv)]{HTY},
\[
\operatorname{Hom}_{\cP}(W_\eta,S_\eta^{\mathrm{pol}})\cong\mathbb K.
\]
Hence $\dim_{\mathbb K}\operatorname{End}_{\cP}(L(\eta))\leq1$.  Scalar
multiples of the identity give a one-dimensional subspace, so equality
holds and the endomorphism ring is $\mathbb K$.
\end{proof}

With Lemmas~\ref{lem:integral-character} and~\ref{lem:endomorphism-field}
in place, we now relate divided powers to simplicity of characters.
For Theorem~\ref{thm:divided-power-simple-characters}, assume that
$p$ is odd.  Let $\Psi:\Lambda_{\mathbb Q}\to\Lambda_{\mathbb Q}$ be a
degree-preserving $\mathbb Q$-linear map that is injective,
satisfies $\Psi(\Lambda)\subseteq\Lambda$, and commutes with every
$e_i$ and $f_i$.  Fix partitions $\alpha,\beta$, a residue
$i\in\mathbb Z/p\mathbb Z$, and an integer $c\geq0$ satisfying
\begin{equation}\label{eq:extremal-string-identities}
f_i^{(c)}s_\alpha=s_\beta,\quad
f_i^{(c+1)}s_\alpha=0,\qquad
e_i^{(c)}s_\beta=s_\alpha,\quad
e_i^{(c+1)}s_\beta=0.
\end{equation}
The first and third identities connect the two Schur functions.
The vanishing of the next divided powers will identify the last
nonzero steps when we pass to representations.

\begin{theorem}[Divided powers and simple characters]
\label{thm:divided-power-simple-characters}
Under these assumptions,
\begin{equation}\label{eq:divided-power-simplicity}
\Psi(s_\alpha)\text{ is a simple character}
\quad\Longleftrightarrow\quad
\Psi(s_\beta)\text{ is a simple character}.
\end{equation}
\end{theorem}

\begin{proof}
For $c=0$, \eqref{eq:extremal-string-identities} gives
$\alpha=\beta$, so the assertion is immediate.  Assume $c\geq1$.
We first realize the divided-power operators on characters by
functors on representations.  We then use
\eqref{eq:extremal-string-identities} to identify the last nonzero
divided powers of simple objects and prove both implications.

For each integer $k\geq1$, the construction of
\cite[Theorem~4]{HTY} gives actions of the symmetric group $\Sym_k$
on the $k$-fold compositions $E_i^k,F_i^k$.
Following \cite[Sections~3.1.4 and~5.2.2]{CR}, we define the
divided-power functors $E_i^{(k)},F_i^{(k)}$ as the images of the
sums of all permutation operators in the respective actions.
For $k=0$, both are the identity functor.
By \cite[Lemma~5.13]{CR}, they are exact direct summands satisfying
\begin{equation}\label{eq:categorical-divided-power-sums}
F_i^k\cong\bigl(F_i^{(k)}\bigr)^{\oplus k!},
\qquad
E_i^k\cong\bigl(E_i^{(k)}\bigr)^{\oplus k!}.
\end{equation}
Here $\oplus k!$ means a direct sum of $k!$ copies.  Taking
characters in \eqref{eq:categorical-divided-power-sums} and using
Lemma~\ref{lem:integral-character} gives, for every $M\in\cP$,
$k!\ch[F_i^{(k)}M]=f_i^k\ch[M]$, and similarly for $E_i$.
Viewing the integral character ring inside $\Lambda_{\mathbb Q}$,
we can cancel $k!$ and obtain
\begin{equation}\label{eq:categorical-divided-power-characters}
\ch[F_i^{(k)}M]=f_i^{(k)}\ch[M],
\qquad
\ch[E_i^{(k)}M]=e_i^{(k)}\ch[M].
\end{equation}
Thus the division takes place in the rational character space;
no inverse of $k!$ in $\mathbb K$ is needed.

Suppose that $\Psi(s_\alpha)=\ch[L]$ for a simple object $L$.
Using \eqref{eq:categorical-divided-power-characters}, commutation
with $f_i$, and \eqref{eq:extremal-string-identities}, we obtain
\begin{equation}\label{eq:divided-power-character-images}
\begin{aligned}
\ch[F_i^{(c)}L]
&=f_i^{(c)}\Psi(s_\alpha)
 =\Psi(f_i^{(c)}s_\alpha)=\Psi(s_\beta)\neq0,\\
\ch[F_i^{(c+1)}L]
&=\Psi(f_i^{(c+1)}s_\alpha)=0.
\end{aligned}
\end{equation}
The nonvanishing uses injectivity of $\Psi$.  An actual object of
$\cP$ has zero character only if it is zero: by
Lemma~\ref{lem:integral-character} its Grothendieck class is then
zero, whereas a nonzero object has a composition series with
positive length.  Hence \eqref{eq:divided-power-character-images}
implies
\[
F_i^{(c)}L\neq0,
\qquad
F_i^{(c+1)}L=0.
\]
By \eqref{eq:categorical-divided-power-sums}, the same nonvanishing
and vanishing hold for $F_i^cL$ and $F_i^{c+1}L$.  Thus $c$ is the
last nonzero step of repeated application of $F_i$ to $L$.
The $F$-version of \cite[Proposition~5.20(a)]{CR} applies to
the functors supplied by \cite[Theorem~4]{HTY}; its condition
on endomorphisms of simple objects follows from
Lemma~\ref{lem:endomorphism-field}.
Thus $F_i^{(c)}L$ is simple.  Its character in
\eqref{eq:divided-power-character-images} is $\Psi(s_\beta)$,
proving one implication.

Conversely, if $\Psi(s_\beta)=\ch[L']$ for a simple object $L'$,
the two identities involving $e_i$ in
\eqref{eq:extremal-string-identities} give
\[
\ch[E_i^{(c)}L']=\Psi(s_\alpha)\neq0,
\qquad
\ch[E_i^{(c+1)}L']=0.
\]
Lemma~\ref{lem:integral-character} and
\eqref{eq:categorical-divided-power-sums} therefore give
$E_i^cL'\neq0$ and $E_i^{c+1}L'=0$.  Thus $c$ is the last
nonzero step for $E_i$ on $L'$.  The $E$-version of
\cite[Proposition~5.20(a)]{CR} makes $E_i^{(c)}L'$ simple, with
character $\Psi(s_\alpha)$.  This proves the reverse implication.
\end{proof}

Proposition~\ref{prop:adjacent-simplicity} will apply
Theorem~\ref{thm:divided-power-simple-characters} with $\Psi=\Theta_p$.
Proposition~\ref{prop:d-chevalley-commutation} already supplies the
required commutation, and the identities
\eqref{eq:extremal-string-identities} will be verified in that
application.

\subsection{Restricted constituents and the remainder}
\label{subsec:restricted-projection}

In this subsection, $p$ may be any prime.  We determine which simple character could equal
$\mathfrak{s}_\lambda$.  Proposition~\ref{prop:restricted-comparison}
identifies its coefficients at all $p$-restricted simple labels with
Specht-module composition multiplicities.  If the Specht module is
reducible, this already excludes simplicity.  Otherwise it leaves one
candidate, and Proposition~\ref{prop:restricted-projection} reduces
the question to the vanishing of a virtual remainder.

A partition is \emph{$p$-regular} if no positive part occurs at least
$p$ times, and \emph{$p$-restricted} if its conjugate is
$p$-regular.  To compare Specht-module multiplicities with
simple-character coefficients, we use $p$-restricted labels for
the simple symmetric group modules.

In this section $S^\tau$ denotes the row Specht module over
$\mathbb K$, in place of the complex modules used in Section~\ref{sec:modular-schur-coefficients}.
For $\tau\vdash d$ which is $p$-regular,
write $D^\tau$ for its simple quotient; these give all simple
$\mathbb K\Sym_d$-modules
\cite[Theorems~4.9, 11.1, and~11.5]{JamesBook}.  Write
$\operatorname{sgn}$ for the sign representation.  For every partition
$\lambda$, set
\begin{equation}\label{eq:dual-specht-convention}
S_\lambda=S^{\lambda'}\otimes\operatorname{sgn}
\cong (S^\lambda)^*.
\end{equation}
The isomorphism is the standard duality isomorphism for Specht modules
\cite[Theorem~8.15]{JamesBook}.
If $\nu$ is $p$-restricted, define
$D_\nu=D^{\nu'}\otimes\operatorname{sgn}$.  Tensoring with the
sign representation preserves simple modules and quotient maps,
and is its own inverse.  Hence $D_\nu$ is the
unique simple quotient of $S_\nu$, and these labels again give all
simple modules.  Thus tensoring with the sign representation and
conjugating the label pass between the row and column conventions.
For a $p$-regular partition $\tau$,
the Mullineux involution $\iota_p$ is characterized by
\begin{equation}\label{eq:mullineux-convention}
D^{\iota_p(\tau)}\cong D^\tau\otimes\operatorname{sgn}.
\end{equation}

For a finite-length module $M$ and a simple module $D$, write
$[M:D]$ for the multiplicity of $D$ in a composition series of $M$.
For $\lambda\vdash d$, define integers $q_{\lambda\nu}$ using the
integral simple-character basis of Lemma~\ref{lem:integral-character}:
\[
 \mathfrak{s}_\lambda=\sum_{\nu\vdash d}q_{\lambda\nu}l_\nu.
\]
Proposition~\ref{prop:restricted-comparison} relates these
polynomial-character coefficients to composition multiplicities
of symmetric group modules.

\begin{proposition}
\label{prop:restricted-comparison}
For every $\lambda\vdash d$ and every $p$-restricted partition $\nu\vdash d$,
\begin{equation}\label{eq:walker-restricted-comparison}
q_{\lambda\nu}=[S_\lambda:D_\nu].
\end{equation}
\end{proposition}

\begin{proof}
By \cite[Proposition~4.2 and the preceding discussion, p.~601]{Walker},
for $p$-restricted $\nu$, the coefficient of $l_\nu$ in
$\mathfrak{s}_\lambda$ equals its coefficient in $s_\lambda$,
which is $[S_\lambda:D_\nu]$.
\end{proof}

The Carter condition requires the $p$-adic valuations of the hook
lengths to be constant down every column.  Such a partition is
necessarily $p$-regular: if a positive part occurred at least $p$
times, the hook lengths in the last column of that block of equal rows
would include $p$ and $1$, contradicting the column condition.
Theorem~\ref{thm:walker-carter-input} gives the simple character in this case and
excludes simplicity for the other $p$-regular partitions.

\begin{theorem}
\label{thm:walker-carter-input}
If $\lambda$ satisfies the Carter condition, then
\begin{equation}\label{eq:walker-carter-character}
\mathfrak{s}_\lambda=l_{\iota_p(\lambda)'}.
\end{equation}
If $\lambda$ is $p$-regular but does not satisfy the Carter condition,
then $\mathfrak{s}_\lambda$ is not the character of a simple strict
polynomial functor.
\end{theorem}

Theorem~\ref{thm:walker-carter-input} is proved in
\cite[Theorem~4.4]{Walker}, whose nonsimplicity assertion is stated
for all sufficiently large ranks.
Indeed, if $\mathfrak{s}_\lambda=l_\mu$ in $\Lambda$, then for every
$n\geq|\lambda|$ its specialization would be the character of the
simple polynomial $GL_n(\mathbb K)$-module $L(\mu)(\mathbb K^n)$.
Taking $n$ sufficiently large
contradicts Walker's conclusion.

We now seek the nonsimplicity conclusion of
Theorem~\ref{thm:walker-carter-input}
without the $p$-regularity hypothesis.
We separate the terms with $p$-restricted labels, already determined
by Proposition~\ref{prop:restricted-comparison}, from the remaining terms.

For a virtual polynomial character
$x=\sum_{\nu\vdash d}c_\nu l_\nu$, define its restricted projection
by retaining only the terms with $p$-restricted labels:
\begin{equation}\label{eq:restricted-projection-operator}
\operatorname{pr}_{\mathrm{res}}(x)
=\sum_{\substack{\nu\vdash d\\\nu\text{ is }p\text{-restricted}}}
c_\nu l_\nu.
\end{equation}

By Proposition~\ref{prop:restricted-comparison}, for every
$\lambda\vdash d$,
\begin{equation}\label{eq:restricted-projection}
\operatorname{pr}_{\mathrm{res}}(\mathfrak{s}_\lambda)
=\sum_{\substack{\nu\vdash d\\\nu\text{ is }p\text{-restricted}}}
[S_\lambda:D_\nu]l_\nu.
\end{equation}
For every partition $\lambda$, define
$R_\lambda=\mathfrak{s}_\lambda-
\operatorname{pr}_{\mathrm{res}}(\mathfrak{s}_\lambda)$.
Thus $R_\lambda$ collects all terms with non-$p$-restricted labels.
It is a virtual character whose simple-character coefficients may be negative.

\begin{proposition}
\label{prop:restricted-projection}
For every partition $\lambda$, if $S_\lambda$ is reducible, then
$\mathfrak{s}_\lambda$ is not a simple character.  If $S_\lambda$
is irreducible, then
\[
\mathfrak{s}_\lambda\text{ is a simple character}
\quad\Longleftrightarrow\quad R_\lambda=0.
\]
\end{proposition}

\begin{proof}
If $S_\lambda$ is reducible, the sum of the nonnegative
multiplicities in \eqref{eq:restricted-projection} is its composition
length, which is at least two.  The projection of a simple character
is either zero or a single $l_\nu$, so
$\mathfrak{s}_\lambda$ cannot be simple.

If $S_\lambda\cong D_\mu$ for a $p$-restricted partition $\mu$, then
\eqref{eq:restricted-projection} gives
\begin{equation}\label{eq:remainder}
\mathfrak{s}_\lambda=l_\mu+R_\lambda,
\qquad
R_\lambda\in
\operatorname{span}_{\mathbb Z}
\{l_\nu:\nu\vdash|\lambda|,\ \nu\text{ is not }p\text{-restricted}\}.
\end{equation}
No term of $R_\lambda$ has label $\mu$.
Thus $R_\lambda=0$ implies $\mathfrak{s}_\lambda=l_\mu$.
Conversely, if $\mathfrak{s}_\lambda$ is simple, its restricted
projection $l_\mu$ forces it to be $l_\mu$ itself, so
$R_\lambda=0$.
\end{proof}

To describe the cases treated next, we introduce two further
conditions on partitions.  A partition is \emph{doubly singular}
if it is neither $p$-regular nor $p$-restricted.  A partition $\lambda$
is a \emph{$p$-JM partition}
if there do not exist three nodes $(a,b),(a,y),(x,b)$ in its diagram
such that $v_p(h_\lambda(a,b))>0$ and this valuation differs from both
$v_p(h_\lambda(a,y))$ and $v_p(h_\lambda(x,b))$.
The term $p$-JM partition agrees with $(p,p)$-JM partition in
\cite[Conjecture~1.1 and the following discussion]{Fayers}.
In that terminology, partitions satisfying the Carter condition and
those whose conjugates satisfy it are called $(p,p)$-Carter and
conjugate $(p,p)$-Carter partitions, respectively.  Equivalently,
they are the $p$-regular and $p$-restricted $p$-JM partitions,
respectively.

For odd $p$, Proposition~\ref{prop:rouquier-reduction} reduces
simplicity for doubly singular $p$-JM partitions to Rouquier blocks,
and Theorem~\ref{thm:rouquier-nonvanishing} proves that $R_\lambda\neq0$
there.  Subsection~\ref{subsec:carter-classification} completes the
classification: Lemma~\ref{lem:initial-reduction} reduces any
counterexample not satisfying the Carter condition for odd $p$ to
this case, while
Proposition~\ref{prop:prime-two} treats $p=2$.

\subsection{Adjacent blocks and reduction to Rouquier blocks}
\label{subsec:rouquier-reduction}

For this subsection and the next, assume that $p$ is odd.
We reduce the simplicity question for doubly singular $p$-JM
partitions to Rouquier blocks, whose quotient structure will be used
in Subsection~\ref{subsec:rouquier-nonvanishing} to prove nonvanishing
of the remainder.  Lemma~\ref{lem:adjacent-regularity} and
Proposition~\ref{prop:adjacent-simplicity} show that passage between
adjacent blocks preserves double singularity and simplicity,
respectively.  Proposition~\ref{prop:rouquier-reduction} then applies
these results along a sequence ending in a Rouquier block.

Here adjacency has the following meaning.  Fix a $p$-block $B$ of
weight $w$ and choose an abacus with $b_j>0$ beads on runner $j$
for $0\leq j<p$.  By adding equally many beads on all runners,
we may take the display large enough to represent every partition in $B$.
Choose $0\leq k<p$, and put
\[
 c=\begin{cases}
 b_{k-1}-b_k,&k>0,\\
 b_{p-1}-b_0+1,&k=0.
 \end{cases}
\]
The additional $1$ when $k=0$ accounts for moving from runner $p-1$
to runner $0$ at the next level of the abacus.
When $c>0$, let $C$ be the weight-$w$ block with runner counts
$b_{k-1}-c$ on runner $k-1$, $b_k+c$ on runner $k$, and the same
counts as $B$ on the other runners; indices are read modulo $p$.
These blocks are called adjacent.  If the partitions of $B$ have
size $d$, those of $C$ have size $d+c$.

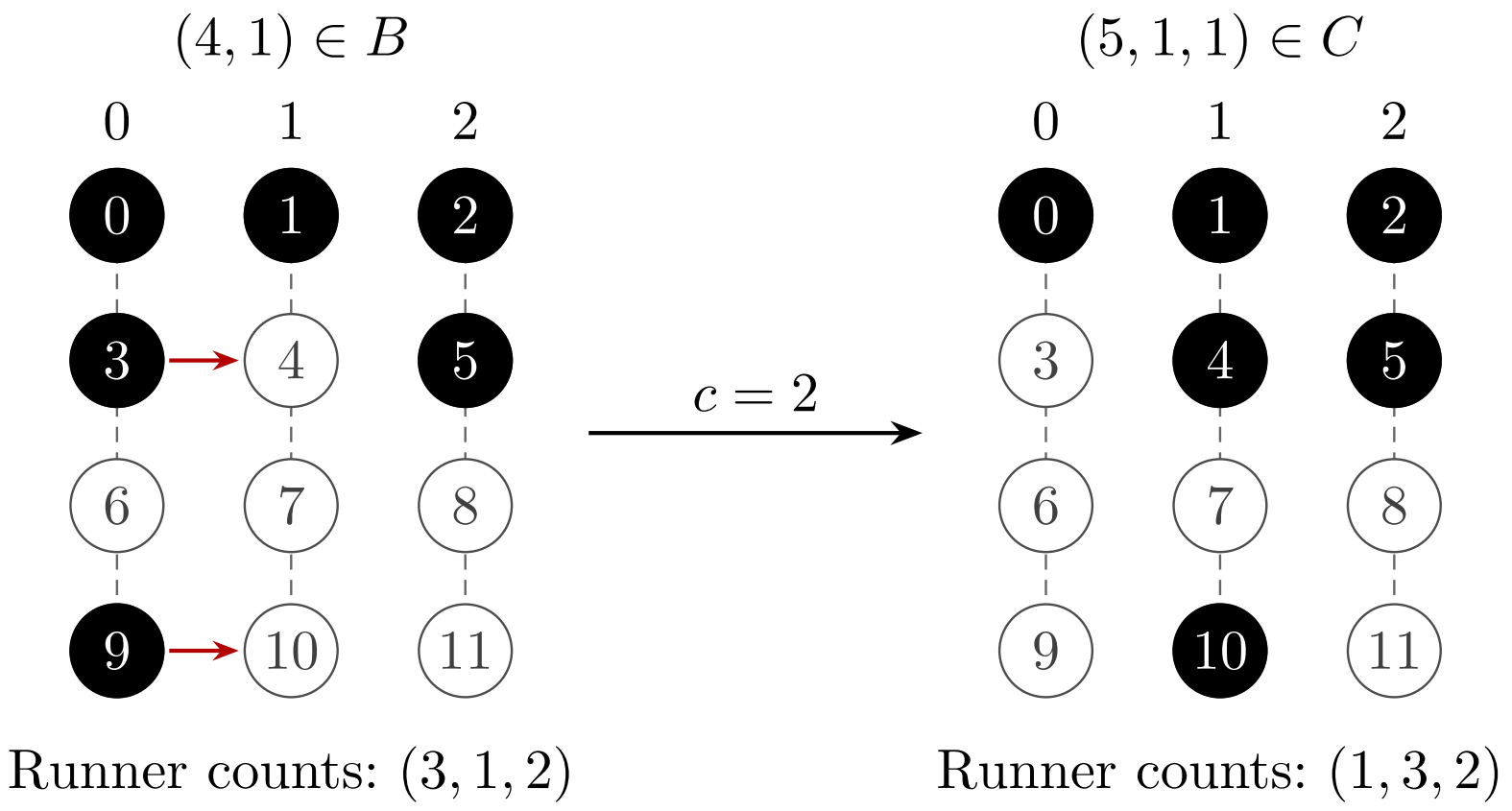
\begin{figure}[!htbp]
\centering
\begin{tikzpicture}[
  x=1cm,y=1cm,
  space/.style={circle,draw=black!70,line width=0.45pt,fill=white,
                text=black!75,minimum size=6.4mm,inner sep=0pt,font=\small},
  bead/.style={circle,draw=black,line width=0.45pt,fill=black,text=white,
               minimum size=6.4mm,inner sep=0pt,font=\small}
]
\node[font=\small] at (1.2,1.2) {$(4,1)\in B$};
\node[font=\small] at (7.6,1.2) {$(5,1,1)\in C$};
\foreach \x/\i in {0/0,1.2/1,2.4/2,6.4/0,7.6/1,8.8/2} {
  \draw[densely dashed,black!60,line width=0.45pt] (\x,0.3)--(\x,-3.3);
  \node[font=\small] at (\x,0.65) {$\i$};
}
\foreach \x/\y/\n in {0/0/0,1.2/0/1,2.4/0/2,
                        0/-1/3,1.2/-1/4,2.4/-1/5,
                        0/-2/6,1.2/-2/7,2.4/-2/8,
                        0/-3/9,1.2/-3/10,2.4/-3/11} {
  \node[space] at (\x,\y) {\n};
  \node[space] at ({\x+6.4},\y) {\n};
}
\foreach \x/\y/\n in {0/0/0,1.2/0/1,2.4/0/2,
                        0/-1/3,2.4/-1/5,0/-3/9} {
  \node[bead] at (\x,\y) {\n};
}
\foreach \x/\y/\n in {6.4/0/0,7.6/0/1,8.8/0/2,
                        7.6/-1/4,8.8/-1/5,7.6/-3/10} {
  \node[bead] at (\x,\y) {\n};
}
\foreach \y in {-1,-3} {
  \draw[-{Stealth[length=1.8mm]},red!70!black,thick]
    (0.36,\y)--(0.84,\y);
}
\draw[-{Stealth[length=2.2mm]},thick] (3.25,-1.5)--(5.55,-1.5)
  node[midway,above,font=\small] {$c=2$};
\node[font=\small] at (1.2,-3.85) {Runner counts: $(3,1,2)$};
\node[font=\small] at (7.6,-3.85) {Runner counts: $(1,3,2)$};
\end{tikzpicture}
\caption{The moves $3\mapsto4$ and $9\mapsto10$ between representatives
of adjacent $3$-blocks of weight $1$.  Filled circles mark beads and
open circles vacant positions; numbers in circles are abacus positions.}
\label{fig:adjacent-block-example}
\end{figure}

For example, take $p=3$, $w=1$, and let $B$ be the block containing
$(4,1)$.  The six-bead display of $(4,1)$ has runner counts $(3,1,2)$.
Choosing $k=1$ gives $c=3-1=2$.  Moving the beads at positions $3$
and $9$ to positions $4$ and $10$ gives the display of $(5,1,1)$,
with runner counts $(1,3,2)$, as shown in
Figure~\ref{fig:adjacent-block-example}.
The two cores are $(1,1)$ and $(2,1,1)$, respectively, so both
partitions have weight $1$ by \eqref{eq:core-weight-size}.
Thus $(5,1,1)$ lies in the adjacent block $C$, and the sizes of the
partitions in $B$ and $C$ are $5$ and $7=5+c$.

Let $B,C$ be adjacent blocks as above, and let $\alpha\in B$ be a
$p$-JM partition.
By \cite[Section~3 and the proof of Lemma~3.2]{Fayers}, there is a
residue $i$ for which $\alpha$ has exactly $c$ addable $i$-nodes and
no removable $i$-node.  Adjoining all these nodes gives a $p$-JM
partition $\beta\in C$, whose removable $i$-nodes are precisely the
$c$ added nodes and which has no addable $i$-node.

In the example of Figure~\ref{fig:adjacent-block-example},
$\alpha=(4,1)$ and $\beta=(5,1,1)$, with $i=1$ and $c=2$.
The two addable $1$-nodes of $\alpha$ are $(1,5)$ and $(3,1)$;
adjoining them gives $\beta$, and removing them from $\beta$
recovers $\alpha$.
For Lemma~\ref{lem:adjacent-regularity} and
Proposition~\ref{prop:adjacent-simplicity}, fix any pair
$\alpha,\beta$ obtained by the construction above.

\begin{lemma}\label{lem:adjacent-regularity}
The partition $\alpha$ is $p$-regular if and only if $\beta$ is
$p$-regular.  Likewise, $\alpha$ is $p$-restricted if and only if
$\beta$ is $p$-restricted.
\end{lemma}

\begin{proof}
We prove a slightly more general elementary fact.  Here all residues
are taken modulo $p$.  Let $\xi$ be any partition with no
removable node of residue $i$, and let $\zeta$ be obtained by adjoining
all its addable nodes of residue $i$.  We claim that
\begin{equation}\label{eq:refaudit-adjoining-regularity}
\xi\text{ is }p\text{-regular}
\quad\Longleftrightarrow\quad
\zeta\text{ is }p\text{-regular}.
\end{equation}

For $L\geq1$, write $m_L(\xi)$ and $m_L(\zeta)$ for the
multiplicities of $L$ in $\xi$ and $\zeta$, respectively.
If $m_L(\xi)>0$, the rows of length $L$ occupy positions $a_L,\ldots,b_L$,
where $b_L=a_L+m_L(\xi)-1$.  The removable node at the end of this group has
residue
\[
r_L=\overline{L-b_L},
\]
whereas the addable node at the end of its first row has residue
\[
q_L=\overline{L+1-a_L}.
\]
The hypothesis says that $r_L\neq i$ for each part $L$ occurring in
$\xi$.  For such $L$, put $\delta_L=1$ if $q_L=i$ and $\delta_L=0$
otherwise; set $\delta_L=0$ if $L$ is absent.
Finally, $\delta_0=1$ if the addable node creating a new last row
has residue $i$, and $\delta_0=0$ otherwise.  Adjoining all these
outer corners gives, for every $L\geq1$,
\begin{equation}\label{eq:multiplicity-change}
m_L(\zeta)=m_L(\xi)-\delta_L+\delta_{L-1}.
\end{equation}
The term $-\delta_L$ records a row changing from length $L$ to $L+1$,
and $+\delta_{L-1}$ records a row changing from $L-1$ to $L$.
Thus a multiplicity can change by at most one, and we need only check
when it could cross the threshold $p$.

Suppose first that $m_L(\xi)\geq p$, so that $\xi$ is not $p$-regular.
A decrease below $p$ could occur
only if $m_L(\xi)=p$ and $\delta_L=1$.  But then
$q_L-r_L=\overline{m_L(\xi)}=0$, so $r_L=q_L=i$, a contradiction.
Hence a non-$p$-regular partition cannot become
$p$-regular.

Conversely, suppose that $\xi$ is $p$-regular but $\zeta$ is not.
By \eqref{eq:multiplicity-change}, some part $L$ must have
$m_L(\xi)=p-1$, $\delta_{L-1}=1$, and $\delta_L=0$.  The node which
raises a row of length $L-1$ to length $L$ is in row $b_L+1$ and has
residue $\overline{L-b_L-1}=r_L-\overline{1}$.
Since $m_L(\xi)=p-1$, one also has $q_L=r_L-\overline{1}$.
Thus this node and the addable node in row
$a_L$ have the same residue.  The equality $\delta_{L-1}=1$ therefore
forces $\delta_L=1$, a contradiction.  This proves
\eqref{eq:refaudit-adjoining-regularity}.

Equation~\eqref{eq:refaudit-adjoining-regularity}
therefore gives the equivalence of $p$-regularity for $\alpha$ and
$\beta$.  Conjugation changes every residue to its negative and
interchanges $p$-regularity with $p$-restrictedness; moreover,
$\beta'$ is obtained from $\alpha'$ by adjoining all addable
$(-i)$-nodes, and $\alpha'$ has no removable $(-i)$-node.  Applying
\eqref{eq:refaudit-adjoining-regularity} to $\alpha'$ and $\beta'$
gives the equivalence of $p$-restrictedness for $\alpha$ and $\beta$.
\end{proof}

For the same partitions $\alpha,\beta$, we next apply
Theorem~\ref{thm:divided-power-simple-characters} with
$\Psi=\Theta_p$ to preserve simplicity as well.

\begin{proposition}
\label{prop:adjacent-simplicity}
One has
\[
\mathfrak{s}_\alpha\text{ is a simple character}
\quad\Longleftrightarrow\quad
\mathfrak{s}_\beta\text{ is a simple character}.
\]
\end{proposition}

\begin{proof}
To apply Theorem~\ref{thm:divided-power-simple-characters} with
$\Psi=\Theta_p$, we first verify the four identities in
\eqref{eq:extremal-string-identities}.  We then check that $\Theta_p$
is an injective degree-preserving $\mathbb Q$-linear map on
$\Lambda_{\mathbb Q}$, maps $\Lambda$ into itself, and commutes
with every $e_i$ and $f_i$.

By \eqref{eq:refaudit-chevalley-operators}, each sequence of $c$
successive additions of $i$-nodes contributes one Schur function to
$f_i^c s_\alpha$.  Adding an $i$-node changes the status of other
corners only at residues $i-1$ or $i+1$.  Hence, after each addition,
the addable $i$-nodes are precisely those of the original $c$
nodes that have not yet been added.  They can therefore be added in
any order, giving exactly $c!$ sequences, all ending at $\beta$.
Since $\beta$ has no addable $i$-node, none of these sequences can
be continued.  Thus
\begin{equation}\label{eq:fayers-fock}
f_i^c s_\alpha=c!s_\beta,\qquad f_i^{c+1}s_\alpha=0.
\end{equation}
For the reverse direction, the removable $i$-nodes of $\beta$ are
exactly the $c$ nodes just added.  Removing one leaves the others
removable and creates no new removable $i$-nodes, so all $c!$
orders of removal lead to $\alpha$.  Since $\alpha$ has no removable
$i$-node, \eqref{eq:refaudit-chevalley-operators} gives
\begin{equation}\label{eq:fayers-fock-reverse}
e_i^c s_\beta=c!s_\alpha,\qquad e_i^{c+1}s_\beta=0.
\end{equation}
Passing to divided powers in $\Lambda_{\mathbb Q}$,
\eqref{eq:fayers-fock} and \eqref{eq:fayers-fock-reverse} give all
four identities in \eqref{eq:extremal-string-identities}.

We now verify the assumptions on $\Theta_p$.  Its definition in
Subsection~\ref{subsec:modular-schur-functions} as a graded ring
endomorphism of $\Lambda$ shows that its extension to
$\Lambda_{\mathbb Q}$ is $\mathbb Q$-linear, preserves degree, and
maps $\Lambda$ into itself.  By \eqref{eq:theta-power-sums}, it
multiplies each power-sum generator by either $1$ or $1-p$.
Since these generators are algebraically independent and both
multipliers are nonzero, these scalings are invertible over
$\mathbb Q$; in particular, $\Theta_p$ is injective.
Finally, \eqref{eq:theta-as-d} identifies
$\Theta_p=\mathcal D_{p,1-p}$, and
Proposition~\ref{prop:d-chevalley-commutation} gives its commutation
with every $e_i$ and $f_i$.

Since $p$ is odd by our standing assumption, all hypotheses of
Theorem~\ref{thm:divided-power-simple-characters} now hold.
It follows that $\Theta_p(s_\alpha)$ is a simple character if and
only if $\Theta_p(s_\beta)$ is.  By
\eqref{eq:modular-schur-definition}, these are
$\mathfrak{s}_\alpha$ and $\mathfrak{s}_\beta$, respectively.
\end{proof}

We now apply Lemma~\ref{lem:adjacent-regularity} and
Proposition~\ref{prop:adjacent-simplicity} along a sequence of
adjacent blocks ending in a Rouquier block.  With the abacus conventions of
Subsection~\ref{subsec:abacus-core-quotient}, a $p$-block $B$ of weight
$w$ is called a \emph{Rouquier block} if it admits an abacus display with
$c_i$ beads on runner $i$ such that
\begin{equation}\label{eq:rouquier-definition}
c_i-c_{i-1}\geq w-1
\qquad(1\leq i\leq p-1).
\end{equation}

Let $\lambda$ be a $p$-JM partition in a block $B_0$.
By \cite[Lemma~3.1]{Fayers}, there is a sequence of adjacent
blocks $B_0,B_1,\ldots,B_t$ with $B_t$ Rouquier.
Applying \cite[Lemma~3.2]{Fayers} at each step supplies $p$-JM partitions
\[
\xi_0=\lambda,\xi_1,\ldots,\xi_t,
\qquad \xi_j\in B_j\quad(0\leq j\leq t).
\]
Fix any sequence supplied by this construction.

\begin{proposition}
\label{prop:rouquier-reduction}
For every $0\leq j\leq t$, the partition $\xi_j$ is doubly singular
if and only if $\lambda$ is doubly singular, and
\begin{equation}\label{eq:rouquier-reduction-simplicity}
\mathfrak{s}_{\xi_j}\text{ is a simple character}
\quad\Longleftrightarrow\quad
\mathfrak{s}_{\lambda}\text{ is a simple character}.
\end{equation}
\end{proposition}

\begin{proof}
Lemma~\ref{lem:adjacent-regularity} preserves both $p$-regularity and
$p$-restrictedness at each step, and hence preserves double
singularity.  Proposition~\ref{prop:adjacent-simplicity} preserves
simplicity at each step.  Induction along the selected sequence
proves Proposition~\ref{prop:rouquier-reduction}.
\end{proof}

By \eqref{eq:rouquier-reduction-simplicity}, excluding simplicity at
the Rouquier endpoint also excludes it for the original doubly
singular $p$-JM partition.  Subsection~\ref{subsec:rouquier-nonvanishing}
does this by proving nonvanishing of the remainder and applying
Proposition~\ref{prop:restricted-projection}.

\subsection{Nonvanishing of the remainder in Rouquier blocks}
\label{subsec:rouquier-nonvanishing}

For a doubly singular $p$-JM partition $\lambda$ in a Rouquier
block, we first prove that $S_\lambda$ is simple and determine its
$p$-restricted label $\mu$.  This identifies $l_\mu$ as the only
possible simple character equal to $\mathfrak{s}_\lambda$, by
\eqref{eq:remainder} in the proof of
Proposition~\ref{prop:restricted-projection}.  We then prove
$R_\lambda\neq0$ by finding a Schur coefficient of
$\mathfrak{s}_\lambda$ which is zero in $l_\mu$.
To use Proposition~\ref{prop:universal-skew-expansion} for this
calculation, we first compare the runner numbering of a Rouquier
abacus with the quotient convention of
Section~\ref{sec:modular-schur-coefficients}.

Fix a Rouquier block $B$ of weight $w$ and a display of its core
satisfying \eqref{eq:rouquier-definition}, with $c_i$ beads on runner $i$.
To represent every partition in $B$, enlarge the display until $c_i\geq w$
for every $i$.  In each enlargement, \eqref{eq:beta-set-b-plus-r}
moves every old bead down one level and fills the new top level with
one bead on each runner.  This gives another display of the same
partition and preserves the differences $c_i-c_{i-1}$ in
\eqref{eq:rouquier-definition}.

Now fix these runner counts and put $b=\sum_i c_i$.
Every partition in $B$ packs to the same core.  On runner $i$ its
quotient component records how many levels each bead lies below its
packed position, reading the beads from bottom to top.  Thus $c_i$
beads can record a component with up to $c_i$ nonzero parts.
Since the sizes of all the components sum to $w$ by
\eqref{eq:core-weight-size}, the chosen counts suffice for every
partition in $B$.  At this fixed bead number, the partitions in $B$
are represented by varying the vertical positions of the beads on
their respective runners.

For $\eta\in B$, let $\eta^{[i]}$ be the component on runner $i$
of $Q_b(\eta)$ in this display.  To read the component
$\eta^{(i)}$ defined in \eqref{eq:standard-quotient-tuple}, we use
a bead number divisible by $p$.  Choose $0\leq t<p$ so that $p\mid b+t$.
By \eqref{eq:quotient-charge-shift}, increasing the bead number by one
for the same partition carries the component on runner $j$ to runner
$j+1$ modulo $p$.  For instance, in the $(3,2,1)$ example following
\eqref{eq:quotient-charge-shift}, the empty component is on runner $1$
with three beads and on runner $2$ with four beads.
Thus the component $\eta^{[j]}$ is read on runner $j+t$ in $Q_{b+t}(\eta)$.

For example, take $p=3$, $w=2$, and a core display with runner counts
$(1,2,4)$.  Adding one full level gives counts $(2,3,5)$ and $b=10$.
We can now represent the partition $\eta\in B$ whose only nonempty
quotient component in this display is $(1,1)$ on runner $0$, by
moving both beads on that runner down one level from the packed display.
Increasing the bead number from $10$ to $12$ rotates the tuple twice, giving
\[
 Q_{10}(\eta)=\bigl((1,1),\emptypar,\emptypar\bigr),\qquad
 Q_{12}(\eta)=\bigl(\emptypar,\emptypar,(1,1)\bigr).
\]
Thus $\eta^{[0]}=\eta^{(2)}=(1,1)$, whereas $\eta^{(0)}=\emptypar$.

For $\lambda,\zeta\in B$, put
\[
 \widetilde{\boldsymbol\alpha}
 =\bigl((\lambda^{[0]})',\ldots,(\lambda^{[p-1]})'\bigr),\qquad
 \widetilde{\boldsymbol\beta}
 =\bigl((\zeta^{[0]})',\ldots,(\zeta^{[p-1]})'\bigr).
\]

\begin{samepage}
\begin{lemma}
\label{lem:change-bead-number}
For every $\eta\in B$,
\begin{equation}\label{eq:charged-canonical-quotient}
 \eta^{(i)}=\eta^{[i-t]}\qquad(i\in\mathbb Z/p\mathbb Z).
\end{equation}
For all $\lambda,\zeta\in B$,
\begin{equation}\label{eq:coefficient-cyclic-shift}
 a_{\lambda\zeta}
 =\kappa^{(p,w)}_{\widetilde{\boldsymbol\alpha},\widetilde{\boldsymbol\beta}}.
\end{equation}
\end{lemma}
\end{samepage}

\begin{proof}
Repeated application of \eqref{eq:quotient-charge-shift} sends the
component on runner $j$ in $Q_b(\eta)$ to runner $j+t$ in
$Q_{b+t}(\eta)$.  Since $p\mid b+t$, $Q_{b+t}(\eta)$ is the quotient tuple
$\bigl(\eta^{(0)},\ldots,\eta^{(p-1)}\bigr)$ of
\eqref{eq:standard-quotient-tuple}.  This proves
\eqref{eq:charged-canonical-quotient}.  The same cyclic permutation
applies to every $\eta\in B$, since the bead number and runner labels
are fixed.
The kernel
$\mathcal K_p(\mathbf X,\mathbf Y)
=\prod_{i\neq j}D(X^{(i)},Y^{(j)})$
is invariant under a simultaneous permutation of the components of
$\mathbf X$ and $\mathbf Y$.  Hence its coefficient is unchanged by this common
cyclic reindexing.  Proposition~\ref{prop:core-independent-matrix}
now gives \eqref{eq:coefficient-cyclic-shift}.  The normalizing
factors in \eqref{eq:def-a-lambda-mu} depend only on the partitions,
so they too are unchanged by the runner relabelling.
\end{proof}

Equation~\eqref{eq:coefficient-cyclic-shift} allows us to apply
Proposition~\ref{prop:universal-skew-expansion} directly to the
$[i]$-indexed quotients in the proof below.

\begin{theorem}
\label{thm:rouquier-nonvanishing}
Suppose $p$ is odd and $\lambda\in B$ is a doubly singular
$p$-JM partition.  Then $S_\lambda\cong D_\mu$ for a unique $p$-restricted partition
$\mu$, and $R_\lambda\neq0$.
\end{theorem}

\begin{proof}
We first identify the label $\mu$ of $S_\lambda$.  We use the fixed
abacus of Lemma~\ref{lem:change-bead-number} and write
$A=\lambda^{[0]}$ and $B_0=\lambda^{[p-1]}$.
Since $\lambda$ is a $p$-JM partition in a Rouquier block,
\cite[Lemma~4.1]{Fayers} gives
\begin{equation}\label{eq:rouquier-jm-form}
 \lambda^{[0]}=A,\qquad \lambda^{[p-1]}=B_0,\qquad
 \lambda^{[i]}=\emptypar\quad(0<i<p-1).
\end{equation}

By \eqref{eq:dual-specht-convention}, we need the row Specht module
$S^{\lambda'}$.  To read the quotient of $\lambda'$, choose an integer
$C$ large enough that all beads of $\lambda$ lie in $[0,Cp-1]$ and
$C-c_i\geq w$ for every $i$.  Reflecting the vacant positions by
$x\mapsto Cp-1-x$ gives the beads of $\lambda'$.  This reverses the
runners and conjugates their quotient components, so $\lambda'$ has
$B_0'$ on runner $0$ and $A'$ on runner $p-1$.  Its runner counts
$d_i=C-c_{p-1-i}$ satisfy
$d_i-d_{i-1}=c_{p-i}-c_{p-1-i}\geq w-1$, so the conjugate block
is again Rouquier by \eqref{eq:rouquier-definition}.  In such a
display a partition is $p$-regular exactly when its component on runner
$0$ is empty \cite[Section~1.1.3]{Fayers}.  Double singularity means
that neither $\lambda$ nor $\lambda'$ is $p$-regular, so $A$ and $B_0$
are nonempty, and
\eqref{eq:core-weight-size} gives $|A|+|B_0|=w$.

Transposition interchanges rows and columns in the hook-length condition
for a $p$-JM partition, so $\lambda'$ is also $p$-JM.  We may therefore
apply \cite[Proposition~4.3 and its proof]{Fayers} to $\lambda'$,
in the symmetric group case $q=1$ with $e=p>2$.  The proposition gives
the irreducibility of $S^{\lambda'}$.  Its proof identifies the
$p$-regular label $\tau$ of $S^{\lambda'}\cong D^\tau$ by conjugating
the component of $\lambda'$ on runner $0$ and placing it on runner $1$,
while retaining the component on runner $p-1$.  Thus the only
nonempty components of $\tau$ in the conjugate
display are $B_0$ on runner $1$ and $A'$ on runner $p-1$.
Tensoring with the sign representation in \eqref{eq:dual-specht-convention}
now gives $S_\lambda\cong D^\tau\otimes\operatorname{sgn}=D_{\tau'}$.
Since $\tau$ is $p$-regular, $\mu=\tau'$ is $p$-restricted.
If $S_\lambda\cong D_\nu$ for a $p$-restricted partition $\nu$, tensoring
again with $\operatorname{sgn}$ gives $D^{\nu'}\cong D^\tau$.
Distinct $p$-regular partitions label nonisomorphic simple modules, so
$\nu'=\tau$ and hence $\nu=\mu$.
Reversing the runners and conjugating the components once more gives,
in the original display,
\begin{equation}\label{eq:mu-quotient}
 \mu^{[0]}=A,\qquad \mu^{[p-2]}=B_0',\qquad
 \mu^{[i]}=\emptypar\quad(i\neq0,p-2).
\end{equation}
By \eqref{eq:remainder}, $R_\lambda=\mathfrak{s}_\lambda-l_\mu$.
It suffices to find a Schur function with different coefficients in
$\mathfrak{s}_\lambda$ and $l_\mu$.

We choose $\rho$ with the same core as $\lambda$ and with its entire
quotient on runner $p-2$, distinct from the two runners in
\eqref{eq:rouquier-jm-form}.  Using componentwise addition, set
\begin{equation}\label{eq:rho-quotient}
 \rho^{[p-2]}=A'+B_0',\qquad
 \rho^{[i]}=\emptypar\quad(i\neq p-2).
\end{equation}
This tuple has total size $w$, so it specifies a unique $\rho\in B$
by Proposition~\ref{prop:core-quotient-correspondence} and the
reindexing \eqref{eq:charged-canonical-quotient}.  By
\eqref{eq:coefficient-cyclic-shift}, we can use these runner indices
in \eqref{eq:universal-skew-expansion}.  For its two skew Schur factors
on runner $i$ to be nonzero, $\gamma^{(i)}$ must lie in
$(\lambda^{[i]})'$ and $(\gamma^{(i)})'$ must lie in $\rho^{[i]}$.
On each runner at least one of $\lambda^{[i]}$ and $\rho^{[i]}$ is
empty, so every $\gamma^{(i)}$ is empty.  The surviving term gives
\begin{equation}\label{eq:nonzero-rho-coefficient}
 a_{\lambda\rho}
 =\inner{s_{A'}s_{B_0'}}{s_{A'+B_0'}}
 =c^{A'+B_0'}_{A',B_0'}=1.
\end{equation}
Indeed, filling row $j$ of $(A'+B_0')/A'$ with $j$ gives a
Littlewood--Richardson tableau of content $B_0'$.  The lattice-word
condition bounds the entries in row $j$ by $j$.  Since the preceding
rows exhaust the required smaller entries, induction forces this
filling in every row.  By \eqref{eq:def-a-lambda-mu}, the coefficient of $s_\rho$
in $\mathfrak{s}_\lambda$ is therefore
$(-1)^w\eps_p(\lambda)\eps_p(\rho)\in\{1,-1\}$.

To see why $s_\rho$ has coefficient zero in $l_\mu$, we compare
$\rho$ and $\mu$ in lexicographic order.  By \eqref{eq:mu-quotient}
and \eqref{eq:rho-quotient}, their bead positions agree off runners
$0$ and $p-2$.  For $\eta\in B$, the $j$th largest bead position on
runner $i$ is
\begin{equation}\label{eq:runner-beta-set}
 i+p\bigl(c_i-j+\eta^{[i]}_j\bigr)\qquad(1\leq j\leq c_i),
\end{equation}
since the $j$th bead lies $\eta^{[i]}_j$ levels below its packed
position.  In particular, the largest bead position of $\rho$ on
runner $p-2$ is
\[
 x=p-2+p\bigl(c_{p-2}-1+A'_1+(B_0')_1\bigr).
\]
It exceeds the largest position of $\mu$ on that runner by
$pA'_1>0$.  On runner $0$, all positions in either display are at most
$p(c_0-1+A_1)$.  Since $B_0$ is nonempty, $A_1\leq |A|\leq w-1$.
Summing \eqref{eq:rouquier-definition} over the first $p-2$ gaps gives
$c_{p-2}-c_0\geq(p-2)(w-1)\geq A_1$, and hence
\[
 p(c_0-1+A_1)\leq p(c_{p-2}-1)<x.
\]
Thus $x$ is the largest position occupied in exactly one display,
and it belongs to $\rho$.  In the decreasing lists of beta-numbers
the first difference is therefore larger for $\rho$.  Subtracting
the same staircase $(b-1,b-2,\ldots,0)$ recovers the two partitions,
so
\begin{equation}\label{eq:rho-lex}
 \mu<_{\mathrm{lex}}\rho.
\end{equation}

On taking characters in \eqref{eq:refaudit-weyl-simple-triangularity},
the change of basis from Schur functions to simple characters is
unitriangular in lexicographic order.  Its inverse has the same
triangularity, so $l_\mu$ is a linear combination of $s_\mu$ and
$s_\nu$ with $\nu<_{\mathrm{lex}}\mu$.  Thus \eqref{eq:rho-lex}
gives $\Coef_{s_\rho}[l_\mu]=0$.  By \eqref{eq:remainder}, the
coefficient of $s_\rho$ in $R_\lambda$ is the same as in
$\mathfrak{s}_\lambda$, namely $(-1)^w\eps_p(\lambda)\eps_p(\rho)$.
This proves $R_\lambda\neq0$.
\end{proof}

Together with Propositions~\ref{prop:restricted-projection}
and~\ref{prop:rouquier-reduction},
Theorem~\ref{thm:rouquier-nonvanishing} excludes simple modular Schur
characters indexed by doubly singular $p$-JM partitions for odd $p$.

\subsection{Classification of simple modular Schur characters}
\label{subsec:carter-classification}

Let $p$ be any prime and fix a partition $\lambda\vdash d$.
Theorem~\ref{thm:walker-carter-input}
gives the simple character when $\lambda$ satisfies the Carter condition.
To prove necessity, we first treat $p=2$, when \eqref{eq:prime-two-transpose}
gives $\mathfrak{s}_\lambda=s_{\lambda'}$.  For odd $p$, it remains
to show that any counterexample must be a doubly singular $p$-JM
partition, so that Proposition~\ref{prop:rouquier-reduction} and
Theorem~\ref{thm:rouquier-nonvanishing} apply.

\begin{proposition}\label{prop:prime-two}
Suppose $p=2$.  If $\lambda$ does not satisfy the Carter condition,
then $\mathfrak{s}_\lambda\neq l_\mu$ for every partition $\mu$.
\end{proposition}

\begin{proof}
By \eqref{eq:prime-two-transpose} and Lemma~\ref{lem:integral-character},
$\mathfrak{s}_\lambda=s_{\lambda'}$ is the character of the Weyl
functor $W_{\lambda'}$.  In \eqref{eq:refaudit-weyl-simple-triangularity},
$L(\lambda')$ occurs in $W_{\lambda'}$ with multiplicity one.  If
$s_{\lambda'}$ were a simple character, linear independence of the
simple characters would force all other composition multiplicities
to vanish.  Thus $W_{\lambda'}$ would be simple.

Since the empty partition satisfies the Carter condition, $d>0$.
Evaluation at $\mathbb K^d$ would then give a simple Weyl module
$\Delta(\lambda')$ by the equivalence in Subsection~\ref{subsec:divided-powers}.
For the ordinary Schur algebra the parameter is $q=1$, so the quantum
characteristic here is $e=2$.  By \cite[Proposition~1.2]{Fayers},
$\Delta(\lambda')$ is simple exactly when the conjugate of its label,
namely $\lambda$, satisfies the Carter condition.  This contradicts the
hypothesis.
\end{proof}

Now suppose $p$ is odd.  To apply Proposition~\ref{prop:rouquier-reduction},
we determine what simplicity of $\mathfrak{s}_\lambda$ would force on
a partition $\lambda$ that fails the Carter condition.

\begin{lemma}
\label{lem:initial-reduction}
If $\mathfrak{s}_\lambda$ is a simple character and $\lambda$ does
not satisfy the Carter condition, then $\lambda$ is a doubly
singular $p$-JM partition.
\end{lemma}

\begin{proof}
Write $\mathfrak{s}_\lambda=l_\mu$ with $\mu\vdash d$.
By Proposition~\ref{prop:restricted-projection}, $S_\lambda$ is
irreducible, and so is its dual $S^\lambda$ by
\eqref{eq:dual-specht-convention}.  For odd $p$, failure of the
$p$-JM condition implies that $S^\lambda$ is reducible
\cite[Section~2, proof of the \emph{if} part of Conjecture~1.1]{FayersReducible}.
Hence $\lambda$ is $p$-JM.  Theorem~\ref{thm:walker-carter-input}
also rules out $p$-regularity, since $\lambda$ fails the Carter condition.

It remains to exclude $p$-restrictedness.  Suppose $\lambda$ is
$p$-restricted.  By the convention following
\eqref{eq:dual-specht-convention}, $D_\lambda$ is the simple quotient
of $S_\lambda$.  Since $S_\lambda$ is irreducible, it is isomorphic
to $D_\lambda$.  Equation~\eqref{eq:walker-restricted-comparison}
therefore gives $q_{\lambda\lambda}=[S_\lambda:D_\lambda]=1$.
Thus $l_\lambda$ occurs in $\mathfrak{s}_\lambda=l_\mu$, forcing
$\mu=\lambda$.

Since $\lambda$ is both $p$-restricted and $p$-JM, its conjugate
satisfies the Carter condition, as recalled after
Proposition~\ref{prop:restricted-projection}.  By
\cite[Proposition~1.2]{Fayers} at $q=1$ and $e=p$, the Weyl module
of label $\lambda$ in rank $d$ is simple.  Under the evaluation
equivalence this gives $W_\lambda\cong L(\lambda)$, and hence
$s_\lambda=l_\lambda$ by Lemma~\ref{lem:integral-character}.
We have obtained $\mathfrak{s}_\lambda=l_\lambda=s_\lambda$, so
Proposition~\ref{cor:walker-conjecture-47} implies that $\lambda$ is a
$p$-core.  Every hook length of a $p$-core has $p$-adic valuation
zero, which gives the Carter condition in every column.  This
contradiction excludes $p$-restrictedness and proves that $\lambda$ is
doubly singular.
\end{proof}

\begin{theorem}\label{thm:conj46-odd}
Suppose $p$ is odd.  If $\lambda$ does not satisfy the Carter condition, then $\mathfrak{s}_\lambda\neq l_\mu$ for every partition $\mu$.
\end{theorem}

\begin{proof}
Suppose $\lambda$ fails the Carter condition and $\mathfrak{s}_\lambda$
is simple.  By Lemma~\ref{lem:initial-reduction}, $\lambda$ is a doubly
singular $p$-JM partition.  We may therefore apply
Proposition~\ref{prop:rouquier-reduction} along its sequence of adjacent
blocks.  The $p$-JM partition $\xi_t$ in the final Rouquier block is
still doubly singular, and \eqref{eq:rouquier-reduction-simplicity}
shows that $\mathfrak{s}_{\xi_t}$ is simple.  By
Theorem~\ref{thm:rouquier-nonvanishing}, $S_{\xi_t}$ is irreducible
and $R_{\xi_t}\neq0$.  Proposition~\ref{prop:restricted-projection}
therefore excludes simplicity of $\mathfrak{s}_{\xi_t}$, a contradiction.
\end{proof}

Combining the two cases gives the necessity of the Carter condition
for every prime $p$.

\begin{corollary}
\label{cor:conj46}
If $\lambda$ does not satisfy the Carter condition, then $\mathfrak{s}_\lambda\neq l_\mu$ for every partition $\mu$.
\end{corollary}

\begin{proof}
Proposition~\ref{prop:prime-two} proves the assertion for $p=2$, and
Theorem~\ref{thm:conj46-odd} proves it for every odd prime.
\end{proof}

Theorem~\ref{thm:walker-carter-input} supplies the converse and
identifies the simple character.

\begin{corollary}[Simple modular Schur characters]
\label{cor:carter-classification}
The modular Schur function $\mathfrak{s}_\lambda$ is a simple
character if and only if $\lambda$ satisfies the Carter condition.
In that case
\begin{equation}\label{eq:carter-simple-label}
\mathfrak{s}_\lambda=l_{\iota_p(\lambda)'}.
\end{equation}
\end{corollary}

\begin{proof}
Necessity follows from Corollary~\ref{cor:conj46}.  For a partition
satisfying the Carter condition, \eqref{eq:walker-carter-character}
gives the simple character in \eqref{eq:carter-simple-label}.
\end{proof}

\begin{remark}
When $p=2$, the sign representation is trivial, so
\eqref{eq:mullineux-convention} gives $\iota_2(\lambda)=\lambda$
for every $2$-regular $\lambda$.  Thus for $\lambda$ satisfying the
Carter condition, \eqref{eq:carter-simple-label} reads
$\mathfrak{s}_\lambda=l_{\lambda'}$.

For any prime $p$ and any $\lambda$ with $S_\lambda$ irreducible,
Proposition~\ref{prop:restricted-projection} and
Corollary~\ref{cor:carter-classification} give $R_\lambda=0$ exactly
when $\lambda$ satisfies the Carter condition.
\end{remark}

\begin{proof}[Proof of Theorem~\ref{thm:intro-simple-classification}]
Let $n\geq d$ be a positive integer.  By the evaluation
equivalence in Subsection~\ref{subsec:divided-powers}, the simple
degree-$d$ polynomial $GL_n(\mathbb K)$-characters are precisely the
specializations of $l_\mu$ for $\mu\vdash d$.  Specialization to $n$
variables is injective on $\Lambda^d$, since every partition of $d$
has at most $n$ parts and the corresponding Schur polynomials are
linearly independent.  Hence $\mathfrak{s}_\lambda(x_1,\ldots,x_n)$
is a simple character exactly when $\mathfrak{s}_\lambda=l_\mu$
for some $\mu\vdash d$.  By Corollary~\ref{cor:carter-classification},
this is equivalent to the Carter condition on $\lambda$.  In particular,
\cite[Conjecture~4.6]{Walker} holds over the infinite field $\mathbb K$
for every positive $n\geq d$.
\end{proof}

\section{Hall--Littlewood modes and normal-order coefficients}
\label{sec:hl-mode-positivity}

We return to characteristic zero and take $q$ to be an indeterminate.
Romero gave an expansion of the compositional Hall--Littlewood
functions in the complete symmetric function basis with coefficients
in $\mathbb N[q-1]$ \cite[Proposition~5.1]{Romero}.
Here we study the action of normalized repeated modes on the power-sum basis.
In Subsection~\ref{c1:subsec:hall-littlewood} we use the half-vertex
operators of Section~\ref{sec:universal-blocks} to construct the
Hall--Littlewood modes.  In Subsection~\ref{c1:subsec:operator-positivity}
we show that normalized repeated modes preserve power-sum positivity
in $q-1$.  Subsection~\ref{c1:subsec:normal-coefficients} isolates the
coefficients in their normal-ordered expansions and gives an application
to fixed-core spaces.  In Section~\ref{sec:c1-positivity} we combine
the positivity and integrality obtained here with a formula for the
Macdonald creation operators to prove Theorem~\ref{thm:intro-c1}.

\subsection{Hall--Littlewood functions and vertex operators}
\label{c1:subsec:hall-littlewood}

We collect the notation and identities needed in the rest of this section.

Let $\Lambda$ be the ring of symmetric functions over $\mathbb Z$ as
before.  The operators below act on $\mathbb Q(q)\otimes_{\mathbb Z}\Lambda$.
Write $P_\lambda(X;q)$ and $Q_\lambda(X;q)$ for the Hall--Littlewood
functions with parameter $q$.
For $j\geq1$, let $m_j(\lambda)$ be the multiplicity of $j$ in
$\lambda$, and put
\[
 b_\lambda(q)=
 \prod_{j\geq1}\prod_{i=1}^{m_j(\lambda)}(1-q^i),
 \qquad Q_\lambda(X;q)=b_\lambda(q)P_\lambda(X;q).
\]
With these conventions, the $P_\lambda(X;q)$ form a $\mathbb Z[q]$-basis
of $\Lambda[q]$.  In $k$ variables, $P_\lambda$ vanishes whenever
$\ell(\lambda)>k$ \cite[Chapter~III, (2.5)--(2.7)]{Macdonald}.

We use the plethystic notation of
\eqref{bridge:eq:half-fields}--\eqref{bridge:eq:half-commutation}.
Thus $p_m[-X]=-p_m[X]$ and
$p_m[(q-1)X]=(q^m-1)p_m[X]$.  A scalar substitution $z\mapsto-z$
instead gives $p_m[(-z)X]=(-z)^mp_m[X]$.  These conventions
distinguish the negative virtual alphabet $-X$ from negating every
letter of $X$.
For a partition $\lambda$ and an integer $k\geq0$, define
\[
 \mathcal R_\lambda(X;q)
 =(-1)^{\ell(\lambda)}P_\lambda[-X;q],
 \qquad
 B_k(q)=\prod_{j=1}^k(q^j-1),
 \qquad B_0(q)=1.
\]
In particular, $\mathcal R_\varnothing=1$.

Using $\partial_{p_m}$ for differentiation with respect to the
algebraically independent power sums, define
\begin{equation}\label{c1:eq:jing-field}
 V(z)=
 \exp\left(\sum_{m\geq1}\frac{q^m-1}{m}p_mz^m\right)
\exp\left(\sum_{m\geq1}\partial_{p_m}z^{-m}\right).
\end{equation}
In terms of \eqref{bridge:eq:half-fields}, this is
\[
 V(z)=\mathcal H_+((q-1)z)\mathcal H_-(z^{-1}),
\]
where $p_m[(q-1)z]=(q^m-1)z^m$.
For an integer $a$, its mode $V_a$ is the coefficient $\Coef_{z^a}[V(z)]$.
For $a\geq1$ and $k\geq0$, set
\begin{equation}\label{c1:eq:normalized-modes}
 E_{a,k}=\frac{V_a^k}{B_k(q)},
 \qquad E_{a,0}=I,
\end{equation}
where $I$ is the identity operator.  On an input of degree $d$, only
finitely many terms of $V_a$ act nontrivially, and the output has degree
$d+a$.

Let $k\geq1$ and let $z_1,\ldots,z_k$ be independent formal variables.
For $Z=(z_1,\ldots,z_k)$, write
$p_m(Z)=\sum_{i=1}^kz_i^m$ and
$Z^{-1}=z_1^{-1}+\cdots+z_k^{-1}$.
Colons place all multiplication operators to the left of all
differentiations.  By \eqref{bridge:eq:half-commutation}, commuting
the second exponential in $V(z_i)$ past the first in $V(z_j)$ for
each $i<j$ gives
\begin{equation}\label{c1:eq:interaction-kernel}
 \begin{split}
 V(z_1)\cdots V(z_k)
 &=K_k(Z;q):V(z_1)\cdots V(z_k):,\\
 K_k(Z;q)&=\prod_{1\leq i<j\leq k}
            \frac{1-z_j/z_i}{1-qz_j/z_i}.
 \end{split}
\end{equation}
Each factor of $K_k$ is expanded in nonnegative powers of $z_j/z_i$.
The vertex-operator realization of Hall--Littlewood functions
\cite[Proposition~3.9, p.~243]{JingVertex}, conjugated by
$f[X]\mapsto f[-X]$, gives
\begin{equation}\label{c1:eq:jing-vacuum}
 V_{\lambda_1}\cdots V_{\lambda_{\ell(\lambda)}}1
 =Q_\lambda[-X;q],
 \qquad
 \mathcal R_\lambda
 =\left(\prod_{a\geq1}^{\mathrm{decr}}
          E_{a,m_a(\lambda)}\right)1.
\end{equation}
Here and below, $\mathrm{decr}$ means that the factors are written
from left to right in decreasing order of the indicated index;
factors with multiplicity zero are identities.  The second equality in
\eqref{c1:eq:jing-vacuum} uses
$b_\lambda(q)=(-1)^{\ell(\lambda)}
\prod_{a\geq1}B_{m_a(\lambda)}(q)$.

To extract a coefficient from the product in \eqref{c1:eq:interaction-kernel},
we use the finite Hall--Littlewood inner product.  For symmetric
Laurent polynomials $f,g$ in the $k$ variables $Z$, put
\begin{equation}\label{c1:eq:finite-inner-product}
 \begin{aligned}
 \langle f,g\rangle'_{q,k}
 &=\frac1{k!}\operatorname{CT}_Z f(Z)g(Z^{-1})\Delta_k(Z;q),\\
 \Delta_k(Z;q)&=\prod_{i\ne j}\frac{1-z_i/z_j}{1-qz_i/z_j}.
 \end{aligned}
\end{equation}
Here $\operatorname{CT}_Z$ denotes the constant term.  We expand each
factor $(1-qz_i/z_j)^{-1}$ in nonnegative powers of $q$ and take the
constant term in $Z$ coefficientwise.
For an integer $d\geq0$, write
\[
 [d]_q!=\prod_{i=1}^d\frac{q^i-1}{q-1}.
\]
For a partition $\lambda$ of length at most $k$, append zeros to
length $k$, let $m_0(\lambda)=k-\ell(\lambda)$, and define
\[
 v_\lambda^{(k)}(q)
 =\prod_{j\geq0}[m_j(\lambda)]_q!.
\]
In particular, the zero parts contribute to this product.
The finite Hall--Littlewood orthogonality relation is
\begin{equation}\label{c1:eq:finite-orthogonality}
 \langle P_\lambda,P_\mu\rangle'_{q,k}
 =\frac{\delta_{\lambda\mu}}{v_\lambda^{(k)}(q)}.
\end{equation}
Equation~\eqref{c1:eq:finite-orthogonality} is the Hall--Littlewood
specialization of
\cite[Chapter~VI, (9.3), (9.5)--(9.7), the remark on p.~372, and \S9, Example~1(c)--(d)]{Macdonald}:
set the two parameters there equal to $0$ and $q$, respectively.

For a partition $\beta$ padded with zeros to length $k$, write
$Z^\beta=z_1^{\beta_1}\cdots z_k^{\beta_k}$.  For a symmetric Laurent
polynomial $F(Z)$, its coefficient after multiplication by $K_k$ is
\begin{equation}\label{c1:eq:coefficient-pairing}
 \mathcal L_\beta(F)
 :=\Coef_{Z^\beta}\bigl[K_k(Z;q)F(Z)\bigr]
 =v_\beta^{(k)}(q)\langle F,P_\beta\rangle'_{q,k}.
\end{equation}
To see the second equality, write the coefficient as
$\operatorname{CT}_Z Z^{-\beta}K_k(Z;q)F(Z)$.  Permuting the variables
in $Z^{-\beta}K_k$ leaves this constant term unchanged, since $F$ is
symmetric.  We can therefore average over all $k!$ permutations.
For each $i<j$, the product for $\Delta_k$ in
\eqref{c1:eq:finite-inner-product} contains factors in both $z_j/z_i$
and $z_i/z_j$.  By \eqref{c1:eq:interaction-kernel}, the factors in
$z_j/z_i$ give $K_k$, so dividing $\Delta_k$ by the product of the
$z_i/z_j$ factors recovers $K_k$.
Since $\Delta_k$ is invariant under permutations, we can take it
outside the sum to obtain
\[
 \sum_{\pi\in\Sym_k}\pi\bigl(Z^{-\beta}K_k(Z;q)\bigr)
 =\Delta_k(Z;q)
   \sum_{\pi\in\Sym_k}\pi\left(
      Z^{-\beta}\prod_{i<j}\frac{1-qz_i/z_j}{1-z_i/z_j}\right).
\]
By the defining formula for $P_\beta$ evaluated at $Z^{-1}$
\cite[Chapter~III, (2.1)]{Macdonald}, the sum on the right equals
$v_\beta^{(k)}(q)P_\beta(Z^{-1};q)$.
Multiplying both sides by $F/k!$ and taking constant terms now gives
\eqref{c1:eq:coefficient-pairing} by \eqref{c1:eq:finite-inner-product}.
For the empty partition, $P_\varnothing=1$ and
$v_\varnothing^{(k)}(q)=[k]_q!$, so in particular
\begin{equation}\label{c1:eq:oriented-constant-term}
 \operatorname{CT}_Z K_k(Z;q)F(Z)
 =[k]_q!\,\langle F,1\rangle'_{q,k}.
\end{equation}
When $X$ also occurs, we take the constant term and the inner product
in $Z$ at each homogeneous degree in $X$.

Let $\rho$ be a partition of length at most $k$, padded with zeros
to length $k$.  For an integer $a\geq0$, adding $a$ to each part
multiplies every summand in the symmetrization formula used in the proof
of \eqref{c1:eq:coefficient-pairing} by $(z_1\cdots z_k)^a$.  Thus
\begin{equation}\label{eq:refaudit-hl-rectangular-shift}
 (z_1\cdots z_k)^aP_\rho(Z;q)=P_{\rho+(a^k)}(Z;q).
\end{equation}
To allow negative parts, let $\beta=(\beta_1,\ldots,\beta_k)$ be a
weakly decreasing integer sequence.  Choose an integer $N\geq0$ with
$\beta_k+N\geq0$ and define
\[
 P_\beta(Z;q)
 =(z_1\cdots z_k)^{-N}
   P_{\beta+(N^k)}(Z;q).
\]
Equation~\eqref{eq:refaudit-hl-rectangular-shift} makes this expression
independent of $N$.  A common shift by $N$ to nonnegative indices preserves
the multiplicities of equal parts, and the factors $(z_1\cdots z_k)^N$ and
$(z_1\cdots z_k)^{-N}$ cancel in \eqref{c1:eq:finite-inner-product}.
Thus \eqref{c1:eq:finite-orthogonality} also holds for these Laurent polynomials.

In the defining sum for $P_\rho$, substituting $Z^{-1}$
reverses the order of the two indices in each factor of the product over $i<j$.
Relabeling $z_i$ as $z_{k+1-i}$ restores the original order and gives the index
$(-\rho_k,\ldots,-\rho_1)$.  For any integer $a$, multiplication by
$(z_1\cdots z_k)^a$ adds $a$ to every part, giving
\begin{equation}\label{c1:eq:rectangle-inversion}
 (z_1\cdots z_k)^aP_\rho(Z^{-1};q)
 =P_{(a-\rho_k,\ldots,a-\rho_1)}(Z;q).
\end{equation}

We evaluate the inner products below by expanding both entries in the
Hall--Littlewood basis.  For the exponential in the first entry, the
Cauchy identity gives
\begin{equation}\label{c1:eq:hl-cauchy}
 \operatorname{Exp}[(q-1)XZ]
 =\sum_{\ell(\lambda)\leq k}
     Q_\lambda[-X;q]P_\lambda(Z;q).
\end{equation}
For the power sum in the second entry, define for every nonempty
partition $\rho$
\[
 n(\rho)=\sum_{i\geq1}(i-1)\rho_i,
 \qquad
 c_\rho(q)
 =q^{\,n(\rho)-\binom{\ell(\rho)}2}
     B_{\ell(\rho)-1}(q).
\]
For every $m\geq1$, \cite[Chapter~III, \S7, Example~2, p.~248]{Macdonald}
gives the power-sum expansion
\begin{equation}\label{c1:eq:power-sum-hl}
 p_m=\sum_{\rho\vdash m}c_\rho(q)P_\rho(X;q),
 \qquad c_\rho(q)\in\mathbb N[q-1].
\end{equation}
The nonnegativity in \eqref{c1:eq:power-sum-hl} follows from
$c_\rho(q)=q^{\,n(\rho)-\binom{\ell(\rho)}2}B_{\ell(\rho)-1}(q)$
and $n(\rho)\geq\binom{\ell(\rho)}2$.

\subsection{Positivity of normalized Hall--Littlewood modes}
\label{c1:subsec:operator-positivity}

We prove that $E_{a,k}$ preserves $\mathbb Q_{\geq0}[q-1][p_1,p_2,\ldots]$
for all $a,k\geq1$ (Theorem~\ref{c1:thm:mode-positivity}).
To this end, we first establish Lemmas~\ref{c1:lem:rectangle-lowering}
and~\ref{c1:lem:positive-commutation}.

For integers $a,k\geq1$ and $\rho\subseteq(a^k)$, append zeros to
$\rho$ to obtain $k$ parts and define its rectangular complement by
\[
 \rho^c=(a-\rho_k,\ldots,a-\rho_1),
\]
omitting zero parts.  Thus $\ell(\rho^c)=k-m_a(\rho)$.

\begin{lemma}\label{c1:lem:rectangle-lowering}
For every $a,k,m\geq1$,
\begin{equation}\label{c1:eq:rectangle-lowering}
 m\partial_{p_m}\mathcal R_{(a^k)}
 =\sum_{\substack{\rho\vdash m\\\rho\subseteq(a^k)}}
    \frac{(q^m-1)c_\rho(q)}{B_{m_a(\rho)}(q)}
    \mathcal R_{\rho^c}.
\end{equation}
Every coefficient on the right belongs to $\mathbb N[q-1]$.
\end{lemma}

\begin{proof}
We first pair the Cauchy identity \eqref{c1:eq:hl-cauchy} with
$P_\lambda(Z;q)$ for a partition $\lambda$ of length at most $k$.
By orthogonality \eqref{c1:eq:finite-orthogonality}, only the term
indexed by $\lambda$ survives.  Writing $s=k-\ell(\lambda)$ for the
number of zero parts after padding, we obtain
\begin{equation}\label{c1:eq:vacuum-norm-ratio}
 \frac{\bigl\langle\operatorname{Exp}[(q-1)XZ],P_\lambda(Z;q)
       \bigr\rangle'_{q,k}}{(q-1)^k}
 =\frac{Q_\lambda[-X;q]}{(q-1)^k v_\lambda^{(k)}(q)}
 =\frac{\mathcal R_\lambda}{B_s(q)}.
\end{equation}
For the last equality, \eqref{c1:eq:jing-vacuum} gives
$Q_\lambda[-X;q]=\mathcal R_\lambda\prod_{j\geq1}B_{m_j(\lambda)}(q)$.
Using $B_i(q)=(q-1)^i[i]_q!$, the factors for the positive parts
cancel against those in $(q-1)^k v_\lambda^{(k)}(q)$, leaving
$(q-1)^s[s]_q!=B_s(q)$ in the denominator.

Take $\lambda=(a^k)$ in \eqref{c1:eq:vacuum-norm-ratio}.
Then $B_0(q)=1$ and $P_{(a^k)}(Z;q)=(z_1\cdots z_k)^a$ by
\eqref{eq:refaudit-hl-rectangular-shift}.
In the ordinary exponential form of the Cauchy kernel, the only term
in the exponent involving $p_m(X)$ is $(q^m-1)p_m(X)p_m(Z)/m$.  Thus
\[
 m\frac{\partial}{\partial p_m(X)}\operatorname{Exp}[(q-1)XZ]
 =(q^m-1)p_m(Z)\,\operatorname{Exp}[(q-1)XZ].
\]
By \eqref{c1:eq:finite-inner-product}, moving
$p_m(Z)$ to the second entry replaces it by $p_m(Z^{-1})$, so
\begin{equation}\label{c1:eq:rectangle-finite-inner-product}
 m\partial_{p_m}\mathcal R_{(a^k)}
 =\frac{q^m-1}{(q-1)^k}
 \bigl\langle\operatorname{Exp}[(q-1)XZ],
       (z_1\cdots z_k)^a p_m(Z^{-1})\bigr\rangle'_{q,k}.
\end{equation}

Expand $p_m(Z^{-1})$ by \eqref{c1:eq:power-sum-hl}.  Terms with
$\ell(\rho)>k$ vanish in $k$ variables.  Padding the remaining $\rho$
with zeros, we obtain from \eqref{c1:eq:rectangle-inversion}
\[
 (z_1\cdots z_k)^a p_m(Z^{-1})
 =\sum_{\substack{\rho\vdash m\\\ell(\rho)\leq k}}
     c_\rho(q)P_{(a-\rho_k,\ldots,a-\rho_1)}(Z;q).
\]
If $\rho_1>a$, the last part $a-\rho_1$ is negative, so by
\eqref{c1:eq:finite-orthogonality} for Laurent polynomials this term
pairs to zero with every $P_\nu$ in \eqref{c1:eq:hl-cauchy}.
Thus only $\rho\subseteq(a^k)$ contribute to
\eqref{c1:eq:rectangle-finite-inner-product}, and the sum is empty if
$m>ak$.  For each such $\rho$, the index is $\rho^c$.  Its padded form
has one zero for each part of $\rho$ equal to $a$, so the denominator
in \eqref{c1:eq:vacuum-norm-ratio} for $\lambda=\rho^c$ is
$B_{m_a(\rho)}(q)$.  This gives \eqref{c1:eq:rectangle-lowering}.

For positivity, put $\ell=\ell(\rho)$.  The coefficient
in \eqref{c1:eq:power-sum-hl} is $c_\rho(q)=q^{n(\rho)-\binom{\ell}{2}}B_{\ell-1}(q)$
with nonnegative exponent.  If $m_a(\rho)<\ell$, the quotient $B_{\ell-1}(q)/B_{m_a(\rho)}(q)$
is a product of factors $q^j-1$, each in $\mathbb N[q-1]$, as is $q^m-1$.
If $m_a(\rho)=\ell$, then $\rho=(a^\ell)$ and $m=a\ell$.  Cancelling $B_{\ell-1}(q)$ leaves,
apart from the power of $q$, the geometric sum
\[
 \frac{q^{a\ell}-1}{q^\ell-1}=1+q^\ell+\cdots+q^{(a-1)\ell}.
\]
In both cases the coefficient belongs to $\mathbb N[q-1]$.
\end{proof}

We next commute $E_{a,k}$ past multiplication by one power sum
$p_m=p_m(X)$.  In the resulting products, the factors are ordered
by decreasing mode index as in \eqref{c1:eq:jing-vacuum}.

\begin{lemma}\label{c1:lem:positive-commutation}
For every $a,k,m\geq1$,
\begin{equation}\label{c1:eq:positive-commutation}
 [E_{a,k},p_m]
 =\sum_{\substack{\rho\vdash m\\\ell(\rho)\leq k}}
 c_\rho(q)
 \left(\prod_{j\geq1}^{\mathrm{decr}}E_{a+j,m_j(\rho)}\right)
 E_{a,k-\ell(\rho)}.
\end{equation}
\end{lemma}

\begin{proof}
Write $F=:V(z_1)\cdots V(z_k):$.  By \eqref{c1:eq:interaction-kernel},
$K_k(Z;q)F=V(z_1)\cdots V(z_k)$.  Taking the coefficient of
$z_1^a\cdots z_k^a$ selects $V_a$ from each factor $V(z_i)$, giving
$V_a^k$.  Since $p_m(X)$ is independent of $Z$, we may take the
commutator with $p_m$ before extracting this coefficient.  Moreover,
$[K_k(Z;q)F,p_m]=K_k(Z;q)[F,p_m]$ because $K_k(Z;q)$ is scalar.
Dividing by $B_k(q)$ as in \eqref{c1:eq:normalized-modes} gives
\begin{equation}\label{c1:eq:rectangle-commutator-extraction}
 [E_{a,k},p_m]
 =\frac{[V_a^k,p_m]}{B_k(q)}
 =\frac{1}{B_k(q)}\Coef_{z_1^a\cdots z_k^a}
      \bigl[K_k(Z;q)[F,p_m]\bigr].
\end{equation}
By \eqref{c1:eq:jing-field}, the normal-ordered product is
\begin{equation}\label{eq:refaudit-jing-normal-product}
 :V(z_1)\cdots V(z_k):
 =\operatorname{Exp}[(q-1)XZ]\,
   \exp\left(\sum_{j\geq1}p_j(Z^{-1})\partial_{p_j}\right).
\end{equation}
The second exponential replaces $p_m(X)$ by $p_m(X)+p_m(Z^{-1})$,
and the first commutes with multiplication by $p_m(X)$.  Hence
\begin{equation}\label{eq:refaudit-jing-normal-commutator}
 [F,p_m]=p_m(Z^{-1})F.
\end{equation}

Substituting \eqref{eq:refaudit-jing-normal-commutator} into
\eqref{c1:eq:rectangle-commutator-extraction} gives the coefficient of
$z_1^a\cdots z_k^a$ in $K_k(Z;q)p_m(Z^{-1})F$, divided by $B_k(q)$.
For fixed input and output degrees, the degrees of differentiation and
multiplication in \eqref{eq:refaudit-jing-normal-product} are bounded,
so each matrix entry of $F$ is a symmetric Laurent polynomial in $Z$.
We can therefore apply \eqref{c1:eq:coefficient-pairing} entrywise to obtain
\[
 [E_{a,k},p_m]
 =\frac{v_{(a^k)}^{(k)}(q)}{B_k(q)}
   \bigl\langle p_m(Z^{-1})F,P_{(a^k)}\bigr\rangle'_{q,k}
 =\frac{1}{(q-1)^k}
   \bigl\langle F,(z_1\cdots z_k)^a p_m(Z)\bigr\rangle'_{q,k}.
\]
Here $P_{(a^k)}(Z;q)=(z_1\cdots z_k)^a$ by
\eqref{eq:refaudit-hl-rectangular-shift}, and
$v_{(a^k)}^{(k)}(q)/B_k(q)=[k]_q!/B_k(q)=(q-1)^{-k}$.
The second entry contains $p_m(Z)$ because
\eqref{c1:eq:finite-inner-product} evaluates that entry at $Z^{-1}$.

Expand $p_m(Z)$ by \eqref{c1:eq:power-sum-hl} and apply
\eqref{eq:refaudit-hl-rectangular-shift} to each term.  In $k$ variables this gives
\[
 (z_1\cdots z_k)^a p_m(Z)
 =\sum_{\substack{\rho\vdash m\\\ell(\rho)\leq k}}
    c_\rho(q)P_{\rho+(a^k)}(Z;q).
\]
For each $\rho$, append zeros to length $k$ and put $\beta=\rho+(a^k)$.
By \eqref{c1:eq:coefficient-pairing} and \eqref{c1:eq:interaction-kernel},
\[
 \frac{\langle F,P_\beta\rangle'_{q,k}}{(q-1)^k}
 =\frac{\Coef_{Z^\beta}\bigl[K_k(Z;q)F\bigr]}{(q-1)^k v_\beta^{(k)}(q)}
 =\frac{V_{\beta_1}\cdots V_{\beta_k}}
        {(q-1)^k v_\beta^{(k)}(q)}.
\]
The parts of $\rho$ equal to $j>0$ give $m_j(\rho)$ copies of $a+j$
in $\beta$, and the appended zeros give $k-\ell(\rho)$ copies of $a$.
Thus $B_i(q)=(q-1)^i[i]_q!$ gives
\[
 (q-1)^k v_\beta^{(k)}(q)
 =B_{k-\ell(\rho)}(q)\prod_{j\geq1}B_{m_j(\rho)}(q).
\]
The modes are already in decreasing order.  Grouping adjacent equal
modes with their respective denominators gives
\eqref{c1:eq:positive-commutation} by \eqref{c1:eq:normalized-modes}.
\end{proof}

The coefficients in \eqref{c1:eq:positive-commutation} belong to
$\mathbb N[q-1]$ by \eqref{c1:eq:power-sum-hl}.  We now use this
commutator to pass from the value of $E_{a,k}$ on $1$ to its value on
any power-sum monomial.  Write
$\mathbb Q_{\geq0}[q-1][p_1,p_2,\ldots]$ for the cone of polynomials
in the power sums whose coefficients are polynomials in $q-1$ with
nonnegative rational coefficients.

\begin{theorem}[Positivity of operator coefficients]
\label{c1:thm:mode-positivity}
For every $a,k\geq1$, the operator $E_{a,k}$ preserves
$\mathbb Q_{\geq0}[q-1][p_1,p_2,\ldots]$.
\end{theorem}

\begin{proof}
By linearity it suffices to show that $E_{a,k}p_\mu$ belongs to the
cone for every partition $\mu$.  We induct on the output degree
$D=ak+|\mu|$.  At each $D$, we first treat the constant input for
all pairs $a,k$ with $ak=D$.  For nonconstant inputs, commuting past
one power sum by \eqref{c1:eq:positive-commutation} produces terms
whose last mode has output degree $D$ and index greater than $a$.
We therefore treat these inputs in decreasing order of $a$.

For $ak=D$, equation \eqref{c1:eq:jing-vacuum} gives
$E_{a,k}1=\mathcal R_{(a^k)}$.  For each $1\leq m\leq D$,
Lemma~\ref{c1:lem:rectangle-lowering} expresses
$m\partial_{p_m}\mathcal R_{(a^k)}$ as an $\mathbb N[q-1]$-linear
combination of functions $\mathcal R_{\rho^c}$ of degree $D-m$.
The empty complement gives $\mathcal R_\varnothing=1$.  For every
other complement, \eqref{c1:eq:jing-vacuum} constructs
$\mathcal R_{\rho^c}$ from $1$ by successive degree-raising modes.
Every intermediate output has degree at most $D-m<D$, so induction
puts each $\mathcal R_{\rho^c}$, and hence each derivative, in the cone.
Since $\deg p_m=m$, the weighted Euler identity gives
\begin{equation}\label{c1:eq:rectangle-euler}
 D\mathcal R_{(a^k)}
 =\sum_{m=1}^{D}m p_m\partial_{p_m}\mathcal R_{(a^k)}.
\end{equation}
Multiplication by $p_m$ and division by $D>0$ preserve the cone.
This proves the constant-input case and starts the induction at $D=1$,
where the only complement is empty.

If $\mu\ne\varnothing$, remove one part $m$ to obtain $\tau$.
Commuting $E_{a,k}$ past $p_m$ by \eqref{c1:eq:positive-commutation}
gives
\begin{equation}\label{c1:eq:mode-positivity-recursion}
 E_{a,k}p_\mu
 =p_mE_{a,k}p_\tau+\sum_{\substack{\rho\vdash m\\\ell(\rho)\leq k}}
 c_\rho(q)
 \left(\prod_{j\geq1}^{\mathrm{decr}}E_{a+j,m_j(\rho)}\right)
 E_{a,k-\ell(\rho)}p_\tau .
\end{equation}
In the first term, $E_{a,k}p_\tau$ has output degree $D-m$, so it
belongs to the cone by induction, as does its product with $p_m$.
In the summand indexed by $\rho$, pad $\rho$ to length $k$.
The mode indices $a+\rho_1,\ldots,a+\rho_k$ sum to $ak+m$, so the
final output on $p_\tau$ has degree $D$.  Read the product from right
to left, omitting identity factors.  Each step before the last has
output degree less than $D$ and hence preserves positivity by induction,
applied termwise in the power-sum basis.  The last mode is
$E_{a+\rho_1,m_{\rho_1}(\rho)}$.  Its output has degree $D$, and
the case with index $a+\rho_1>a$ has already been proved at this
degree.  Since $c_\rho(q)\in\mathbb N[q-1]$ by
\eqref{c1:eq:power-sum-hl}, every summand in
\eqref{c1:eq:mode-positivity-recursion} belongs to the cone.
\end{proof}

\subsection{Normal-order coefficients}
\label{c1:subsec:normal-coefficients}

We use the normal-ordered product \eqref{eq:refaudit-jing-normal-product}
to read off the coefficient of each differential operator in $E_{a,k}$.
Corollary~\ref{c1:cor:normal-coefficients} gives both
the power-sum positivity and the integrality needed in
Section~\ref{sec:c1-positivity}.

Fix integers $a,k\geq1$.  For a symmetric polynomial $f$ in the $k$ variables
$Z$ over $\mathbb Q(q)$, define
\begin{equation}\label{c1:eq:normal-coefficient-map}
 \mathcal S_{a,k}(f)
 =\frac1{B_k(q)}\Coef_{z_1^a\cdots z_k^a}
     \bigl[K_k(Z;q)\operatorname{Exp}[(q-1)XZ]f(Z^{-1})\bigr].
\end{equation}
When $f$ is a symmetric function, we use its specialization to $Z$.
For $f=p_\sigma$, the value $\mathcal S_{a,k}(p_\sigma)$ is the
multiplication factor in front of $p_\sigma^\perp/z_\sigma$, where
$p_\sigma^\perp$ is the ordinary Hall adjoint of multiplication by $p_\sigma(X)$.
Indeed, expanding the differentiation exponential in
\eqref{eq:refaudit-jing-normal-product} and using
$p_j^\perp=j\partial_{p_j}$ gives
\[
 \exp\left(\sum_{j\geq1}p_j(Z^{-1})\partial_{p_j}\right)
 =\sum_{\sigma\in\mathscr P}
     p_\sigma(Z^{-1})\frac{p_\sigma^\perp}{z_\sigma}.
\]
Choosing $m_j(\sigma)$ copies of $\partial_{p_j}$ contributes the
denominator $m_j(\sigma)!$, accounting for the factorials in $z_\sigma$.
Multiplying by $K_k(Z;q)\operatorname{Exp}[(q-1)XZ]$ and extracting
the coefficient in \eqref{c1:eq:normal-coefficient-map} gives
\begin{equation}\label{c1:eq:normal-expansion}
 E_{a,k}
 =\sum_{\sigma\in\mathscr P}
    \mathcal S_{a,k}(p_\sigma)\frac{p_\sigma^\perp}{z_\sigma}.
\end{equation}
Only terms with $|\sigma|\leq d$ contribute on an input of degree $d$.

For use in Section~\ref{sec:c1-positivity}, we now consider
$\mathcal S_{a,k}(p_\sigma h_b)$.

\begin{corollary}\label{c1:cor:normal-coefficients}
For every $a,k\geq1$, every partition $\sigma$, and every integer $b\geq0$,
\begin{equation}\label{c1:eq:normal-coeff-positive}
 \mathcal S_{a,k}(p_\sigma h_b)
 \in\Lambda[q]\cap\mathbb Q_{\geq0}[q-1][p_1,p_2,\ldots].
\end{equation}
\end{corollary}

\begin{proof}
We first extract $\mathcal S_{a,k}(p_\sigma)$ from
\eqref{c1:eq:normal-expansion}.  For $m\geq1$, write
$\delta_m(A)=[A,p_m]$, with $p_m$ acting by multiplication.
The multiplication coefficients commute with $p_m$.  For the
differential factor indexed by $\tau$ we use
\[
 \frac{p_\tau^\perp}{z_\tau}
 =\prod_{j\geq1}\frac{\partial_{p_j}^{\,m_j(\tau)}}{m_j(\tau)!},
 \qquad
 \left[\frac{\partial_{p_m}^{\,u}}{u!},p_m\right]
 =\frac{\partial_{p_m}^{\,u-1}}{(u-1)!}\quad(u\geq1).
\]
Thus $\delta_m$ removes one occurrence of $\partial_{p_m}$ if
$m$ is a part of $\tau$ and gives zero otherwise.
Take one commutator for each part of $\sigma$ and apply the result
to $1$.  A missing derivative makes the commutator zero, and any
remaining derivative kills $1$, so only $\tau=\sigma$ survives.
For $\tau=\sigma$ the differential factor becomes the identity, giving
\begin{equation}\label{c1:eq:iterated-commutator}
 \mathcal S_{a,k}(p_\sigma)
 =\delta_{\sigma_1}\cdots
   \delta_{\sigma_{\ell(\sigma)}}(E_{a,k})1.
\end{equation}
For $\sigma=\varnothing$ the right side is $E_{a,k}1$.

By \eqref{c1:eq:positive-commutation}, $\delta_m(E_{a,k})$ is a sum
of products of normalized modes with coefficients in $\mathbb N[q-1]$.
Since $\delta_m(I)=0$ and $\delta_m(AB)=A\delta_m(B)+\delta_m(A)B$,
each further commutator gives a sum of the same kind.
Theorem~\ref{c1:thm:mode-positivity} shows that each such product
maps $1$ into the power-sum cone.  Thus \eqref{c1:eq:iterated-commutator}
proves positivity for $p_\sigma$.  For general $b$, expand
$h_b=\sum_{\tau\vdash b}p_\tau/z_\tau$ and use linearity.  Each
$p_\sigma p_\tau$ is again a power-sum monomial, with positive
weight $1/z_\tau$ in this expansion.

For integrality, we compute $\mathcal S_{a,k}$ on the Hall--Littlewood
basis.  Applying \eqref{c1:eq:coefficient-pairing} with index $(a^k)$
to \eqref{c1:eq:normal-coefficient-map} gives
\begin{equation}\label{eq:refaudit-normal-map-inner-product}
 \mathcal S_{a,k}(f)=\frac1{(q-1)^k}
   \bigl\langle\operatorname{Exp}[(q-1)XZ],
                   (z_1\cdots z_k)^a f(Z)\bigr\rangle'_{q,k}.
\end{equation}
Here $P_{(a^k)}=(z_1\cdots z_k)^a$,
$v_{(a^k)}^{(k)}(q)/B_k(q)=(q-1)^{-k}$, and
\eqref{c1:eq:finite-inner-product} evaluates the second entry at $Z^{-1}$.
Take $f=P_\rho$ with $\ell(\rho)\leq k$ and pad $\rho$ to length $k$.
By \eqref{eq:refaudit-hl-rectangular-shift}, the second entry becomes
$P_{\rho+(a^k)}$.  All $k$ parts of $\rho+(a^k)$ are positive, so
\eqref{c1:eq:vacuum-norm-ratio} from the proof of
Lemma~\ref{c1:lem:rectangle-lowering} applies with $s=0$ and gives
\begin{equation}\label{c1:eq:normal-integral-map}
 \mathcal S_{a,k}(P_\rho)=\mathcal R_{\rho+(a^k)}.
\end{equation}
The right side belongs to $\Lambda[q]$ since plethystic substitution
of $-X$ acts in degree $d$ by $(-1)^d\omega$.
Now $p_\sigma h_b\in\Lambda$ has coefficients in $\mathbb Z[q]$ in
the Hall--Littlewood basis of Subsection~\ref{c1:subsec:hall-littlewood}.
The terms with $\ell(\rho)>k$ vanish in $k$ variables, so linearity
and \eqref{c1:eq:normal-integral-map} give
$\mathcal S_{a,k}(p_\sigma h_b)\in\Lambda[q]$.
\end{proof}

We apply the case $b=0$ of Corollary~\ref{c1:cor:normal-coefficients}
in the quotient coordinates of Section~\ref{sec:universal-blocks}.
Fix an integer
$r\geq2$ and an $r$-core $\delta$, and extend
$\mathcal F_\delta$ and the isomorphism $\Psi_\delta$ in
\eqref{eq:untransposed-quotient-map} to $\mathbb Q(q)$.
On $\mathcal F_\delta$, define the vertex operator $\mathscr V^{(r)}(z;q)$ by
replacing multiplication by $p_m$ and differentiation by
$\partial_{p_m}$ in \eqref{c1:eq:jing-field} with multiplication
by $p_{rm}$ and the operator $p_{rm}^{\perp}/(rm)$, respectively.
For integers $a,k\geq1$, its normalized repeated mode is
\[
 \mathscr E^{(r)}_{a,k}
 =\frac{\bigl(\Coef_{z^a}[\mathscr V^{(r)}(z;q)]\bigr)^k}{B_k(q)}.
\]
Writing $p_m^{(i)}=p_m(X^{(i)})$ for the power sum in the $i$th
quotient alphabet, consider the cone
\begin{equation}\label{bridge:eq:core-cone}
 \Psi_\delta^{-1}\!\left(
   \mathbb Q_{\geq0}[q-1][p_m^{(i)}:m\geq1,\ 0\leq i<r]\right).
\end{equation}

\begin{corollary}\label{bridge:cor:fixed-core}
For all $a,k\geq1$, the operator $\mathscr E^{(r)}_{a,k}$
preserves the cone \eqref{bridge:eq:core-cone}.
\end{corollary}

\begin{proof}
By \eqref{eq:untransposed-heisenberg-intertwining} and
\eqref{eq:untransposed-heisenberg-adjoint}, $\Psi_\delta$ carries
$p_{rm}$ to multiplication by $P_m=\sum_{i=0}^{r-1}p_m^{(i)}$ and
$p_{rm}^{\perp}/(rm)$ to the average
$d_m=r^{-1}\sum_{i=0}^{r-1}\partial/\partial p_m^{(i)}$.
Each of the $r$ partial derivatives sends $P_n$ to $\delta_{mn}$,
so $d_m(P_n)=\delta_{mn}$.  By the product rule, moving $d_m$ past
multiplication by $P_n$ gives the same correction as moving
$\partial_{p_m}$ past $p_n$.  Since the $P_m$ commute among themselves,
as do the $d_m$, normal ordering gives \eqref{c1:eq:normal-expansion}
for $\Psi_\delta\mathscr E^{(r)}_{a,k}\Psi_\delta^{-1}$ with $p_m$
replaced by $P_m$ and $\partial_{p_m}$ by $d_m$.

For the term indexed by $\sigma$, Corollary~\ref{c1:cor:normal-coefficients}
with $b=0$ gives
$\mathcal S_{a,k}(p_\sigma)\in\mathbb Q_{\geq0}[q-1][p_1,p_2,\ldots]$.
Substituting $P_m=\sum_i p_m^{(i)}$ for each $p_m$ therefore gives a
multiplication factor with nonnegative coefficients in the quotient power
sums and $q-1$.  The differential factor is
$\prod_{j\geq1}d_j^{m_j(\sigma)}/m_j(\sigma)!$.  Each partial derivative
sends a power-sum monomial to zero or to a positive integer multiple of
another monomial, so averaging and dividing by the factorials preserve
nonnegativity.  On an input of total quotient degree $w$, only terms
with $|\sigma|\leq w$ contribute.  Summing these terms proves that
$\mathscr E^{(r)}_{a,k}$ preserves the cone \eqref{bridge:eq:core-cone}.
\end{proof}

The coefficient of $z^a$ raises the total quotient degree by $a$,
so $\mathscr E^{(r)}_{a,k}$ raises the $r$-weight by $ak$ and the
ordinary degree by $rak$.

In Section~\ref{sec:c1-positivity}, the Macdonald mode formula
\eqref{c1:eq:full-modes} takes the coefficient of $z_1\cdots z_k$,
so the proof of Theorem~\ref{c1:thm:conjecture-one} uses
Corollary~\ref{c1:cor:normal-coefficients} with $a=1$.

\section{Macdonald creation operators and character positivity}
\label{sec:c1-positivity}

In Subsection~\ref{c1:subsec:macdonald-modes} we express the Macdonald
creation modes in terms of the Hall--Littlewood kernel $K_k(Z;q)$
from \eqref{c1:eq:interaction-kernel}.
Applying this formula to power-sum monomials in
Subsection~\ref{c1:subsec:conjecture-one} gives the coefficients
studied in Subsection~\ref{c1:subsec:normal-coefficients}.
Their positivity and integrality prove Theorem~\ref{thm:intro-c1},
from which we deduce a positivity result for Macdonald characters.

\subsection{Modes of the Macdonald creation operators}
\label{c1:subsec:macdonald-modes}

Proposition~\ref{c1:prop:full-modes} expresses each Macdonald creation
mode as a sum of normal-ordered products.  We work over $\mathbb Q(q,t)$
with $q,t$ algebraically independent.
Write $J_\lambda^{(q,t)}=J_\lambda(X;q,t)$ for the integral Macdonald
symmetric function used in the introduction.  Let $c$ be an auxiliary
formal parameter for the diagonal operator $\boldsymbol\Delta_c$.
For a node $(i,j)$, set
\[
 \chi(i,j)=q^{j-1}t^{-(i-1)},\qquad
 \boldsymbol\Delta_c J_\lambda^{(q,t)}
 =\prod_{\square\in\lambda}(1-c\chi(\square))J_\lambda^{(q,t)}.
\]
For a formal alphabet $A$, recall the plethystic exponential
$\operatorname{Exp}[A]=\exp\bigl(\sum_{m\geq1}p_m[A]/m\bigr)$,
where $\exp$ denotes the ordinary formal exponential.
Let $\mathcal P_A$ denote multiplication by $\operatorname{Exp}[AX]$.
We introduce the operator
\begin{equation}\label{c1:eq:gamma-definition}
 \boldsymbol\Gamma(u,v)
 =\boldsymbol\Delta_{1/v}
   \mathcal P_{uv(1-t)/(1-q)}
   \boldsymbol\Delta_{1/v}^{-1}.
\end{equation}
Fix a partition $\lambda$ and an integer $n\geq0$.  To see how the
$u^n$ coefficient specializes at $v=0$, we compute
$\Coef_{u^n}[\boldsymbol\Gamma(u,v)J_\lambda^{(q,t)}]$ first over
$\mathbb Q(q,t,v)$.  Only the middle factor in
\eqref{c1:eq:gamma-definition} depends on $u$.  It acts by multiplication by
\[
 \operatorname{Exp}[uvX(1-t)/(1-q)]
 =\sum_{j\geq0}u^jv^j h_j[X(1-t)/(1-q)],
\]
so taking the $u^n$ coefficient replaces this factor by multiplication
by $v^n h_n[X(1-t)/(1-q)]$.

Apply the three factors in \eqref{c1:eq:gamma-definition} from right to
left.  First, $\boldsymbol\Delta_{1/v}^{-1}$ divides $J_\lambda^{(q,t)}$
by $\prod_{\square\in\lambda}(1-\chi(\square)/v)$.  Next, expand
$h_n[X(1-t)/(1-q)]J_\lambda^{(q,t)}$ in the Macdonald basis and
multiply each coefficient by $v^n$.  For $n\geq1$, the Pieri rule
gives a sum over partitions $\xi$ obtained by adding a horizontal strip
of $n$ boxes to $\lambda$.
For $n=0$, we have $h_0=1$ and $\xi=\lambda$.
Finally, $\boldsymbol\Delta_{1/v}$ acts on each such $J_\xi^{(q,t)}$
by multiplication by $\prod_{\square\in\xi}(1-\chi(\square)/v)$.
Thus the coefficient of $J_\xi^{(q,t)}$ in the Pieri expansion acquires
the factor
\begin{equation}\label{eq:refaudit-gamma-polynomial-factor}
 v^n\frac{\prod_{\square\in\xi}(1-\chi(\square)/v)}
               {\prod_{\square\in\lambda}(1-\chi(\square)/v)}
 =v^n\prod_{\square\in\xi\setminus\lambda}(1-\chi(\square)/v)
 =\prod_{\square\in\xi\setminus\lambda}(v-\chi(\square)).
\end{equation}
The factors for the boxes of $\lambda$ cancel, and each of the $n$ new
boxes absorbs one factor of $v$.  Hence the $u^n$ coefficient raises
degree by $n$ and has matrix entries polynomial in $v$, so it can be
evaluated at $v=0$.

Let $R$ be the algebra automorphism given by
$Rp_m=(1-t^m)(1-q^m)^{-1}p_m$, and put
$\widehat\Delta_c=R^{-1}\boldsymbol\Delta_cR$ and
$\widehat\Gamma(u,v)=R^{-1}\boldsymbol\Gamma(u,v)R$.
Let $H(s)$ denote multiplication by $\sum_{j\geq0}h_j[X]s^j$.
Then
\begin{equation}\label{c1:eq:gamma-conjugation}
 \widehat\Gamma(u,v)
 =\widehat\Delta_{1/v}H(uv)\widehat\Delta_{1/v}^{-1}.
\end{equation}
Define
\begin{equation}\label{c1:eq:G-definition}
 \begin{aligned}
 G_n(v)&=(-1)^n\Coef_{u^n}\bigl[\widehat\Gamma(u,v)\bigr],\\
 G(u,v)&=\sum_{n\geq0}G_n(v)u^n
 =\widehat\Gamma(-u,v)
 =\widehat\Delta_{1/v}H(-uv)
   \widehat\Delta_{1/v}^{-1}.
 \end{aligned}
\end{equation}
In particular, $G_0(v)=I$.
Both $R$ and $R^{-1}$ are self-adjoint for the ordinary Hall inner product,
since they are diagonal in the orthogonal power-sum basis.  Thus,
for $n\geq0$ and $|\nu|=|\mu|+n$,
\begin{equation}\label{c1:eq:original-pairing}
 \begin{aligned}
 (-1)^n t^{|\mu|}
 \left\langle \Coef_{u^n}\bigl[\boldsymbol\Gamma(u,v)Rp_\mu\bigr],
              R^{-1}p_\nu\right\rangle
 &=(-1)^n t^{|\mu|}
   \left\langle \Coef_{u^n}\bigl[R^{-1}\boldsymbol\Gamma(u,v)Rp_\mu\bigr],
                p_\nu\right\rangle\\
 &=t^{|\mu|}z_\nu\Coef_{p_\nu}\bigl[G_n(v)p_\mu\bigr].
 \end{aligned}
\end{equation}
The last equality uses \eqref{c1:eq:G-definition} and
$\langle F,p_\nu\rangle=z_\nu\Coef_{p_\nu}[F]$.
The left-hand side of \eqref{c1:eq:original-pairing}
is the expression in Conjecture~1 of \cite[Section~4.5]{BD26}.
For $n<0$ the coefficient vanishes since $\boldsymbol\Gamma(u,v)$ has
only nonnegative powers of $u$.

Introduce the vertex operator
\begin{equation}\label{c1:eq:eta-definition}
 \eta(z)=
 \exp\left(\sum_{m\geq1}\frac{q^m-1}{mt^m}p_mz^m\right)
 \exp\left(\sum_{m\geq1}(t^m-1)\partial_{p_m}z^{-m}\right).
\end{equation}
For $Z=z_1+\cdots+z_k$, write $Z^{-1}=z_1^{-1}+\cdots+z_k^{-1}$
and use the kernel $K_k(Z;q)$ of
\eqref{c1:eq:interaction-kernel}, expanded in nonnegative powers of
$z_j/z_i$ for $i<j$.  Put
\[
 c_k(q)=\frac{q^{\binom{k}{2}}}{B_k(q)},\qquad c_0(q)=1.
\]
The operator $\Coef_{z_1\cdots z_k}$ extracts exponent one in
every variable, and $\operatorname{CT}_Z$ means extraction of their
constant term.  Proposition~\ref{c1:prop:full-modes} expresses $G_n(v)$
in terms of the kernel $K_k(Z;q)$ and normal-ordered products of the
operators $\eta(z_i)$ with complete symmetric functions.  The colons
place all multiplication factors to the left of the derivatives, so the
derivatives act only on the input symmetric function.

\begin{samepage}
\begin{proposition}\label{c1:prop:full-modes}
For every $n\geq0$,
\begin{equation}\label{c1:eq:full-modes}
 \begin{aligned}
 G_n(v)
 =\sum_{k=0}^n(-v)^{n-k}t^nc_k(q)
 \Coef_{z_1\cdots z_k}\bigl[K_k(Z;q)\,
 :\eta(z_1)\cdots\eta(z_k)
           h_{n-k}[X/t+Z^{-1}]:\bigr].
 \end{aligned}
\end{equation}
\end{proposition}
\end{samepage}

\begin{proof}
We first express $G(u,v)$ in terms of $G(u,0)$.  To compute $G_n(0)$,
we start from a known constant-term formula and extract the
coefficient of $z_1\cdots z_n$.  We then recover the coefficient of
$u^n$ in $G(u,v)$ using the generating function for the complete
symmetric functions.

The reduction to $v=0$ uses the five-term relation in the modified
Macdonald basis with parameters $q,t^{-1}$.  Let $U$ be the algebra
automorphism given by $Up_m=(1-q^m)p_m$.  Since
$\widetilde H_\lambda^{(q,t^{-1})}
 =t^{-n(\lambda)}J_\lambda^{(q,t)}[X/(1-t)]$, where
$n(\lambda)=\sum_{i\geq1}(i-1)\lambda_i$, we have
\begin{equation}\label{c1:eq:modified-basis}
 U\widetilde H_\lambda^{(q,t^{-1})}
       =t^{-n(\lambda)}R^{-1}J_\lambda^{(q,t)}.
\end{equation}
By \eqref{c1:eq:modified-basis}, $U$ carries the modified $\Delta_c$ to $\widehat\Delta_c$.
The nabla operator acts on the modified Macdonald basis by
\[
 \nabla_{q,t^{-1}}\widetilde H_\lambda^{(q,t^{-1})}
 =(-1)^{|\lambda|}\Bigl(\prod_{\square\in\lambda}\chi(\square)\Bigr)
     \widetilde H_\lambda^{(q,t^{-1})}.
\]
Let $\widehat N$ act on $R^{-1}J_\lambda^{(q,t)}$ by the eigenvalue
$\prod_{\square\in\lambda}\chi(\square)$, and let $\sigma$ act by
$(-1)^d$ on degree $d$.  By \eqref{c1:eq:modified-basis},
$U\nabla_{q,t^{-1}}U^{-1}$ and $\sigma\widehat N$ have the same eigenvalue
on each $R^{-1}J_\lambda^{(q,t)}$, so $U\nabla_{q,t^{-1}}U^{-1}=\sigma\widehat N$.
Write $P(s)$ for multiplication by $\operatorname{Exp}[sX/(1-t^{-1})]$ and
put $F(s)=\widehat N P(s)\widehat N^{-1}$.  Conjugation by $U$
takes multiplication by $\operatorname{Exp}[sX/((1-q)(1-t^{-1}))]$ to $P(s)$.
We use the Garsia--Mellit five-term relation \cite[Theorem~1.1]{GM19}
in the form given in \cite[equation~(11)]{BD26}, taking $u=cs$ and $v=1/c$ there.
Conjugating by $U$ and using
$\sigma P(cs)\sigma^{-1}=P(-cs)$, with ordinary negation of the
parameter, gives
\begin{equation}\label{c1:eq:five-term}
 \widehat\Delta_cP(s)\widehat\Delta_c^{-1}=F(-cs)P(s).
\end{equation}

Set $s=-uv$.  The product $P(s/t)^{-1}P(s)$ is $H(s)$, since its
logarithm is multiplication by $\sum_{m\geq1}p_m s^m/m$.
In \eqref{c1:eq:G-definition}, conjugate
these two factors separately using \eqref{c1:eq:five-term} with
$c=1/v$.  Since $-cs=u$, their product is
\[
 G(u,v)=(F(u/t)P(s/t))^{-1}F(u)P(s)
       =P(s/t)^{-1}F(u/t)^{-1}F(u)P(s).
\]
By \eqref{eq:refaudit-gamma-polynomial-factor}, each coefficient of
$G(u,v)$ is polynomial in $v$, so we can set $v=0$ here.  Both outer
factors then equal $I$.  To restore them, write
$C(a)=\mathcal P_{a/(1-t)}=P(a/t)^{-1}$.  Since
$P(s/t)^{-1}=H(s)C(ts)$ and $P(s)=C(ts)^{-1}$, we obtain
\begin{equation}\label{c1:eq:recover-modes}
 \begin{aligned}
 G(u,0)&=F(u/t)^{-1}F(u),\\
 G(u,v)&=H(s)C(ts)G(u,0)C(ts)^{-1},\qquad s=-uv.
 \end{aligned}
\end{equation}

To compute $G_n(0)$, we introduce the operators $D_k$.
For $k\geq1$ and a symmetric function $f$, apply the normal-ordered product
of $\eta(z_1),\ldots,\eta(z_k)$ from \eqref{c1:eq:eta-definition} to $f$,
multiply by $c_k(q)K_k(Z;q)$, and take the constant term in $z_1,\ldots,z_k$.
Thus
\begin{equation}\label{c1:eq:zero-modes}
 D_k f=c_k(q)\operatorname{CT}_Z
 \bigl[K_k(Z;q)\,:\eta(z_1)\cdots\eta(z_k):\,f\bigr],
 \qquad D_0=I.
\end{equation}
For $k\geq1$, $D_k$ is the $k$th constant-term operator of \cite[Definition~3.10, (3.9)]{FHHSY}
conjugated by $R^{-1}(\,\cdot\,)R$.
Indeed, this conjugation multiplies $p_m$ by $(1-q^m)/(1-t^m)$
and $\partial_{p_m}$ by its reciprocal, taking the vertex operator
in their equation~(1.7) to $\eta(z)$.
The prefactor in the first normal-ordering formula in the proof of
their Proposition~3.25 (p.~095215-28) is
$(-1)^k q^{\binom{k}{2}}/\prod_{j=1}^k(1-q^j)=c_k(q)$.
The Macdonald eigenvectors in their Proposition~3.17 may be rescaled from
$P_\lambda(X;q,t)$ to $J_\lambda^{(q,t)}$, so $D_k$ is diagonal
on the conjugated basis $R^{-1}J_\lambda^{(q,t)}$.

Let $d_\lambda(u)$ be the eigenvalue of $D(u)=\sum_{k\geq0}D_ku^k$
on $R^{-1}J_\lambda^{(q,t)}$.  For $s_i^\lambda=q^{\lambda_i}t^{-i}$,
the eigenvalue formula \cite[equation~(3.18) and Proposition~3.17]{FHHSY}
gives
\begin{equation}\label{c1:eq:zero-eigenvalues}
 d_\lambda(u)=\exp\left(\sum_{m\geq1}
         \frac{1-t^m}{1-q^m}p_m(s^\lambda)\frac{u^m}{m}\right).
\end{equation}
Compare this with the eigenvalue at the empty partition.  Summing
over the boxes in each row gives the finite difference
\[
 p_m(s^\lambda)-p_m(s^\varnothing)
 =\sum_{i=1}^{\ell(\lambda)}(q^{m\lambda_i}-1)t^{-im}
 =(q^m-1)t^{-m}\sum_{\square\in\lambda}\chi(\square)^m.
\]
Substituting in \eqref{c1:eq:zero-eigenvalues} and exponentiating
identifies the ratio of eigenvalues with the eigenvalue of
$\widehat\Delta_{u/t}\widehat\Delta_u^{-1}$ on the same basis vector.
Thus
\begin{equation}\label{c1:eq:zero-generating}
 \frac{d_\lambda(u)}{d_\varnothing(u)}
 =\prod_{\square\in\lambda}\frac{1-u\chi(\square)/t}{1-u\chi(\square)},
 \qquad
 D(u)=d_\varnothing(u)\widehat\Delta_{u/t}\widehat\Delta_u^{-1}.
\end{equation}

To extract the coefficient of $z_1\cdots z_n$ from the constant-term
formula \eqref{c1:eq:zero-modes}, we conjugate $D(u)$ by $C(a)$.
Setting $s=a/t$ in \eqref{c1:eq:five-term} and using $C(a)=P(a/t)^{-1}$ gives
$C(a)^{-1}\widehat\Delta_c C(a)=F(-ca/t)^{-1}\widehat\Delta_c$.
Use this identity with $c=u/t$ for the first diagonal factor in
\eqref{c1:eq:zero-generating}.  For the second, take the inverse of
the identity at $c=u$, which places $F(-au/t)$ to the right of
$\widehat\Delta_u^{-1}$.  Since $d_\varnothing(u)$ is a scalar, we obtain
\begin{equation}\label{c1:eq:translated-zero-modes}
 \begin{aligned}
 C(a)^{-1}D(u)C(a)
 &=F(-au/t^2)^{-1}
   \bigl(d_\varnothing(u)\widehat\Delta_{u/t}\widehat\Delta_u^{-1}\bigr)F(-au/t)\\
 &=F(-au/t^2)^{-1}D(u)F(-au/t).
 \end{aligned}
\end{equation}

We compare the coefficients of $a^nu^n$ in \eqref{c1:eq:translated-zero-modes}.
On the left, the coefficient of $u^n$ is $C(a)^{-1}D_nC(a)$.
To find its coefficient of $a^n$, recall from \eqref{c1:eq:eta-definition}
that the differentiation part of $\eta(z)$ translates
$p_m$ by $(t^m-1)z^{-m}$.  Applying this translation to the
multiplication factor $C(a)=\mathcal P_{a/(1-t)}$ gives
\begin{equation}\label{c1:eq:eta-conjugation}
 C(a)^{-1}\eta(z)C(a)
 =\exp\left(-\sum_{m\geq1}\frac{(a/z)^m}{m}\right)\eta(z)
 =(1-a/z)\eta(z).
\end{equation}
Thus the $n$ vertex operators in \eqref{c1:eq:zero-modes} contribute
the factor $\prod_{i=1}^n(1-a/z_i)$.  To obtain $a^n$ we must choose
$-a/z_i$ from every factor.  This supplies $(z_1\cdots z_n)^{-1}$,
which shifts the $z_1\cdots z_n$ term to the constant term.
Including the sign, the coefficient of $a^nu^n$ on the left is
\begin{equation}\label{eq:refaudit-translated-mode-coefficient}
 (-1)^nc_n(q)\operatorname{CT}_Z
 (z_1\cdots z_n)^{-1}K_n(Z;q):\eta(z_1)\cdots\eta(z_n):.
\end{equation}

On the right of \eqref{c1:eq:translated-zero-modes}, the two $F$-factors
depend only on $au$.  If together they contribute $(au)^r$, then
choosing $D_j u^j$ from $D(u)$ gives $a^r u^{r+j}$.
To obtain $a^nu^n$ we must have $r=n$ and $j=0$, so only $D_0=I$
contributes.  By the first identity in \eqref{c1:eq:recover-modes},
the remaining product is
$F(-au/t^2)^{-1}F(-au/t)=G(-au/t,0)$.
Its term of degree $n$ in $au$ is $G_n(0)(-au/t)^n$, so the coefficient
on the right is $(-1)^nt^{-n}G_n(0)$.  Equating this with
\eqref{eq:refaudit-translated-mode-coefficient} and cancelling $(-1)^n$ gives
\begin{equation}\label{c1:eq:nabla-modes}
 G_n(0)=t^nc_n(q)\Coef_{z_1\cdots z_n}
       \bigl[K_n(Z;q):\eta(z_1)\cdots\eta(z_n):\bigr].
\end{equation}

We now return to $G(u,v)$.  Reversing the conjugation in
\eqref{c1:eq:eta-conjugation} multiplies $\eta(z)$ by $(1-a/z)^{-1}$.
Thus in \eqref{c1:eq:recover-modes}, where $a=ts=-uvt$, substituting
\eqref{c1:eq:nabla-modes} gives
\begin{equation}\label{c1:eq:all-mode-generating}
 \begin{aligned}
 G(u,v)
 =H(-uv)\sum_{k\geq0}(ut)^kc_k(q)
 \operatorname{CT}_Z &(z_1\cdots z_k)^{-1}K_k(Z;q)\,
 :\eta(z_1)\cdots\eta(z_k):\\[-2pt]
 &\hspace{-12mm}\cdot
       \prod_{i=1}^k(1+uvt/z_i)^{-1}.
 \end{aligned}
\end{equation}
To extract $G_n(v)$, take the coefficient of $u^n$ in
\eqref{c1:eq:all-mode-generating}.  The factor $(ut)^k$ restricts
the sum to $k\leq n$.  In the geometric series for
$\prod_i(1+uvt/z_i)^{-1}$, each inverse power of a $z_i$ contributes
one power of $u$, so at most $n-k$ such powers can occur in total.
On an input of degree $d$, the derivatives in \eqref{c1:eq:eta-definition}
supply at most $d$ further inverse powers.  Since $(z_1\cdots z_k)^{-1}$
has degree $-k$ and $K_k(Z;q)$ has total $Z$-degree zero, the constant
term uses at most $d+n$ positive powers from the multiplication part
of the vertex operators.  Once these terms are fixed, the exponents
needed from $K_k$ are fixed as well.  Giving $z_i$ weight $i$, each
ratio $z_j/z_i$ in \eqref{c1:eq:interaction-kernel} has positive weight
$j-i$, so only finitely many terms of the kernel can supply those exponents.

We now combine the two factors that supply the remaining powers of $u$.
By homogeneity, $h_j[X]=t^jh_j[X/t]$, while the product over $i$
is the generating function for $h_b[Z^{-1}]$.  Thus
\[
 \begin{aligned}
 \left(\sum_{j\geq0}(-uv)^j h_j[X]\right)\prod_{i=1}^k(1+uvt/z_i)^{-1}
 &=\left(\sum_{j\geq0}(-uvt)^j h_j[X/t]\right)
   \left(\sum_{b\geq0}(-uvt)^b h_b[Z^{-1}]\right)\\
 &=\sum_{r\geq0}(-uvt)^r h_r[X/t+Z^{-1}],
 \end{aligned}
\]
where the last equality uses
$\sum_{j+b=r}h_j[X/t]h_b[Z^{-1}]=h_r[X/t+Z^{-1}]$.
The $X$-dependent factor comes from $H(-uv)$, which acts to the left
of all derivatives.  Thus $h_r[X/t+Z^{-1}]$ enters the normal-ordered
product as a multiplication factor.

In the $k$th summand, the power of $u$ is $k+r$, so we take $r=n-k$.
The scalar factor is then
$(ut)^k(-uvt)^{n-k}=u^n(-v)^{n-k}t^n$.
Finally, the prefactor $(z_1\cdots z_k)^{-1}$ turns the constant term
into the coefficient of $z_1\cdots z_k$ in the remaining product.
Summing over $0\leq k\leq n$ gives \eqref{c1:eq:full-modes}.
\end{proof}

We remark that when $k=0$ in \eqref{c1:eq:full-modes}, the products are
empty and coefficient extraction is the identity.  The summand is multiplication by
$(-v)^n t^n h_n[X/t]=(-v)^n h_n[X]$.

\subsection{Proof of Theorem~\ref{thm:intro-c1}}
\label{c1:subsec:conjecture-one}

To apply Proposition~\ref{c1:prop:full-modes} to a power-sum monomial
$p_\mu$, we expand $h_{n-k}[X/t+Z^{-1}]$ in the two alphabets.
The resulting factors $h_b[Z^{-1}]$ combine with $p_\sigma[Z^{-1}]$,
where $\sigma$ records the parts of $\mu$ selected by the derivatives.
For $k\geq1$, collecting powers of $t$ then reduces the coefficient
extraction to $\mathcal S_{1,k}(p_\sigma h_b)$, to which
Corollary~\ref{c1:cor:normal-coefficients} applies.

For partitions $\sigma,\mu$, the notation $\sigma\subseteq\mu$ in
this subsection means inclusion of their multisets of parts, and
$\mu\setminus\sigma$ is the complementary multiset.  We write
\[
 \binom\mu\sigma
   =\prod_{j\geq1}\binom{m_j(\mu)}{m_j(\sigma)}.
\]
We use the ordinary Hall inner product of \eqref{c1:eq:original-pairing}.

\begin{theorem}[Macdonald coefficient positivity]
\label{c1:thm:conjecture-one}
For every $n\geq0$ and partitions $\mu,\nu$ with
$|\nu|=|\mu|+n$,
\[
 t^{|\mu|}z_\nu\Coef_{p_\nu}\bigl[G_n(v)p_\mu\bigr]
       \in\mathbb N[q-1,t-1,-v].
\]
By \eqref{c1:eq:original-pairing}, this proves
Conjecture~1 of \cite{BD26}.
\end{theorem}

\begin{proof}
We keep track of the powers of $t$ in the mode expansion, then prove
nonnegativity in the power-sum basis and integrality by Hall pairing.
Set $m=|\mu|$.  In the normal product in
\eqref{c1:eq:full-modes}, the annihilation part acts on $p_\mu$ by
substituting $p_j\mapsto p_j+(t^j-1)p_j[Z^{-1}]$.
Expanding this substitution gives
\[
 \sum_{\sigma\subseteq\mu}\binom\mu\sigma
       \prod_{s\in\sigma}(t^s-1)
       p_{\mu\setminus\sigma}[X]p_\sigma[Z^{-1}].
\]
Also,
\[
 h_{n-k}[X/t+Z^{-1}]
 =\sum_{b=0}^{n-k}t^{-(n-k-b)}
       h_{n-k-b}[X]h_b[Z^{-1}].
\]
For the term indexed by $\sigma$ and $b$, the factor
$p_\sigma[Z^{-1}]h_b[Z^{-1}]$ has total $Z$-degree
$-|\sigma|-b$, while $K_k(Z;q)$ has total $Z$-degree zero.
Extracting $z_1\cdots z_k$, of total degree $k$, therefore selects
degree $k+b+|\sigma|$ from the creation factor
$\operatorname{Exp}[(q-1)XZ/t]$.  Its component of that degree is
$t^{-(k+b+|\sigma|)}$ times the corresponding component of
$\operatorname{Exp}[(q-1)XZ]$.  Combining this with the prefactor
$t^{m+n}$ and the factor $t^{-(n-k-b)}$ gives
\[
 t^{m+n}t^{-(k+b+|\sigma|)}t^{-(n-k-b)}
       =t^{m-|\sigma|}.
\]
For $k\geq1$, the remaining coefficient extraction, including the factor
$B_k(q)^{-1}$ from $c_k(q)$, is $\mathcal S_{1,k}(p_\sigma h_b)$ by
\eqref{c1:eq:normal-coefficient-map}.  Writing the $k=0$ term separately
gives the finite identity
\begin{equation}\label{c1:eq:positive-expansion}
 \begin{aligned}
 t^mG_n(v)p_\mu
 ={}&(-v)^nt^mh_np_\mu\\
 &+\sum_{k=1}^n(-v)^{n-k}q^{\binom{k}{2}}
   \sum_{\sigma\subseteq\mu}\binom\mu\sigma
       t^{m-|\sigma|}\prod_{s\in\sigma}(t^s-1)
       p_{\mu\setminus\sigma}\\
 &\hspace{24mm}\cdot
       \sum_{b=0}^{n-k}h_{n-k-b}[X]\,
                    \mathcal S_{1,k}(p_\sigma h_b).
 \end{aligned}
\end{equation}
Since the parts of $\sigma$ are selected from those of $\mu$, the
remaining exponent $m-|\sigma|$ is nonnegative.  Thus the factor
$t^m$ clears all negative powers of $t$ in \eqref{c1:eq:positive-expansion}.

By Corollary~\ref{c1:cor:normal-coefficients} with $a=1$,
$\mathcal S_{1,k}(p_\sigma h_b)$ has nonnegative rational power-sum
coefficients in $q-1$.  Each $h_j=\sum_{\tau\vdash j}p_\tau/z_\tau$
also has nonnegative rational power-sum coefficients, while
$q^{\binom{k}{2}}$, $t^{m-|\sigma|}$, and $t^s-1$ have nonnegative
integer coefficients in $q-1$ or $t-1$, as appropriate.
Every summand in \eqref{c1:eq:positive-expansion} therefore has
nonnegative rational power-sum coefficients in $q-1,t-1,-v$.
Hence
\begin{equation}\label{eq:refaudit-macdonald-rational-positivity}
 t^mz_\nu\Coef_{p_\nu}\bigl[G_n(v)p_\mu\bigr]
       \in\mathbb Q_{\geq0}[q-1,t-1,-v].
\end{equation}

For integrality, Corollary~\ref{c1:cor:normal-coefficients} also gives
$\mathcal S_{1,k}(p_\sigma h_b)\in\Lambda[q]$.
Since $p_{\mu\setminus\sigma}$ and each $h_j$ also belong to
$\Lambda$, \eqref{c1:eq:positive-expansion} lies in
$\Lambda[q,t,v]$.  Both $t^mG_n(v)p_\mu$ and $p_\nu$
have Schur coefficients in $\mathbb Z[q,t,v]$, so orthonormality
of the Schur basis makes their Hall pairing integral.
Power-sum orthogonality gives
\begin{equation}\label{eq:refaudit-macdonald-integral-pairing}
 t^mz_\nu\Coef_{p_\nu}\bigl[G_n(v)p_\mu\bigr]
 =\left\langle t^mG_n(v)p_\mu,p_\nu\right\rangle
       \in\mathbb Z[q,t,v]
       =\mathbb Z[q-1,t-1,-v].
\end{equation}
Combining \eqref{eq:refaudit-macdonald-rational-positivity} and
\eqref{eq:refaudit-macdonald-integral-pairing} proves membership in
$\mathbb N[q-1,t-1,-v]$.
\end{proof}

We now apply Theorem~\ref{c1:thm:conjecture-one} to Macdonald
characters, whose definition uses the Macdonald inner product
\[
 \langle p_\rho,p_\tau\rangle_{q,t}
 =\delta_{\rho\tau}z_\rho
       \prod_i\frac{1-q^{\rho_i}}{1-t^{\rho_i}}.
\]
For a partition $\mu$ and $\ell\geq1$, we use the character of
\cite[Definition~1.4]{BD26}, with the shifts $t^{1-j}$ as in the
proof of Proposition~4.7 there:
\begin{equation}\label{c1:eq:character-definition}
 \begin{aligned}
 \widetilde{\boldsymbol\theta}_\mu^{(q,t)}(v_1,\ldots,v_\ell)
 =\bigl\langle p_\mu,
 \Coef_{u^{|\mu|}}\bigl[\boldsymbol\Gamma(u,v_1)
  \boldsymbol\Gamma(t^{-1}u,v_2)\cdots
  \boldsymbol\Gamma(t^{-(\ell-1)}u,v_\ell)1\bigr]
                 \bigr\rangle_{q,t}.
 \end{aligned}
\end{equation}
The coefficient extraction in \eqref{c1:eq:character-definition}
selects homogeneous degree $|\mu|$, as does pairing with $p_\mu$.

\begin{corollary}\label{c1:cor:character-positivity}
For every partition $\mu$ and every $\ell\geq1$,
\[
 (-1)^{|\mu|}t^{2(\ell-1)|\mu|}
   \widetilde{\boldsymbol\theta}_\mu^{(q,t)}(v_1,\ldots,v_\ell)
 \in\mathbb Q_{\geq0}[q-1,t-1,-v_1,\ldots,-v_\ell].
\]
\end{corollary}

\begin{proof}
Expand \eqref{c1:eq:character-definition} in the basis $\{Rp_\lambda\}$.
For $1\leq j\leq\ell$, use $\widehat\Gamma=R^{-1}\boldsymbol\Gamma R$
and substitute $-t^{1-j}u$ for $u$ in \eqref{c1:eq:G-definition} to obtain
\[
 \begin{aligned}
 \boldsymbol\Gamma(t^{1-j}u,v_j)Rp_\lambda
   &=R\widehat\Gamma(t^{1-j}u,v_j)p_\lambda\\
   &=\sum_{r\geq0}(-t^{1-j}u)^r\,R G_r(v_j)p_\lambda\\
   &=\sum_{r\geq0}(-u)^r t^{-(j-1)r}\,R G_r(v_j)p_\lambda .
 \end{aligned}
\]
The Macdonald inner product gives
$\langle p_\mu,Rp_\lambda\rangle_{q,t}=\delta_{\mu\lambda}z_\mu$.
Sum over sequences of partitions $\rho=(\rho(0),\ldots,\rho(\ell))$ with
$\rho(\ell)=\varnothing$, $\rho(0)=\mu$, and
$n_j=|\rho(j-1)|-|\rho(j)|\geq0$ for $1\leq j\leq\ell$.
Applying the operators in \eqref{c1:eq:character-definition} from right
to left, we take the coefficient of $Rp_{\rho(j-1)}$ in
$\boldsymbol\Gamma(t^{1-j}u,v_j)Rp_{\rho(j)}$ for $j=\ell,\ldots,1$.
Then $\sum_j n_j=|\mu|$ and $|\rho(j)|=\sum_{i>j}n_i$, so
\[
 2(\ell-1)|\mu|-\sum_j(j-1)n_j-\sum_j|\rho(j)|
 =2\sum_j(\ell-j)n_j\geq0.
\]
Taking the coefficient of $u^{|\mu|}$ and grouping the powers of $t$ gives
\[
 \begin{aligned}
 &(-1)^{|\mu|}t^{2(\ell-1)|\mu|}
   \widetilde{\boldsymbol\theta}_\mu^{(q,t)}(v_1,\ldots,v_\ell)\\
 &\quad=z_\mu\sum_{\rho}
   t^{2\sum_j(\ell-j)n_j}
   \prod_{j=1}^{\ell}\left(
      t^{|\rho(j)|}\Coef_{p_{\rho(j-1)}}
      \bigl[G_{n_j}(v_j)p_{\rho(j)}\bigr]\right).
 \end{aligned}
\]
For each $1\leq j\leq\ell$, Theorem~\ref{c1:thm:conjecture-one}
with input partition $\rho(j)$ and output partition $\rho(j-1)$ gives
\[
 t^{|\rho(j)|}z_{\rho(j-1)}\Coef_{p_{\rho(j-1)}}
 \bigl[G_{n_j}(v_j)p_{\rho(j)}\bigr]\in\mathbb N[q-1,t-1,-v_j].
\]
Since $z_{\rho(j-1)}>0$ and $t=1+(t-1)$, each summand has
nonnegative rational coefficients in the stated variables.
\end{proof}

\appendix
\section{Background and proof of Lemma~\ref{lem:linear-presentation-multiplication}}
\label{sec:appendix-lemma}

We recall the background needed for the vector bundles in
Section~\ref{subsec:cyclic-constant-rank} and for the proof of
Lemma~\ref{lem:linear-presentation-multiplication}.

\subsection{Background on vector bundles and cohomology}
\label{subsec:appendix-background}

Throughout this subsection, fix a field $K$ and an integer $n\geq1$,
and put $X=\mathbb P^n_K$.
Its homogeneous coordinates are $[c_0:\cdots:c_n]$; a $K$-rational
point is a nonzero vector in $K^{n+1}$ up to nonzero scalar multiplication.
The standard affine charts are $U_i=\{c_i\ne0\}$, for $0\leq i\leq n$.
On such a chart, regular functions are locally quotients of polynomials
in the coordinates $c_j/c_i$ with denominator nonvanishing on the
neighborhood in question. The structure sheaf $\mathcal O=\mathcal O_X$
assigns to each open set its ring of regular functions.

A sheaf $\mathcal F$ of $\mathcal O$-modules assigns to an open set $W$
an $\mathcal O(W)$-module $\mathcal F(W)$ of sections, with restriction maps to smaller
open sets. Restrictions compose as expected, and sections on an open
cover that agree on overlaps glue to a unique section.
A vector bundle of rank $b$ is identified with its sheaf of sections,
which is locally isomorphic to $\mathcal O^b$; such a sheaf is called
finite locally free, and a local identification with $\mathcal O^b$
is called a local trivialization. Its fiber at a point is the
$b$-dimensional vector space obtained by evaluating local sections
there, over the residue field of that point: the field in which regular
functions take their values there (over $K$ at a $K$-rational point).
A section over an open set assigns a vector to each fiber, with
regular coordinates in local trivializations.
Thus a fiber records values at one point, whereas the sheaf records
sections and their restrictions. A line bundle has rank one.
For these definitions, see \cite[Chapter~II, \S\S1 and~5]{Hartshorne}.

The tautological line bundle $\mathcal O(-1)$ has fiber $K\vec c$
at the point represented by $\vec c\in K^{n+1}\setminus\{0\}$.
Write $\mathcal F^\vee$ for the dual of a vector bundle $\mathcal F$,
obtained by taking the dual vector spaces in local trivializations.
Let $\mathcal O(1)=\mathcal O(-1)^\vee$ and, for an integer $d$,
let $\mathcal O(d)$ be the $d$th tensor power of $\mathcal O(1)$,
using dual powers for $d<0$ and $\mathcal O(0)=\mathcal O$.
Write $\mathcal F(d)=\mathcal F\otimes\mathcal O(d)$ for its twist.
The notation $\mathcal O(d)^b$ means a direct sum of $b$ copies.
Write $H^0(X,\mathcal F)=\Gamma(X,\mathcal F)=\mathcal F(X)$ for
global sections. For $\mathcal O(d)$, we have
\begin{equation}\label{app:eq:line-bundle-sections}
 H^0(X,\mathcal O(d))=\Gamma(X,\mathcal O(d))
       \cong K[c_0,\ldots,c_n]_d\qquad(d\geq0),
\end{equation}
where the subscript denotes the homogeneous component of degree $d$.
In particular, linear forms define maps $\mathcal O(-1)\to\mathcal O$,
and homogeneous polynomials of degree $d$ define maps
$\mathcal O\to\mathcal O(d)$.

Bundle maps are matrices of regular functions after local
trivializations. On $X$, a map of constant fiber rank can locally be
put in a matrix form with an identity block and a zero block:
choose an invertible minor and perform row and column operations.
Its kernel, image, and cokernel are therefore vector bundles.
An exact sequence means that each image equals the next kernel as
sheaves; this can be checked locally. A short exact sequence of vector
bundles splits locally, so ranks add and determinants multiply.
For a rank-$b$ bundle $\mathcal F$, its determinant is the line bundle
$\det\mathcal F=\bigwedge^b\mathcal F$.
Local splittings need not give a global splitting, and a surjective
map of bundles need not be surjective on global sections.

Sheaf cohomology extends global sections to groups $H^q(X,\mathcal F)$,
for integers $q\geq0$. In our setting it has a concrete description.
For vector bundles and the kernels and quotients of their sheaves of
sections, use the ordered \v Cech complex of the cover
$\mathfrak U=(U_0,\ldots,U_n)$. Write $\check C^q(\mathfrak U,\mathcal F)$
for its term in degree $q\geq0$, and $\delta^q$ for its differential,
the alternating sum of restrictions obtained by omitting one index.
The terms and their cohomology are given by
\begin{equation}\label{app:eq:cech-complex}
 \begin{aligned}
 \check C^q(\mathfrak U,\mathcal F)
   &=\bigoplus_{0\leq i_0<\cdots<i_q\leq n}
       \Gamma(U_{i_0}\cap\cdots\cap U_{i_q},\mathcal F),\\
 H^q(X,\mathcal F)&\cong\ker\delta^q/\operatorname{im}\delta^{q-1}.
 \end{aligned}
\end{equation}
Here the preceding term is zero when $q=0$, so
$\operatorname{im}\delta^{-1}=0$. For example, $\delta^0$ sends a family of
sections $(s_i)$ to their differences $s_j-s_i$ on $U_i\cap U_j$.
Thus its kernel consists of local sections agreeing on overlaps,
hence of global sections. For every $q$, the second line of
\eqref{app:eq:cech-complex} is the usual kernel modulo image.
The sheaves considered here are quasi-coherent, meaning that locally
they arise from modules on affine charts, with restriction given by
localization. Since the intersections of the $U_i$ are affine,
\eqref{app:eq:cech-complex} computes sheaf cohomology
\cite[Tag~01XD]{Stacks}.

For example, a short exact sequence
$0\to\mathcal F'\to\mathcal F\to\mathcal F''\to0$
gives an exact sequence beginning
\begin{equation}\label{app:eq:sections-long-exact}
 0\longrightarrow H^0(X,\mathcal F')
 \longrightarrow H^0(X,\mathcal F)
 \longrightarrow H^0(X,\mathcal F'')
 \longrightarrow H^1(X,\mathcal F').
\end{equation}
A section of $\mathcal F''$ lifts to one of $\mathcal F$ precisely
when its image in $H^1(X,\mathcal F')$ is zero
\cite[Chapter~III, Theorem~1.1A(c) and Proposition~2.6]{Hartshorne}.
This explains how sheaf cohomology detects the failure of surjectivity
on global sections.
We will use the following consequences of the line-bundle cohomology
calculation in \cite[Chapter~III, Theorem~5.1]{Hartshorne}:
\begin{equation}\label{eq:projective-line-bundle-vanishing}
 \begin{aligned}
 H^q(X,\mathcal O(-j))&=0
       &&(q\geq0,\ 1\leq j\leq n),\\
 H^q(X,\mathcal O(d))&=0
       &&(q>0,\ d\geq0).
 \end{aligned}
\end{equation}

A bounded complex $A$ of sheaves has finitely many nonzero terms $A^i$
and differentials $d_A^i:A^i\to A^{i+1}$ with $d_A^{i+1}d_A^i=0$.
Its cohomology sheaves are
\begin{equation}\label{app:eq:complex-cohomology}
 \mathcal H^i(A)=\ker d_A^i/\operatorname{im}d_A^{i-1}.
\end{equation}
For a complex of vector spaces, $H^i$ denotes the same kernel modulo
image construction. For an integer $j$, the shift $A[j]$ moves degree
$i+j$ to degree $i$. If the terms of $A$ are vector bundles, its dual
complex is obtained by dualizing the terms and reversing degrees.
Our differential conventions are
\begin{equation}\label{app:eq:complex-shift-dual}
 \begin{aligned}
 (A[j])^i&=A^{i+j},& d_{A[j]}^i&=(-1)^j d_A^{i+j},\\
 (A^\vee)^i&=(A^{-i})^\vee,
 &d_{A^\vee}^i&=(-1)^{i+1}(d_A^{-i-1})^\vee.
 \end{aligned}
\end{equation}
Tensor products with a line bundle and twists of a complex are taken
term by term. A chain map is a family of maps commuting with the
differentials. It is a quasi-isomorphism if it induces isomorphisms
on all cohomology sheaves. The derived category formally makes these
maps invertible; for complexes, $A\simeq B$ means an isomorphism in
this category. See \cite[Tag~05RR]{Stacks}.

For a chain map $f:A\to B$, the mapping cone has terms
$\operatorname{Cone}(f)^i=B^i\oplus A^{i+1}$ and differential
$(b,a)\mapsto(d_Bb+f(a),-d_Aa)$ on local sections.
It gives a distinguished triangle
\begin{equation}\label{app:eq:cone-triangle}
 A\xrightarrow{f}B\longrightarrow\operatorname{Cone}(f)
       \longrightarrow A[1].
\end{equation}
A distinguished triangle is a triangle isomorphic in the derived
category to one of this form. It gives a long exact sequence of
cohomology sheaves, continuing
$\mathcal H^i(A)\to\mathcal H^i(B)\to
\mathcal H^i(\operatorname{Cone}(f))\to\mathcal H^{i+1}(A)$;
see \cite[Tags~014E and~05RR]{Stacks}.
For bounded complexes of vector bundles, dualizing the cone
construction shows that dualization reverses distinguished triangles;
tensoring by a line bundle preserves them.

The standard truncation $\tau_{\leq-1}A$ keeps the terms below degree
$-1$, replaces $A^{-1}$ by $\ker d_A^{-1}$, and sets all terms in
nonnegative degrees to zero. It retains precisely the negative
cohomology sheaves. If $A$ is concentrated in nonpositive degrees,
there is a distinguished triangle
\begin{equation}\label{app:eq:truncation-triangle}
 \tau_{\leq-1}A\longrightarrow A\longrightarrow\mathcal H^0(A)
       \longrightarrow(\tau_{\leq-1}A)[1],
\end{equation}
where $\mathcal H^0(A)$ is placed in degree $0$;
see \cite[Tag~08J5, Remark~13.12.4]{Stacks}.
Conversely, a distinguished triangle $B\to A\to\mathcal F\to B[1]$
with $B$ having only negative cohomology and $\mathcal F$ a sheaf
in degree $0$ identifies $B\simeq\tau_{\leq-1}A$.
Truncation is functorial and commutes with tensoring by a line bundle.
The same construction in other degrees shows that, if $A$ has just
one nonzero cohomology sheaf $\mathcal F$ in degree $i$, then
$A\simeq\mathcal F[-i]$.

For a bounded complex $A$ of the sheaves considered above, its
hypercohomology $\mathbb H^k(X,A)$ combines the differential of $A$
with the local gluing calculation. More precisely, form the
\v Cech complexes of its terms and their total complex, whose term
in degree $k$ is
$\bigoplus_{p+q=k}\check C^q(\mathfrak U,A^p)$, with differential
$d_A+(-1)^p\delta$ on $\check C^q(\mathfrak U,A^p)$.
The cohomology of this total complex is $\mathbb H^k(X,A)$;
see \cite[Tags~01XD and~01FP]{Stacks}.
For a single sheaf in degree $0$, this recovers ordinary sheaf
cohomology. Thus $\mathcal H^i(A)$ is a sheaf computed from the
differentials, whereas $H^q(X,\mathcal F)$ and $\mathbb H^k(X,A)$
are vector spaces computed using local sections and restrictions.
For bounded complexes $A,B$, quasi-isomorphism invariance gives
\begin{equation}\label{app:eq:hypercohomology-invariance}
 A\simeq B\quad\Longrightarrow\quad
 \mathbb H^k(X,A)\cong\mathbb H^k(X,B)\qquad(k\in\mathbb Z).
\end{equation}

Only two calculation rules are needed below.
Write $\Gamma(X,A)$ for the vector-space complex obtained by taking
global sections of each term, with the induced differentials.
If $H^q(X,A^p)=0$ for every $p$ and every $q>0$, then
\begin{equation}\label{app:eq:hypercohomology-global-sections}
 \mathbb H^k(X,A)\cong H^k(\Gamma(X,A))
       \qquad(k\in\mathbb Z).
\end{equation}
This standard calculation is given in \cite[Tag~05TA]{Stacks}.
In the \v Cech description, every column then has cohomology only
in degree zero, where it equals $\Gamma(X,A^p)$.
In particular, if $H^q(X,A^p)=0$ for every $p$ and every $q\geq0$,
the remaining vector-space complex is zero. Hence
\begin{equation}\label{app:eq:hypercohomology-zero-terms}
 \mathbb H^k(X,A)=0\qquad(k\in\mathbb Z).
\end{equation}

\subsection{Proof of Lemma~\ref{lem:linear-presentation-multiplication}}
\label{subsec:appendix-lemma-proof}

We retain the setting of \eqref{eq:general-linear-presentation}
and \eqref{eq:section-multiplication-map}, and assume that
\eqref{eq:general-linear-presentation} is exact.
In particular, $n\geq2$, $a\geq1$, and $K$ has characteristic zero.
We use the conventions and calculations of
Appendix~\ref{subsec:appendix-background}.

\begin{proof}
Choose a basis of $V$. The assertion is that the products of $n$ of
the corresponding degree-$a$ forms, indexed by multisets of basis
elements, form a basis in degree $na$. We follow the hypercohomology
argument of \cite[Section~4.1]{RSWY}.

\smallskip
\noindent\textit{Step 1. Express the kernel and cokernel of $\alpha_n$ as hypercohomology.}
Put $\mathcal U=U\otimes\mathcal O(-1)$,
$\mathcal V=V\otimes\mathcal O$, and
$C=[\mathcal U\xrightarrow{\varphi}\mathcal V]$ in degrees $-1,0$.
By \eqref{eq:general-linear-presentation},
$\mathcal H^{-1}(C)=E$ and $\mathcal H^0(C)=L$.
For $s\geq0$, the symmetric-power complex $S_s=\operatorname{Sym}^sC$
has terms
\begin{equation}\label{eq:symmetric-power-complex}
 S_s^{-i}=\bigwedge^i\mathcal U\otimes
                    \operatorname{Sym}^{s-i}\mathcal V
       \qquad(0\leq i\leq s).
\end{equation}
For $1\leq i\leq s$, its differential is given below.
Here $u_1,\ldots,u_i$ are local sections of $\mathcal U$,
$f$ is a local section of $\operatorname{Sym}^{s-i}\mathcal V$,
and a hat indicates omission of a factor:
\[
 d_S(u_1\wedge\cdots\wedge u_i\otimes f)
 =\sum_{h=1}^i(-1)^{h-1}
   u_1\wedge\cdots\widehat{u_h}\cdots\wedge u_i
      \otimes\varphi(u_h)f.
\]
Locally, \eqref{eq:general-linear-presentation} splits $C$ into $E$
in degree $-1$, $L$ in degree $0$, and a two-term complex with
identity differential. The positive symmetric powers of this last
summand are exact (they are the positive-degree pieces of the Koszul
complex on independent variables). Only the exterior factors from
$E$ and the symmetric factors from $L$ contribute to cohomology, giving
\begin{equation}\label{eq:power-complex-cohomology}
 \mathcal H^{-i}(S_s)=\bigwedge^iE\otimes L^{s-i}
       \qquad(0\leq i\leq s).
\end{equation}
This identification is independent of the local splitting. To recover
the cohomology class of a product on the right, take the wedge $e$
of its $E$-factors and lift its $L$-factors to $\mathcal V$. Changing
one lift by $\varphi(u)$ changes the representative by
$e\otimes\varphi(u)f=d_S(u\wedge e\otimes f)$, where $f$ is the
product of the other lifts. The resulting classes agree on overlaps.

Append $L^s$ to $S_s$ in degree $1$ using the quotient
$\operatorname{Sym}^s\mathcal V\to L^s$ induced by $\varepsilon$.
The resulting complex $\widehat S_s$ retains precisely the negative
cohomology sheaves in \eqref{eq:power-complex-cohomology}.
For $s=n$, its negative-degree terms are
\[
 \bigwedge^iU\otimes\operatorname{Sym}^{n-i}V
                         \otimes\mathcal O(-i),
             \qquad 1\leq i\leq n,
\]
whose cohomology vanishes in every degree by
\eqref{eq:projective-line-bundle-vanishing}. The remaining terms,
$\operatorname{Sym}^nV\otimes\mathcal O$ and $L^n=\mathcal O(na)$,
have no higher cohomology by the same formula. We can therefore
compute hypercohomology by taking global sections, as in
\eqref{app:eq:hypercohomology-global-sections}. In degrees $0$ and $1$
this leaves
\[
 \operatorname{Sym}^nV\xrightarrow{\alpha_n}
                  H^0(\mathbb P^n_K,L^n).
\]
The quotient sends $v_1\cdots v_n$ to
the product of the sections $\varepsilon(v_j\otimes1)$, so its map on
global sections is $\alpha_n$ in \eqref{eq:section-multiplication-map}. Hence
\begin{equation}\label{eq:section-multiplication-kernel-cokernel}
 \begin{aligned}
 \mathbb H^0(\mathbb P^n_K,\widehat S_n)&\cong\ker\alpha_n,\\
 \mathbb H^1(\mathbb P^n_K,\widehat S_n)&\cong\operatorname{coker}\alpha_n.
 \end{aligned}
\end{equation}
It now suffices to find a complex representing $\widehat S_n$ whose
terms have vanishing cohomology in every degree. For $n=2$ we can
already see such a complex. Taking ranks and determinants in
\eqref{eq:general-linear-presentation} gives
$\operatorname{rank}E=1$ and $E\simeq\mathcal O(-a-1)$.
After augmentation, \eqref{eq:power-complex-cohomology} leaves the
single cohomology sheaf $E\otimes L\simeq\mathcal O(-1)$ in degree
$-1$. Thus $\widehat S_2\simeq\mathcal O(-1)[1]$, and the required
vanishing follows from \eqref{eq:projective-line-bundle-vanishing}.

\smallskip
\noindent\textit{Step 2. Compare successive symmetric powers.}
For general $n$, we first lower the symmetric power from $n$ to $n-1$.
Comparing once more will isolate the top exterior power of $E$.
For these comparisons we adapt the map in the proof of
\cite[Proposition~3.6]{RSWY}.
For an integer $k\geq1$ and local sections $v_1,\ldots,v_k$ of
$\mathcal V$, define
\[
 \delta_\varepsilon(v_1\cdots v_k)
  =\sum_{j=1}^k
       v_1\cdots\widehat{v_j}\cdots v_k\otimes\varepsilon(v_j).
\]
Set $\delta_\varepsilon(1)=0$. The map $\delta_\varepsilon$ is
$\mathcal O$-linear with values in
$\operatorname{Sym}^{k-1}\mathcal V\otimes L$.
For local sections $u$ of $\mathcal U$ and $f$ of
$\operatorname{Sym}^{k-1}\mathcal V$, the identity
$\varepsilon\varphi=0$ gives
\begin{equation}\label{eq:refaudit-delta-commutation}
 \delta_\varepsilon(\varphi(u)f)
 =\varphi(u)\delta_\varepsilon(f).
\end{equation}
For $s\geq1$, use the same normalization $1/s$ in every degree and set
\begin{equation}\label{eq:symmetric-comparison-map}
 D_s^{-i}=\frac1s\operatorname{id}_{\wedge^i\mathcal U}
                       \otimes\delta_\varepsilon
       \qquad(0\leq i\leq s).
\end{equation}
Equation~\eqref{eq:refaudit-delta-commutation} shows that these maps
commute with the differentials and define a chain map
$D_s:S_s\to S_{s-1}\otimes L$.
It also verifies the chain condition at $i=s$, where the target
term is zero. For $0\leq i<s$, the source and target cohomology
sheaves are both $\bigwedge^iE\otimes L^{s-i}$ by
\eqref{eq:power-complex-cohomology}. On a product representing such
a class, $\delta_\varepsilon$ selects one of the $s-i$ factors from
$L$. Each choice gives the same product in $L^{s-i}$, so
\begin{equation}\label{eq:symmetric-comparison-cohomology}
 \mathcal H^{-i}(D_s)
 =\frac{s-i}{s}\operatorname{id}_{\wedge^iE\otimes L^{s-i}}
       \qquad(0\leq i<s).
\end{equation}

Taking ranks in \eqref{eq:general-linear-presentation} gives
$\operatorname{rank}E=n-1$. Writing
$\det E=\bigwedge^{n-1}E$ for its determinant line bundle, we also obtain
\begin{equation}\label{eq:linear-presentation-determinant}
 \det E\simeq\mathcal O(-(n-1)(a+1)),\qquad
 \det E\otimes L^{n-1}\simeq\mathcal O(-(n-1)).
\end{equation}
Indeed, $\det(\ker\varepsilon)\simeq L^{-1}$ and
$\det E\simeq\det\mathcal U\otimes L$.
We use \eqref{eq:symmetric-comparison-cohomology} first with $s=n$.
Since $\bigwedge^iE=0$ for $i\geq n$, every nonzero cohomology sheaf
has $i\leq n-1$, and the scalar $(n-i)/n$ is nonzero in $K$.
Consequently $D_n$ is a quasi-isomorphism, giving
\begin{equation}\label{eq:refaudit-symmetric-quasi-isomorphism}
 S_n\simeq S_{n-1}\otimes L.
\end{equation}
Recall that augmentation removes $\mathcal H^0(S_s)=L^s$.
In the notation of \eqref{app:eq:truncation-triangle}, this reads
\begin{equation}\label{eq:augmentation-truncation}
 \widehat S_s\simeq
 \operatorname{Cone}(S_s\to L^s)[-1]
 \simeq\tau_{\leq-1}S_s.
\end{equation}
To extend $D_n$ to the augmented complexes, take the identity on
$L^n$ in degree $1$. Indeed, applying $D_n^0$ to $v_1\cdots v_n$
and then multiplying in $L^n$ gives $n$ copies of
$\varepsilon(v_1)\cdots\varepsilon(v_n)$, each with coefficient $1/n$.
This is exactly the original augmentation. The extended map induces
the same isomorphisms in negative cohomological degrees as $D_n$,
and both augmented complexes have zero cohomology sheaves elsewhere. Hence
\begin{equation}\label{eq:augmented-symmetric-comparison}
 \widehat S_n\simeq\widehat S_{n-1}\otimes L.
\end{equation}

Now put $m=n-1$ and use the same comparison with $s=m$.
By \eqref{eq:symmetric-comparison-cohomology}, $D_m$ is an isomorphism
on every cohomology sheaf in degrees $-m+1,\ldots,0$.
The only remaining source cohomology is
$\mathcal H^{-m}(S_m)=\det E$, while the target has no term in degree $-m$.
The long exact sequence of cohomology sheaves for the mapping cone
therefore shows that $\operatorname{Cone}(D_m)[-1]$ has just one
nonzero cohomology sheaf, namely $\det E$ in degree $-m$.
It is thus isomorphic to $\det E[m]$, and the mapping-cone triangle becomes
\begin{equation}\label{eq:symmetric-top-triangle}
 \det E[m]\longrightarrow S_m
 \xrightarrow{D_m}S_{m-1}\otimes L
 \longrightarrow\det E[m+1].
\end{equation}

\smallskip
\noindent\textit{Step 3. Use self-duality to prove the vanishing.}
We now seek a model for the negative-degree part of $S_m$, which
the augmentation retains. Duality reverses the arrows in
\eqref{eq:symmetric-top-triangle}. We choose the twist and shift
that turn the leftmost term $\det E[m]$ into $L^m$ on the right
and identify the middle term with $S_m$.
Indeed, the differential of $C$ has constant rank, its kernel $E$
has rank $m$, and its cokernel $L$ is a line bundle. Since divided
powers agree with symmetric powers in characteristic zero,
\cite[Theorem~3.2]{RSWY} gives
\begin{equation}\label{eq:derived-self-duality}
 S_m\simeq S_m^\vee\otimes\det E\otimes L^m[m].
\end{equation}
Dualize \eqref{eq:symmetric-top-triangle} and tensor by
$\det E\otimes L^m[m]$. Using \eqref{eq:derived-self-duality} for the
middle term, we obtain a distinguished triangle beginning with
\begin{equation}\label{eq:symmetric-dual-triangle}
 S_{m-1}^\vee\otimes\det E\otimes L^{m-1}[m]
 \longrightarrow S_m\longrightarrow L^m.
\end{equation}
Here the rightmost term is $L^m$ because
$(\det E[m])^\vee\otimes\det E\otimes L^m[m]=L^m$.
By \eqref{eq:symmetric-power-complex}, dualizing $S_{m-1}$ moves
its degrees to $0,\ldots,m-1$, and shifting by $[m]$ moves them
to $-m,\ldots,-1$. Thus the first term of
\eqref{eq:symmetric-dual-triangle} is concentrated in negative degrees,
while $L^m$ is concentrated in degree $0$.
The long exact sequence now shows that the left arrow accounts for
all the negative cohomology sheaves of $S_m$, and the right arrow
induces an isomorphism on $\mathcal H^0$. By the characterization
following \eqref{app:eq:truncation-triangle}, the left term represents
$\tau_{\leq-1}S_m$, which is $\widehat S_m$ by
\eqref{eq:augmentation-truncation}. Thus
\begin{equation}\label{eq:augmented-dual-comparison}
 \widehat S_m\simeq
 S_{m-1}^\vee\otimes\det E\otimes L^{m-1}[m].
\end{equation}
Combining \eqref{eq:augmented-symmetric-comparison},
\eqref{eq:augmented-dual-comparison}, and
\eqref{eq:linear-presentation-determinant} yields
\begin{equation}\label{eq:augmented-acyclic-model}
 \begin{aligned}
 \widehat S_n
 &\simeq S_{n-2}^\vee\otimes\det E\otimes L^{n-1}[n-1]\\
 &\simeq S_{n-2}^\vee(-(n-1))[n-1].
 \end{aligned}
\end{equation}

We can now prove the required vanishing.
By \eqref{eq:symmetric-power-complex}, the term of
$S_{n-2}$ in degree $-i$ is a direct sum of $\mathcal O(-i)$,
for $0\leq i\leq n-2$.
After dualizing and twisting by $\mathcal O(-(n-1))$, these become
direct sums of $\mathcal O(i-(n-1))$.
The twists range from $-(n-1)$ to $-1$, so all their sheaf cohomology
vanishes by \eqref{eq:projective-line-bundle-vanishing}.
Apply \eqref{app:eq:hypercohomology-zero-terms} to this complex.
Then \eqref{eq:augmented-acyclic-model} and the invariance in
\eqref{app:eq:hypercohomology-invariance} give
\begin{equation}\label{eq:augmented-hypercohomology-vanishing}
 \mathbb H^q(\mathbb P^n_K,\widehat S_n)=0
       \qquad(q\in\mathbb Z).
\end{equation}

Combining \eqref{eq:augmented-hypercohomology-vanishing} and
\eqref{eq:section-multiplication-kernel-cokernel} shows that
$\alpha_n$ is injective and surjective, proving the lemma.
\end{proof}

\section*{Acknowledgments}
N.J. is partially supported by the Simons Foundation
(No. MP-TSM-00002518) and National Natural Science Foundation of China
(No. 12171303).  N.L. is supported by Beijing Natural Science Foundation (No. 1264053) and the Fundamental Research Funds for the Central Universities (No. JZ2026HGTS0237).

\end{document}